\documentclass{amsart} \usepackage{verbatim} %\input{sch1.sty}
\usepackage{amssymb} \usepackage{amscd} \usepackage{epic}
\usepackage{eepic} \usepackage{epsfig} \usepackage{pb-diagram}
\usepackage{bm}

\renewcommand{\Im}{\operatorname{Im}} \renewcommand{\Re}{\operatorname{Re}}

\newcommand{\tphi}{{\tilde{\phi}}}

\newcommand{\vecx}{\vec{x}} \newcommand{\vecy}{\vec{y}}
\newcommand{\vecnu}{\vec{\nu}} \newcommand{\vecmu}{\vec{\mu}}
\newcommand{\vecv}{{\vec{v}}}

\newcommand\half{\frac1{2}}

\renewcommand\MR[1]{} \newcommand\evo{\lambda_0} %changed this back!; JW
\newcommand\evosq{\lambda_0^2}
\newcommand\evosr{\lambda_0}
\newcommand\sigmasch{\sigma_{\SCH}} \newcommand\XXbh{X}
\newcommand\XXb{X^2_{\text{b}}} 
\newcommand\CIdot{\dot C^\infty} 
\newcommand\Wt{\tilde W} \newcommand\Vr{\mathcal{V}_r}
\newcommand\Vc{\mathcal{V}_c} \newcommand\gtcl{\overline{\gamma^2}}
\newcommand\Nub{\mathcal{V}_b} \renewcommand\ss{{\mathbf{\sigma}}}
\newcommand\xx{{\mathbf{x}}} 
\newcommand\xxb{{\mathbf{{x_b}}}} 
\newcommand\qx{{\mathbf{x_q}}}
\newcommand\id{\operatorname{id}}
\newcommand\phit{\tilde \phi} \newcommand\Ghat{\hat G}
\newcommand\Gsharp{G^\sharp} 
\newcommand\Nullset{\mathcal{N}}

\theoremstyle{plain} \newtheorem{theorem}{Theorem}[section]
\newtheorem{proposition}[theorem]{Proposition}
\newtheorem{lemma}[theorem]{Lemma}
\newtheorem{corollary}[theorem]{Corollary}

\theoremstyle{definition} \newtheorem{defn}[theorem]{Definition}
\newtheorem{example}[theorem]{Example}

\theoremstyle{remark} \newtheorem*{remark}{Remark}

\newlength{\KDlen} \newlength{\sclen} \newlength{\tlen} \newlength{\slen}
\newlength{\qsclen}  \newlength{\omegalen}

\numberwithin{equation}{section}

\DeclareMathSymbol{\leqslant}{\mathrel}{AMSa}{"36}
\DeclareMathSymbol{\geqslant}{\mathrel}{AMSa}{"3E}
\DeclareMathSymbol{\gtrless}{\mathrel}{AMSa}{"3F}
\DeclareMathSymbol{\lessgtr}{\mathrel}{AMSa}{"37}

\renewcommand{\leq}{\leqslant} \renewcommand{\geq}{\geqslant}

\def\sch/{{Schr\"odinger}} % to use macro, type \sch/ --- 
\def\psido/{{$\Psi$DO}}

\newcommand{\SC}{\ensuremath{\mathrm{sc}}}
\newcommand{\SF}{\ensuremath{\mathrm{s}\Phi}}
\newcommand{\QSF}{\ensuremath{\mathrm{qs}\Phi}}
\newcommand{\QSC}{\ensuremath{\mathrm{qsc}}}
\newcommand{\B}{\ensuremath{\mathrm{b}}}
\newcommand{\SCH}{\ensuremath{\mathrm{sc,h}}}

\newcommand\restrictedto{\!\!\upharpoonright}

 \newcommand{\RR}{{\mathbb{R}}}

\newcommand\CI{{\mathcal{C}}^{\infty}}

 \newcommand{\h}{{\frac{1}{2}}}
\newcommand{\ep}{{\epsilon}} 
\DeclareMathOperator{\sgn}{{sgn}}

 \newcommand{\Tsf}{{^{\SF}T}}
\newcommand{\Tsfstar}[1][\mbox{}]{{^{\SF}T^*_{#1}}}
\newcommand{\Tscstar}[1][\mbox{}]{{^{\SC}T^*_{#1}}}
\newcommand{\Tscstarbar}{\Tbarscstar}
\newcommand{\Tqsfstar}[1][\mbox{}]{{^{\QSF}T^*_{#1}}}
\newcommand{\Nsfstar}{{^{\SF}N^*}}
\newcommand{\Nqsfstar}{{^{\QSF}N^*}}
\newcommand{\Nqsf}{{^{\QSF}N^*}}

\newcommand\near{{\operatorname{{near}}}}
\newcommand\far{{\operatorname{{far}}}}
\newcommand\innt{{\operatorname{{int}}}}

\newcommand{\Tsc}[1][\mbox{}]{{\settowidth{\tlen}{$T$}
    \settowidth{\sclen}{{\scriptsize \SC}} \hspace{\sclen}
    T_{#1}^{\hspace{-\sclen}\hspace{-\tlen} \SC\hspace{\tlen}}}}

\newcommand{\Tbarscstar}[1][\mbox{}]{{\settowidth{\tlen}{$\overline{T}$}
    \settowidth{\sclen}{{\scriptsize \SC}}
 \hspace{1.1\sclen}\overline{T}_{#1}^{\hspace{-1.1\sclen}\hspace{-\tlen}
      \SC
 \hspace{\tlen}\hspace{.1\sclen} \ast}}}

\newcommand{\nqsc}[1][\cdot]{N_{\QSC}}

 \DeclareMathOperator{\supp}{{supp}}

     \newcommand{\pa}{{\partial}}

    \newcommand{\abs}[1]{{\left\lvert{#1}\right\rvert}}
    \newcommand{\norm}[1]{{\left\lVert{#1}\right\rVert}}
    
    \newcommand{\ang}[1]{{\left\langle{#1}\right\rangle}}

\newcommand\Vsc{\mathcal{V}_{\SC}} 
\newcommand\Vqsc{\mathcal{V}_{\QSC}} \newcommand\Vb{\mathcal{V}_{\B}}
\newcommand\Vsf{\mathcal{V}_{\SF}} \newcommand\Vqsf{\mathcal{V}_{\QSF}}
\def\Nsf{{}^{\SF}N^*} \def\CIsf{C^{\infty}_\Phi}

\newcommand{\mar}[1]{{\marginpar{\sffamily{\scriptsize #1}}}}
\newcommand{\jw}[1]{{\mbox{}}} \newcommand{\ah}[1]{{\mbox{}}}

\newcommand{\Lap}{\Delta} \def\face{\operatorname{Fa}}
\def\ilabel#1{{\label{#1}}}

\def\dbyd#1#2{\frac{\partial #1}{\partial #2}}

\newcommand\Lt{\tilde L} \newcommand\mf{\operatorname{mf}}

\newcommand\scH{{}^{\operatorname{sc}}H} 
 \newcommand\Id{\operatorname{Id}}
\newcommand\sub{\operatorname{sub}} 

\newcommand\fibre{{\text{fibre}}}

\newcommand{\bface}{\text{bf}} \newcommand{\bfc}{\text{bf}}
\newcommand{\rb}{\text{rb}} \newcommand{\lb}{\text{lb}}
 
\renewcommand{\sf}{\text{sf}} \newcommand\finfty{\operatorname{fi}}
\newcommand\sinfty{\operatorname{scl}}

\newcommand\XXsc{X^2_{\SC}} 
\newcommand\MMb{M^2_{\B}}
\newcommand\MMsc{M^2_{\SC}}

\newcommand\Gh{\hat G} \newcommand\Gt{\tilde G} 
 \def\Lsharp{{L^\sharp}} 
\newcommand\Jsharp{J^{\sharp}} \newcommand\Jt{\tilde J}

\def\sfOh{{}^{\sf} \Omega^\half}
\def\scOh{{}^{\operatorname{sc}} \Omega^\half} 
\def\qsfOh{{}^{\QSF} \Omega^\half} 
\newcommand\Diag{\Delta} %%changed this-- JW
 
 \newcommand\Psisch{\Psi_{\SCH}}
\newcommand\rest{\restriction} \def\sfOh{{}^{\SF} \Omega^\half}
\newcommand\Diffsf{{}^{\SF}\operatorname{Diff}}
\newcommand\SR{\operatorname{SR}} \newcommand\mubar{\overline{\mu}}
\newcommand\spn{\operatorname{span}} \newcommand\cZ{\mathcal{Z}}
\newcommand\nubar{\overline{\nu}} \newcommand\muhat{{\hat \mu}}
\newcommand\Lc{\overline{L}}

\newcommand\spnn{\operatorname{span}}

\title[Semiclassical resolvent]{The semiclassical resolvent and the
propagator for nontrapping scattering metrics} \author{Andrew Hassell}
\address{Department of Mathematics, Australian National University,
Canberra, 0200 ACT Australia} \author{Jared Wunsch} \address{Department of
Mathematics, Northwestern University, Evanston IL USA} 
\begin{document}
\maketitle

%\centerline{\today}

\begin{abstract} Consider a compact manifold with boundary $M$ with
  a scattering metric $g$ or, equivalently, an asymptotically conic
  manifold $(M^\circ, g)$. (Euclidean $\RR^n$, with a compactly supported
  metric perturbation, is an example of such a space.)  Let $\Delta$ be the
  positive Laplacian on $(M,g)$, and $V$ a smooth potential on $M$ which
  decays to second order at infinity.  In this paper we construct the
  kernel of the operator $(h^2 \Delta + V - (\evo \pm i0)^2)^{-1}$, at a
  nontrapping energy $\evo > 0$, uniformly for $h \in (0, h_0)$, $h_0 > 0$
  small, within a class of Legendre distributions on manifolds with
  codimension three corners. Using this we construct the kernel of the
  propagator, $e^{-it(\Delta/2 + V)}$, $t \in (0, t_0)$ as a quadratic
  Legendre distribution.  We also determine the global semiclassical structure of the spectral projector, Poisson operator and scattering matrix.
   \end{abstract}

\tableofcontents

\part{Introduction}
\section{Overview}

In this paper we analyze the structure of the semiclassical resolvent  on a
class of noncompact manifolds with asymptototically conic ends.  The class
of asymptototically conic, or `scattering,' 
manifolds, introduced by Melrose \cite{Melrose:Spectral}, consists of those
Riemannian manifolds that can be described as the interior of a manifold
$M$ with boundary, such that in terms of some boundary defining function
$x,$ we can write the metric $g$ near $\pa M$ as
$$ g = \frac{dx^2}{x^4} + \frac{k}{x^2}
$$ where $k$ is a smooth 2-cotensor on $M$ with $k\restrictedto_{\pa M}$ a nondegenerate metric on $\pa M;$ there is no loss of generality in assuming that $k$ has no $dx$ component, so that $k = k(x, y, dy)$ \cite{Joshi-SaBarreto2}. 
In terms of $r = 1/x$ this reads
$$
g = dr^2 + r^2 k(\frac1{r}, y, dy)
$$
and is thus asymptotic to the exact conic metric $dr^2 + r^2 k(0, y, dy)$ as $r \to \infty$. The interior $M^\circ$ of $M$ is thus metrically complete, with  the boundary of $M$ `at infinity'. An
important class of examples is that of asymptotically Euclidean spaces,
pictured in a radial compactification: here $M$ is the unit ball, and
$k\restrictedto_{S^{n-1}}$ is the standard metric on the sphere.  More
generally, collar neighborhoods of boundary components are large conic ends
of the scattering manifold.  

We are concerned here with the operator $H = h^2 \Delta + V$, where $\Delta =
\Delta_g$ is the Laplacian on $M$ with respect to the metric $g$, $h \in
(0, h_0)$ is a small parameter (`Planck's constant') and $V \in x^2\CI(M)$
is a real potential function, smooth on $M$ and vanishing to second order
at the boundary (hence,  $V$ is $O(r^{-2})$ and thus short-range). The bulk
of this paper is concerned with the analysis of the semiclassical
resolvent, i.e.\ the operator $(h^2 \Delta + V - \evosq)^{-1}$, for $\evosr$
real, or more precisely the limit of this as $\evosr$ approaches the real
axis from above or below, denoted $(h^2 \Delta + V - (\evosr \pm
i0)^2)^{-1}$. 

For $\evosr$ non-real, the resolvent $(H - \evosq)^{-1}$ is a relatively
simple object, as $H - \evosq$ is then an elliptic operator in the
`semiclassical scattering calculus' of pseudodifferential operators, hence
a parametrix, and indeed the inverse itself, lies in this calculus
\cite{Melrose:Spectral}, \cite{Wunsch-Zworski}; also see
Section~\ref{sec:hpseudos}. In the limit as $\Im \evosr\to 0,$ ellipticity,
in the strengthened sense required by the scattering calculus, fails, and
the resolvent becomes more complicated.  Hassell and Vasy
\cite{hasvas1,hasvas2} analyzed the resolvent in this regime for a fixed $h
> 0$.  In this paper, we analyze the resolvent $(H - (\evosr \pm
i0)^2)^{-1}$ uniformly as $h \to 0$.  We assume throughout that the energy
level $\evosq$ is \emph{nontrapping}. That is, we assume that every
bicharacteristic of the operator $H - \evosq$ reaches the boundary
$\partial M$ in both directions, or equivalently, every bicharacteristic
eventually leaves each compact set $K \subset M^\circ$. In the case $V
\equiv 0$, bicharacteristics are simply geodesics and the condition is that
there are no trapped geodesics: every maximally extended geodesic reaches
infinity in both directions.

Our main result is the identification of the Schwartz kernel $R_\pm$ of $(H - (\evosr  \pm i0)^2)^{-1}$ as a
\emph{Legendrian distribution}.
We now informally describe Legendrian distributions, and how
those arising in the Schwartz kernel of $R_\pm$ are associated to the
underlying geometry of the problem.  First, a Legendrian distribution on a
manifold $N$ with boundary is a smooth function on the interior of $N$ with
singular, oscillatory behavior at $\pa N.$ It can locally be written as a
sum of oscillatory integrals of the form $\int a(x,y,v) e^{i \phi(y,v)/x}
\, dv$ where $x$ is a boundary defining function, $y$ are variables in $\pa N,$ $\phi$ satisfies a nondegeneracy condition, and the variable $v$ ranges over a compact set in some Euclidean
space.  Associated to such a distribution (indeed, parametrized by $\phi$
much as in H\"ormander's theory of Lagrangian distributions) is a
Legendrian manifold in the \emph{scattering cotangent bundle}, a rescaled version of the cotangent bundle of $N,$
restricted to $\pa N;$ this bundle has a natural contact structure (see Definition~\ref{sc-def}). 
Legendre distributions, introduced in \cite{MZ}, were generalized to the
setting of manifolds with codimension-two corners, with fibered boundaries,
by Hassell-Vasy \cite{hasvas1}.  Here we further generalize to
codimension-three boundaries with fibrations; these arise naturally as the
Schwartz kernel of $R_\pm$ lies in the manifold with codimension-three corners $M \times M
\times [0,h_0)_h,$ while the fibrations on the various faces arise from
projection operators on this product.  Indeed, we need a  further
refinement: a class of `Legendrian conic pairs' generalizing that
constructed by Melrose-Zworski \cite{MZ} and  Hassell-Vasy
\cite{hasvas1}; these distributions are associated to pairs of Legendrian
manifolds, one of which is allowed to have a conic singularity at its
intersection with the other.

The manifold $M \times M \times [0, h_0)$ is too crude a space on which to
describe the structure of the resolvent kernel.  In particular, the
asymptotic behaviour of the kernel at the corner $\partial M \times
\partial M \times [0, h_0)$ will depend in a complicated way on the angle
of approach.  We work on the space $X$ which is obtained from $M \times M
\times [0, h_0)$ by blowing up\footnote{This operation amounts analytically
to the introduction of polar coordinates in the transverse coordinates, or
geometrically to the introduction of a new boundary hypersurface replacing
the corner; see \cite{M} or \S\ref{Lconicpts} for details.} $\partial M
\times \partial M \times [0, h_0)$. That is, $X = \MMb \times [0, h_0)$
where $\MMb$, the b-double space, introduced in \cite{tapsit}, is the blowup of $M^2$ at $(\partial
M)^2$. The space $X$ has four boundary hypersurfaces: the `main face' $\MMb
\times \{ 0 \}$, denoted $\mf$; the left and right boundaries $\lb$ and
$\rb$, which are $\partial M \times M \times [0, h_0)$ and $M \times
\partial M \times [0, h_0)$, respectively, lifted to $X$; and the boundary
hypersurface created by the blowup, which we denote $\bfc$, which is a
quarter-circle bundle over $(\partial M)^2 \times [0, h_0)$. (See
Figure~\ref{fig:X} in Section~\ref{fir}.)

The resolvent kernel is most naturally described
in\jw{Notation check: $\XXbh$ is blown up triple space?}\ah{We seem to have given up on $\XXbh$; I am just using $X$ instead}  two pieces,
$R_\pm = K_\psi+K'.$ The first of these, $K_\psi$, is a semiclassical
scattering pseudodifferential operator, and thus has the same microlocal
structure as the (true) resolvent kernel when $\Im \evosr \neq 0$. It
captures the diagonal singularity of $R_\pm$ in the interior of $X$, but
not uniformly as the boundary of $X$ is approached.

The other piece, $K'$, is Legendrian in nature; in particular, it is smooth
in the interior of $X$, but is oscillatory as the boundary is
approached. There are in fact \emph{three} Legendrian submanifolds
associated to $K'$. The first is the conormal bundle at the diagonal,
denoted $N^* \Delta_b$, reflecting the fact that the pseudodifferential
part $K_\psi$ cannot capture the singularities at the diagonal near its
intersection with the characteristic variety $\Sigma(H - \evosq)$; this
intersection is nontrivial at the boundary of the diagonal, both at $\mf$
and at $\bfc$.  The second is what we call the `propagating Legendrian'
$L$, which is obtained by flowout from the intersection of the
characteristic variety and $N^* \Delta_b$ under the Hamilton vector field
associated to $H - \evosq$.  In fact, $L$ is divided into two halves, $L =
L_+ \cup L_-$ by $N^* \Delta_b$, and the incoming ($-$)/outgoing ($+$)
resolvent $R_\pm$ is singular only at $L_\pm$. The geometry of $(N^*
\Delta_b, L_\pm)$ is that of a pair of cleanly intersecting Legendre
submanifolds, and $K'_\pm$ microlocally lies in a calculus of `intersecting
Legendre distribution' associated to this pair, analogous to the class of
intersecting Lagrangian distributions of Melrose and Uhlmann \cite{MU}.
The propagating Legendrian $L$ turns out to have conic singularities at
$\bfc$, and another Legendrian, $L_2^\sharp$, appears,
to `carry off' the singularities at the conic intersection.  This
latter Legendrian consists of those points in phase space over $\bfc$ which
point in pure outgoing/incoming directions in both variables.  We thus
state our first main theorem as follows (relevant classes of Legendrian
distributions are defined below in \S\ref{sec:legendrian}--\S\ref{qsfs},
while the particular Legendrian manifolds referred to are discussed in
\S\ref{sec:prop}):

\begin{theorem}\ilabel{main1} The semiclassical outgoing resolvent kernel
  $R_+$, multiplied by the density factor $|dh|^{1/2}$, is the sum of a
  semiclassical pseudodifferential operator of order $(-2,0,0)$, an
  intersecting Legendre distribution associated to the diagonal Legendrian
  $N^* \Delta_b$ and the propagating Legendrian $L_+,$ and a conic
  Legendrian pair associated to $L_+$ and the outgoing Legendrian
  $L_2^\sharp$. The orders of the Legendrians are $3/4$ at $N^* \Delta_b$,
  $1/4$ at $L_+$, $(2n-3)/4$ at $L_2^\sharp$, $-1/4$ at $\bfc$, and
  $(2n-1)/4$ at $\lb$ and $\rb$.
\end{theorem}

This is rather similar in nature to the main result of
\cite{hasvas1}.\ah{Omit this para?} The main difference is that in
\cite{hasvas1} only propagation inside $\bfc$ needed to be considered; this
is closely related to geodesic flow `at infinity', and only involves
exact conic geometry. Here, by contrast (and in the case $V \equiv 0$) the
geodesic flow over the entire manifold $M$ is relevant.

\mbox{}From Theorem~\ref{main1} we can obtain analogous results for other
fundamental operators in scattering theory, including the  spectral
projections, Poisson operator and scattering matrix, since their kernels
can be obtained from the resolvent in a straightforward way. The simplest
one is the spectral measure $dE(h^{-2})$ which is $1/2\pi i$ times the
difference between the incoming and outgoing resolvents. In taking this
difference the diagonal singularity disappears and we obtain, in the
notation of Section~\ref{Legdist-conic},

\begin{corollary}\ilabel{sp-cor} The\ah{Orders?} spectral measure $dE(h^{-2})$ times $|dh/h^2|^{-1/2}$ is an intersecting Legendre distribution associated to the conic pair $(L, L_2^\sharp)$:
$$
dE(h^{-2}) \otimes |dh/h^2|^{-1/2} \in I^{1/4, n/2 - 3/4; n/2 - 1/4, -1/4}((L, L_2^\sharp), X; \sfOh).
$$
\end{corollary}

This generalizes results (in the globally nontrapping setting) of Vainberg \cite{Vainberg2} and Alexandrova \cite{Alexandrova2}. 

The Poisson operator $P(h^{-1})$, as defined in \cite{MZ}, takes a function $f$ on the boundary $\partial M$ and maps it to that generalized eigenfunction with eigenvalue $h^{-2}$ with outgoing data $f$. It can be obtained from the resolvent kernel by restriction to $\rb$ after the removal of an oscillatory factor. The Legendrians $L$ and $L_2^\sharp$ themselves have `boundary values' at $\rb$, which are denoted $\SR$ and $G^\sharp$ respectively. Here $\SR$ stands for `sojourn relation'; it is the twisted graph of a contact transformation identified in \cite{HW}, and is related to the sojourn time of Guillemin \cite{Gu} (see Section~\ref{srt} for further discussion). 

\begin{corollary}\ilabel{Poisson-cor} The Poisson kernel $P(h^{-1})$, times $|dh/h^2|^{1/2}$, is  a Legendre distribution associated to the conic pair $(\SR, G^\sharp)$ of order $(0, (n-1)/2; 0)$. 
\end{corollary}

The scattering matrix $S(h^{-1})$ is, in turn, obtained by restricting the
kernel of $P(h^{-1})$ to the boundary, now at $\rb \cap \bfc$, although the
limit here is more subtle, compared to that for the Poisson operator, as it
only makes sense distributionally---this was explained in \cite{MZ}. The
sojourn relation $\SR$ has a `boundary value' at $\bfc$ which we denote $T$
and call the `total sojourn relation.'
We obtain a global characterization of the S-matrix as an oscillatory
function. It has
two kinds of behavior: for fixed $h>0$ it was shown by Melrose-Zworski to be
a Fourier integral operator on $\pa M,$ i.e.\ a Lagrangian distribution on
$\pa M \times \pa M.$ On the other hand, \emph{away from these
singularities,} it has been shown by Alexandrova to be a semiclassical
FIO \cite{Alexandrova}, \cite{Alexandrova1}.  Our structure theorem is
that the semiclassical scattering matrix globally lies in a calculus of
`Legendrian-Lagrangian' distributions (defined in Section~\ref{qsfs}) that combine these two different behaviors.

\begin{theorem}\ilabel{sc-thm} The scattering matrix $S(h^{-1})$, times $|dh/h^2|^{1/2}$, is a Legendrian-Lagrangian distribution of order $(-1/4, -1/4)$ associated to the total sojourn relation $T$. 
\end{theorem}

In a prior paper, the authors constructed a partial parametrix for the
Schr\"odinger propagator on nontrapping scattering manifolds; this
parametrix was valid in regions where one variable may range out to $\pa M$
(i.e.\ out to `infinity') but the other is restricted to lie in a compact
set in $M^\circ.$ Here, by integrating over the spectrum and using Corollary~\ref{sp-cor}, 
we are able to extend our description of the Schr\"odinger propagator to a
global one. To state the theorem, we note that, based e.g.\ on the form of
the free propagator $(2\pi it)^{-n/2} e^{i|z-z'|^2/2t}$ on $\RR^n$, we
expect the propagator to be Legendrian, with semiclassical parameter $t$,
but with  \emph{quadratic} oscillations at spatial infinity. We define such
a class of quadratic scattering-fibred Legendre
distributions. Corresponding to the Legendre submanifolds $L$, $L_2^\sharp$
introduced earlier are quadratic Legendre submanifolds $Q(L)$, $Q(L_2^\sharp)$ (see Section~\ref{qsfs}). Our result is

\begin{theorem}\ilabel{prop-thm} The Schr\"odinger propagator $e^{-it(\Delta/2+V)}$ is for $0 < t < t_0 < \infty$ a quadratic Legendre distribution associated to the conic pair $(\tilde L, \tilde G_2^\sharp)$:
\begin{equation}
e^{-it((1/2)\Delta + V)} \in I^{3/4, n/2+1/4; 1/4,-n/2+1/4}(\MMb \times [0,
  t_0), (Q(L), Q(L_2^\sharp)); \qsfOh).
\ilabel{propspace}\end{equation}
\end{theorem}

\
 
The resolvent construction described here is a direct generalization of
work of Hassell-Vasy \cite{hasvas1,hasvas2} on the fixed-energy
resolvent. This work was in turn motivated by the paper \cite{MZ} of
Melrose-Zworski on the Poisson operator and scattering matrix for
scattering metrics. All these works are ultimately based on the original
paper of Melrose \cite{Melrose:Spectral}.  The construction is also related
to the parametrix construction of Isozaki-Kitada \cite{IK}, which is valid
in the outgoing region.\jw{Check this thoroughly!  Get right citation,
too.}

Our results on the scattering matrix have many antecedents.  The
description of the behavior of the scattering matrix in the semiclassical
regime, away from singularities of the kernel (which occur at the diagonal
in $\RR^n$ with the usual normalizations) originates with Majda
\cite{Majda} for the case of obstacles and for compactly-supported metric
perturbations of $\RR^n$ by Guillemin \cite{Gu}.  The semiclassical limit
on $\RR^n$ with potential has been studied by Protas \cite{Protas},
Vainberg \cite{Vainberg1}, Yajima \cite{Yajima}, Robert-Tamura
\cite{Robert-Tamura}, and Alexandrova \cite{Alexandrova}, in varying
degrees of generality.  (See \cite{Alexandrova} for a clear summary of this
literature.)

Numerous authors have studied the structure of the Schr\"odinger kernel on
flat space (with a potential).  In this setting parametrices have been
constructed by Fujiwara \cite{Fujiwara1}, Zelditch \cite{Zelditch2},
Tr\`eves \cite{Treves1} and Yajima \cite{Yajima1}.  For a
compactly-supported nontrapping perturbation, Kapitanski-Safarov
\cite{KapSaf1} have constructed a parametrix modulo $\CI(\RR^n),$ but
without control over asymptotics at infinity.  More recently, Tataru
\cite{Tataru} has completed a construction of a frequency-localized
outgoing parametrix, valid for $C^2$ time-dependent coefficients that are
only rather weakly asymptotically flat; this construction, while not giving
a global description of the Schwartz kernel, suffices for obtaining
global-in-time Strichartz estimates.

\

The paper is divided into four parts. In the following section we give some
heuristic motivation for our geometric approach, particularly for the
choice of the space $X$ and the `scattering fibred structure' on it.  The
fundamental mathematical objects involved in this structure, namely the Lie
algebra of vector fields, scattering-fibred tangent and cotangent bundles,
and contact structures at the boundary, are introduced more formally in Section~\ref{two}.

In Part 2, we give the definitions of Legendre distributions of various
sorts. Unfortunately, although this follows a well-worn path (via
\cite{Ho:FIO1}, \cite{Melrose:Spectral}, \cite{MZ}, \cite{hasvas1} and
\cite{hasvas2}), there is little we can use directly from previous
literature, since we need to generalize to manifolds with codimension three
corners, so this part is rather long and technical. Each section follows a
similar template: we define the relevant Legendre submanifolds, explain how
to parametrize them, show that parametrizations always exist and the
equivalence of parametrizations, define Legendre distributions and give a
symbol calculus.  The reader should perhaps skip this part on a first
reading and return to it as needed.

In Part 3, we construct the semiclassical resolvent, thereby proving
Theorem~\ref{main1}, using the machinery from part 2, following
\cite{hasvas2} rather closely.

In Part 4, we prove Corollary~\ref{sp-cor}, Theorem~\ref{prop-thm},
Corollary~\ref{Poisson-cor} and Theorem~\ref{sc-thm}.

\emph{Acknowledgements.} We thank Andr\'as Vasy and Nicolas Burq for
illuminating discussions; we are grateful to Vasy for allowing some of the
fruits of his joint work \cite{Legdist} with A.H.\ to appear here.  This
research was supported in part by a Fellowship, a Linkage and a Discovery
grant from the Australian Research Council (A.H.) and by NSF grants
DMS-0100501 and DMS-0401323 (J.W.).

%%%%%%%%%%%%%%%%%%%%%%%%%%%%%%%%%

\section{Geometric motivation}\ilabel{fir}

Before getting into details we make some additional motivational remarks
about the geometric ingredients of this paper.

\subsection{The space $X$} 
\begin{figure}\centering
\epsfig{file=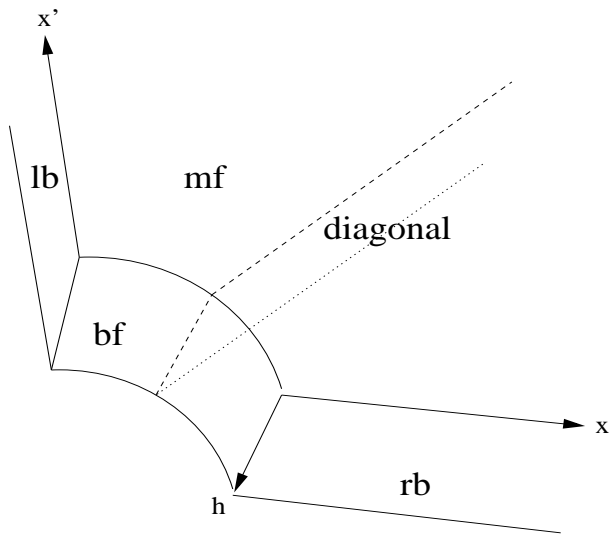,width=8cm,height=7cm}
\caption{The space $X$; in this figure dimensions in the direction of $\partial M$, in either factor, are not shown.}
\label{fig:X}
\end{figure}

In the Overview we introduced the space $X$, which is the blowup of $M \times M \times [0, h_0)$ at the corner $\partial M \times \partial M \times [0, h_0)$, or  in other words, 
$X = \MMb \times [0, h_0)$. The space $\MMb$ has boundary hypersurfaces $\lb = \partial M \times M$ (the left boundary), $\rb
  = M \times \partial M$ (the right boundary) and the blowup face $\bfc$
  (the `b-face'), which is a quarter-circle bundle over $(\partial M)^2$.
%We ignore the lack of compactness as $h \to h_0$ since we shall only be interested in functions on $X$ supported near $h=0$, say $h \leq h_0/2$.  
The boundary hypersurfaces of $X$ are
  then $\bfc \times [0, h_0), \lb \times [0, h_0), \rb \times [0, h_0)$ and
  $\MMb \times \{ 0 \}$. We shall denote these hypersurfaces (by abuse of
  notation) $\bfc, \lb, \rb$ and $\mf$. The diagonal in $X$ is the submanifold $\Delta_b \times [0, h_0)$, where $\Delta_b \subset \MMb$ is the lift of the diagonal in $M^2$ to $\MMb$. In a further abuse of notation we shall denote $\Delta_b \times \{ 0 \} \subset X$ simply by $\Delta_b$. See Figure~\ref{fig:X}. 
  
A total boundary defining function for a manifold with corners is, by definition, a product of defining functions for each boundary hypersurface. 
  The total boundary defining function for $X$ can be taken to be
  $\xx = h \rho$, where $\rho$,   a  total boundary defining function for $\MMb$, is given by $\rho^{-2} = x^{-2} + (x')^{-2}$. Here $x$ is a boundary defining function for $M$, lifted to $M^2$ by the left
  projection and then to $\MMb$ by the blowdown map, while $x'$ is the same boundary defining function on $M$ lifted via the right
projection. 

In this subsection we give some motivation for the choice of $X$ as the
space on which to analyze the kernel of the resolvent $(H - (\evosq \pm
i0))^{-1}$. We first point out that it allows us to decouple the diagonal
singularities from the long-range behaviour far from the diagonal, i.e.\
the lack of decay in the kernel at spatial infinity (at $\bfc$, $\lb$ and
$\rb$) and as $h \to 0$. Indeed, on $X$ the diagonal is separated from
$\lb$ and $\rb$, while it meets $\bfc$ transversally. This allows us to
solve for the resolvent kernel by first determining the conormal
singularity at the diagonal using standard pseudodifferential techniques,
and then solving away the remaining error as a separate step.

Consider the free semiclassical resolvent kernel on $\RR^n$, say for $n=3$,
$$
\frac1{4\pi h^2} \frac{e^{i\evosr |z-z'|/h}}{|z-z'|}.
$$ Let us ignore the diagonal singularity in the remainder of this section,
in view of the remarks above, for example by multiplying by a function on
$X$ that vanishes in a neighbourhood of $\Delta_b$.\ah{Define this
notation} Considered as a function on $X$, the resulting kernel is the
product of a function \emph{conormal} at the boundary of $X$, times an
explicit oscillatory factor $e^{i\evosr |z-z'|/h}$.  This would not be true
if the diagonal were not blown up, i.e.\ $\abs{z-z'}^{-1}$ is not conormal
at the boundary on the space $M^2 \times [0, h_0),$ meaning that it does
not have stable regularity under repeated application of vector fields
tangent to all boundary faces.  In this sense the singularities of the
resolvent kernel at the boundary (and away from the diagonal) are
`resolved' when lifted to the blowup space $X$.

More crucially, the blowup is needed so that we can analyze the resolvent kernel as a Legendre distribution at spatial infinity. A Legendre distribution of the simplest sort is given by an oscillatory function 
$$
e^{i\Phi/\xx} a
$$ where $\xx$ is the total boundary defining function for $X$ as above and
$a$ is conormal on $X$. The phase function $\Phi$ should be smooth (and
have certain properties with respect to fibrations at the boundary---see
the following subsection, and Section~\ref{leg-param}). The function $\Phi
= \Phi(z,z')$ is given, loosely speaking (and for $V \equiv 0$), by the
geodesic distance between $z$ and $z'$, at least in the region where this
is smooth; thus we want a compactification of $(M^\circ)^2$ where $\xx
d(z,z')$ is a smooth function up to the boundary (at least in this region,
and away from the diagonal). The b-double space $\MMb$ has this property
\cite{HTW}, and the blowup is essential here.

\subsection{Scattering-fibred structure}\ilabel{2sfs} The space $X$ comes
equipped with fibrations on its boundary hypersurfaces, and a corresponding
Lie algebra of vector fields, which dictate the type of Legendre
distributions we expect to find comprising the semiclassical
resolvent. This is dealt with in detail in Section~\ref{two}, but we give
an informal motivation here. We begin by noting the vector fields out of
which our operator is built. Near the boundary of $M$, the vector fields of
unit length with respect to our metric $g$ are $\CI(M)$-linear combinations
of the vector fields
$$ x^2 \partial_x \text{ and } x \partial_{y_i}.
$$
(Note that in polar coordinates on Euclidean space,  $\partial_r = -x^2 \partial_x$ and $\partial_{\omega}/r = x \partial_\omega$ are of approximately unit length as $r \to \infty$.) These vector fields generate the scattering Lie algebra of vector fields introduced by Melrose \cite{Melrose:Spectral}. In the semiclassical setting, we multiply each derivative by $h$, so we can think of our operator $H = h^2 \Delta$ (acting in either the left or the right set of variables) as being `built' out of the vector fields 
$$
h x^2 \partial_x , \quad h x \partial_{y_i}, \quad h (x')^2 \partial_{x'} , \quad h x' \partial_{y'_i}
$$
where the left set of variables is indicated without, and the right set with, a prime. 

Motivated by the program proposed by Melrose \cite{Melrose:Kyoto}, we should add one more vector field to this set in order to obtain $N = \dim X$ vector fields, so that it can generate a vector bundle which can be taken to replace the tangent bundle of $X$. It is not obvious what this extra vector field should be, but in hindsight we can observe that the vector field
$$
h (x \partial_x + x' \partial_{x'} - h \partial_h )
$$
fits the bill. In fact, on Euclidean space, the  self-adjoint operator corresponding to this is
$$
-ih \big( r \partial_r + r' \partial_{r'} + h \partial_h + n)
$$
and this \emph{annihilates} the semiclassical resolvent kernel on $\RR^n$ (this follows immediately from the fact that it is $h^{-n}$ times a function of $(z-z')/h$). 

We now have a set of vector fields generating a Lie algebra, and \emph{we can expect that the semiclassical resolvent has fixed regularity under the repeated application of these vector fields} (away from the diagonal). At $\mf \subset X$, i.e.\ at the interior of the $h=0$ face, we obtain from these vector fields all the scattering vector fields, i.e.\ those of the form
$$
h^2 \partial_h, \quad h \partial_{z_i}, \quad h \partial_{z'_j}
$$ where we use $z = (z_1, \dots, z_n)$ as a local coordinate on $M^\circ$,
here and throughout this paper. Thus the resolvent can be expected to be
Legendre at the interior of $\mf$ (equivalently, a semiclassical Lagrangian
distribution). At the other boundary hypersurfaces, the situation is a
little different. Our vector fields do not vanish at $\bfc$, $\lb$ or
$\rb$; rather they are tangent to the leaves of a fibration on each of
these boundaries. At $\bfc$, all the vector fields vanish except the last
one introduced above, which restricts to $h^2 \partial_h$ at $\bfc$.  At
$\rb$, the vector fields $h (x')^2 \partial_{x'}$ and $h x'
\partial_{y'_j}$ vanish, but the others restrict to $h \partial_{z_i}$ and
$h (x \partial_x - h \partial_h)$, which do not. These statements can be
rephrased by saying that on $\bfc$, the vector fields are constrained to be
tangent to the leaves of the fibration that projects off the $h$ factor,
while at $\rb$ the vector fields are constrained to be tangent to the
leaves of the fibration $\rb = M \times \partial M \times [0, h_0) \to
\partial M$ which projects to the second factor. We finally end up with a
characterization of our vector fields in terms of these boundary fibrations
and the total boundary defining function $\xx$ (see Definition~\ref{2.3},
and also Example~\ref{rem:ex}).  \emph{Our ansatz in this paper---justified
by Theorem~\ref{main1}---is that the semiclassical resolvent is Legendrian
with respect to this Lie algebra structure}, which we call the
scattering-fibred structure, on $X$.

\subsection{Sojourn relations}\ilabel{srt}
The sojourn time was introduced by Guillemin \cite{Gu}, motivated by a
result of Majda \cite{Majda}, in connection with metric or obstacle
scattering on $\RR^n$. Let $\gamma$ be a geodesic with asymptotic incoming
direction $y$ and asymptotic outgoing direction $y'$, ($y, y' \in
S^{n-1}$), and suppose that it is nondegenerate, meaning that locally it is
the only such geodesic (in a quantitative sense, so that the corresponding
Jacobian is nonzero). Guillemin defined the sojourn time $T(y, y')$ to be
the limit $l(R) - 2R$ where $l(R)$ is the length of the part of the
geodesic lying inside $B(R, 0)$. He then showed that the scattering matrix
locally took the form
\begin{equation}
S(\lambda, y, y') =  \sigma(y, y')^{-1/2} \lambda^{(n-1)/2} e^{i\lambda T(y, y')}  + O(\lambda^{(n-3)/2})
\ilabel{Gui}\end{equation}
(or a sum of such terms if there are finitely many such geodesics)
where $\sigma$ is a Jacobian factor. This has been generalized by Alexandrova, who removed the nondegeneracy assumption and proved that the scattering matrix is a semiclassical Fourier integral operator away from the diagonal. The Lagrangian to which the scattering matrix is associated, which Alexandrova calls the scattering relation, is parametrized by Guillemin's sojourn time whenever it is projectable, i.e.\ whenever $(y, y')$ locally form coordinates on it. Thus, the sojourn time is better thought of as a \emph{Lagrangian submanifold} rather than as a function.

In our Corollary~\ref{Poisson-cor} and Theorem~\ref{sc-thm} we see the sojourn time show up naturally, in two different guises. For simplicity we explain this in the case of zero potential. First we consider a geodesic emanating from a point $(z, \hat \zeta)$ in the cosphere bundle of $M^\circ$. By the nontrapping assumption, this geodesic $\gamma(s)$ tends to infinity as $s \to \infty$, and it does so in such away that the limit 
$$
\nu = \lim_{s \to \infty} s - r(\gamma(s))
$$ exists, where $s$ is arc-length along the geodesic and $r = 1/x$ is the
radial coordinate. In \cite{HW} we showed that there is a contact
transformation, which we called the sojourn relation, taking the point $(z,
\hat \zeta) \in S^* M^\circ$ to $(y_0', \nu, \mu)$, where $y_0'$ is the
asymptotic direction of $\gamma$ as $s \to \infty$ and $\mu$ is the
limiting value of $s^{-2} dy'/ds$ as $s \to \infty$. The image point
$(y_0', \nu, \mu)$ can be taken to lie in the boundary of the scattering
cotangent bundle\footnote{To be completely invariant it should be thought
of as lying in an affine bundle identified in \cite{HW}.} (see
Definition~\ref{sc-def}). We show in Corollary~\ref{Poisson-cor} that the
Poisson operator is a Legendre distribution associated to a Legendre
submanifold $\SR$ which is the twisted graph of the sojourn relation.  Just
as the Poisson operator is a boundary value of the resolvent kernel
(divided by $e^{ir'/h}$), so the sojourn relation appears as the `boundary
value' of the Legendrian $L$ associated to the resolvent (see
Section~\ref{subsec:sojourn}). The function $\nu$ appears as the boundary
value of $\psi - r'$ where $\psi$ is the function parametrizing $L$, with
the renormalizing term $r'$ coming directly from the removal of the
oscillatory factor $e^{ir'/h}$. Moreover, whenever $(z, y')$ locally form
coordinates on $\SR$, the function $\nu(z, y')$ locally parametrizes $\SR$.

When the point $z$ itself tends to infinity, say $z = \gamma(s')$ along a fixed geodesic $\gamma$, with $s' \to -\infty$,  the coordinate $\nu$ itself diverges as $1/r$ and we can take a limit
$$
\tau = \lim_{s' \to \infty} \nu - s' = \lim_{s, s' \to \infty} s - s' - r(\gamma(s)) - r(\gamma(s'))
$$
which is precisely Guillemin's sojourn time. We obtain the kernel of the scattering matrix as a boundary value of the Poisson operator, divided by $e^{ir/h}$, and in doing so, we find the total sojourn relation $T$ appearing  as the `boundary value' of the sojourn relation, with the sojourn time $\tau$ as the (renormalized) limit of $\nu$. Whenever $(y, y')$ locally form coordinates on $T$ (the nondegeneracy condition of Guillemin) then $\tau(y, y')$ locally parametrizes $T$, and we recover the description \eqref{Gui} of the scattering matrix. 

Our Theorem~\ref{sc-thm} improves upon results already in the literature in two ways. First, we treat (nontrapping)  asymptotically \emph{conic}, rather than flat, metrics, and second it is completely global. In particular, we do not need to localize away from the geodesics which are uniformly close to infinity (corresponding to the localization away from the diagonal in Alexandrova's result). Indeed it is this limiting regime which provides the transition between the Legendre (or semiclassical Lagrangian) behaviour of the scattering matrix in the limit $h \to 0$ and the Lagrangian behaviour of the scattering matrix for fixed $h$ as proved by Melrose-Zworski \cite{MZ}, since the latter is related to the geodesics `at infinity'. Our class of Legendrian-Lagrangian distributions unifies these two regimes into a single microlocal object.

%%%%%%%%%%%%%%%%%%%%%%%%%%%%%%%%%%

\section{Scattering-fibred structure}\ilabel{two}

In this section we shall define the scattering-fibred structure on
manifolds with corners. Although we only need the case of manifolds with
corners of codimension at most three, this structure can be defined on
manifold with corners of arbitrary codimension, and there is some
conceptual gain in considering the general case. So we shall give the basic
definitions for corners of arbitrary codimension, but rapidly specialize to
the case of codimension three corners for most of the exposition.  The
basic definitions are based on unpublished work \cite{Legdist} by the
first-named author and Andr\'as Vasy, and we thank him for permission to
use this material.  Note that the case of corners of
codimension two has been explicitly worked out in \cite{hasvas1}.  To
begin, we review the scattering structure on a manifold with boundary.

\subsection{Scattering structures}
Let $X$ be an $n$-dimensional manifold with boundary, and let  $x$ denote a boundary defining function on $X$. Denote by $\Vb(X)$  the Lie algebra of vector fields on $X$ tangent to $\pa X.$
\begin{defn}
The Lie algebra of scattering vector fields $\Vsc$ is defined by
\begin{equation}
V \in \Vsc(X) \text{ iff } V \in \Vb(X) \text{ and } V(x) \in x^2 \CI(X).
\end{equation}
\end{defn}
It is easy to verify that if $y$ are coordinates in $\pa X,$ extended to a collar neighbourhood of the boundary,  we may write a
scattering vector field locally near the boundary as a $\CI(X)$-linear combination of $x^2\pa_x$ and
$x\pa_{y_i}$, whilst away from the boundary a scattering vector field is simply a smooth vector field. It follows that $\Vsf(X)$ is the space of sections of a vector bundle over $X$.

\begin{defn}\ilabel{sc-def} We define $\Tsc(X)$, the \emph{scattering tangent
    bundle over} $X$, to be the vector bundle of which $\Vsc(X)$ is the
  space of sections; explicitly, the fibre $\Tsc_p(X)$ at $p \in X$ is
  given by $\Vsc(X)/I_p \cdot \Vsc(X)$, where $I_p(\Vsc(X))$ is the set of
  vector fields of the form $f V$, where $f \in \CI(X)$ vanishes at $p$
  and $V \in \Vsc(X)$.  We define $\Tscstar(X)$, the
  \emph{scattering cotangent bundle over} $X$, to be the dual vector
  bundle to $\Tsc(X)$.
\end{defn}

Locally near the boundary, the scattering cotangent space is spanned by the sections $d(1/x) = -dx/x^2$ and $dy_i/x$. Thus any point in $\Tsfstar X$ can be written 
$$
\nu d\big( \frac1{x} \big) + \sum_i \mu_i \frac{dy_i}{x}
$$
and this defines linear coordinates $(\nu, \mu_i)$ on each fibre of $\Tsfstar X$. In these coordinates, the natural symplectic form on $\Tsfstar X$ takes the form
$$
\omega = d\nu d\big( \frac1{x} \big) + \sum_i d \big( \frac{\mu_i}{x} \big) dy_i.
$$
There is a natural structure on $\Tscstar[\pa X] X$ defined by
contracting the symplectic form with $x^2 \pa x$ and restricting to the
boundary, taking the form
$$
d\nu - \sum_i \mu_i dy_i
$$
in these coordinates. 

Further details about scattering structures, and in particular of the
``scattering algebra'' of pseudodifferential operators that microlocalize
the scattering vector fields, can be found in \cite{Melrose:Spectral}.

\subsection{Scattering-fibred structures on manifolds with corners}\ilabel{subsection:horrific}
  Let $X$ be a compact manifold with corners of codimension $d$. %Let $M_1(X)$ denote the set of boundary hypersurfaces (i.e.\ codimension $1$ faces) of $X$.%, and let $\face(p), p \in X$, denote the highest codimension face containing $p$. 

\begin{defn}\ilabel{2.3} A scattering-fibred structure on $X$ consists of

(a) an ordering of the boundary hypersurfaces $\{ H_1, H_2, \dots , H_d \}$ of $M$, where we allow $H_i$ to be disconnected, i.e.\ to be a union of a disjoint collection of connected boundary hypersurfaces; 

(b) fibrations $\phi_{H_i} : H_i \to Z_i$, $1 \leq i \leq d$, to a compact manifold  $Z_i$ with corners of codimension $i-1$, and 

(c) a total boundary defining function $\xx$ (that is, a product of boundary defining functions $\prod_i \rho_i$ where $\rho_i$ is a boundary defining function for $H_i$)
which is distinguished up to multiplication by positive $\CI$ functions
which are constant on the fibres of $\pa X$.

The fibrations $\phi_i$ are assumed to satisfy the following conditions:

(i) if $i<j$, then $H_i \cap H_j$ is transverse to the fibres of $\phi_i$, and thus $\phi_i$
is a fibration from $H_i \cap H_j$ to $Z_i$, and

(ii) $H_i \cap H_j$ is a union of fibres of $\phi_j$ and thus $\phi_j$ is a
fibration from $H_i \cap H_j$ to $\pa_i Z_j\equiv \phi_j(H_i \cap H_j)$, where $\partial_i Z_j$ is a boundary hypersurface of $Z_j$. In
addition, 

(iii) there is a fibration $\phi_{ij} : \pa_i Z_j \to Z_i$ such that when
restricted to $H_i \cap H_j$, $\phi_i = \phi_{ij} \circ \phi_j$; in other words, there is a commutative diagram

\begin{equation}
\begin{diagram}
\node{H_i} \arrow{s,l}{\phi_i}
%\arrow{s,l}{c} \arrow{ese}
\node{H_i \cap H_j} \arrow{w,t}{inc}\arrow{e,t}{inc} \arrow{s,l}{\phi_j} 
\arrow{sw,t}{\phi_i}
\node{H_j} \arrow{s,r}{\phi_j} \\ %\arrow{wsw}\\
\node{Z_i} \node{\pa_i Z_j} \arrow{w,b}{\phi_{ij}}\arrow{e,b}{inc}
\node{Z_j}
\ilabel{CD}\end{diagram} .
\end{equation}

In this paper, we shall always assume the following additional condition:

(iv) The manifold $Z_d$ coincides with $H_d$ and the fibration $\phi_d$ is the identity map. 

The hypersurface $H_d$ will often be denoted $\mf$ (for `main face').
\end{defn}

We have a local model for this structure. Let $p \in M$ be a point on the codimension $d$ corner of $M$. 

\begin{proposition}\ilabel{X-coords}
Near $p$ there are local coordinates $x_1, \dots , x_d, y_1,
\dots, y_d$ where $x_i \geq 0$ is a boundary defining function for $H_i$ and 
$y_i$ lies in a neighbourhood of zero in $\RR^{d_i}$, such that $p = (0,
\dots, 0)$, and there are coordinates $(x_1, \dots, x_{i-1}, y_1, \dots,
y_i)$ on $Z_i$ near $\phi_i(p)$ such that,
locally, each $\phi_i$ takes the form
\begin{equation}
(x_1, \dots, x_{i-1}, x_{i+1}, \dots, x_d, y_1, \dots, y_d) \mapsto 
(x_1, \dots, x_{i-1}, y_1, \dots, y_i) .
\ilabel{local-coords}
\end{equation}
Moreover, the coordinates can be chosen so that $x_i$ is constant on the fibres of $H_j$ for $j > i$, and
    $\prod_i x_i = \xx$. 
\end{proposition}

\begin{proof}
We begin by choosing coordinates on the $Z_i$, in a neighbourhood of $\phi_i(p)$. We start with coordinates $y_1$ for $Z_1$, where $y_1$ lies in a neighbourhood of $0$ in $\RR^{k_1}$ and $y_1(\phi_i(p)) = 0$. Using the implicit function theorem, we may  choose coordinates $(y_1, y_2)$ on $\partial_1 Z_2$ so that the projection from $\partial_1 Z_2$ to $Z_1$ takes the form $(y_1, y_2) \to y_1$. We choose an arbitrary boundary defining function $x_1$ for $Z_2$ and extend the coordinates $(y_1, y_2)$ to a neighbourhood of $\partial_1 Z_2$, and in this way have coordinates $x_1, y_1, y_2$ on a neighbourhood of $\phi_2(p)$ in $Z_2$. Inductively, given coordinates $(x_1, \dots, x_{j-1}, y_1, \dots y_j)$ near $\phi_j(p)$ in $Z_j$, we choose coordinates on $\partial_{j} Z_{j+1}$ of the form $(x_1, \dots, x_{j-1}, y_1, \dots y_j, y_{j+1})$ so that the projection from $\partial_{j} Z_{j+1}$ to $Z_j$ is the coordinate projection off $y_{j+1}$. We then choose an arbitrary boundary defining function for $\partial_{j} Z_{j+1}$ and extend the coordinates from the boundary into the interior, and 
in this way have coordinates on a neighbourhood of $\phi_{j+1}(p)$ in $Z_{j+1}$. 

We can lift the coordinates from $Z_i$ to $H_i$ by the fibration $\phi_i$ in a neighbourhood of $p$. Thus $x_j$ and $y_{j}$ are defined on the union of $H_j, \dots, H_d$. These functions agree on intersections $H_j \cap H_k$ due to the way they are defined on $Z_i$ and due to the commutativity of the diagram \eqref{CD}. Hence they extend  to smooth functions on a neighbourhood of $p$. Finally we define
$x_d = \xx/(x_1 \dots x_{d-1})$ and all conditions are satisfied. 

\end{proof}

Thus, in the codimension three case, there are local coordinates near the corner in which the fibrations take the form
\begin{equation}\begin{aligned}
\phi_1 : (x_2, x_3, y_1, y_2, y_3) &\mapsto y_1 \\ \phi_2 : (x_1, x_3, y_1,
y_2, y_3) &\mapsto (x_1, y_1, y_2) \\ \id = \phi_3 : (x_1, x_2, y_1, y_2,
y_3) &\mapsto (x_1, x_2, y_1, y_2, y_3) \end{aligned}\ilabel{coords-3}\end{equation}
We proceed to give the main example of the scattering-fibred structure for the purposes of this paper.

\begin{example}\ilabel{23} Let $Y$ be a scattering-fibred manifold with codimension 2 corners. Thus $Y$ has two boundary hypersurfaces $K_1$ and $K_2$ with boundary defining functions $x_1, x_2$ together with fibrations $\psi_i : K_i \to Z_i$; moreover, $Z_2 = K_2$ and $\psi_2$ is the identity, while $Z_1$ is a manifold without boundary and the fibres of $\psi_1$ are transverse to the boundary. 
 
 Then the space $X = Y \times [0, \epsilon)_{x_3}$ is, in a natural way, a scattering-fibred manifold with codimension 3 corners. The boundary hypersurfaces are now $H_1 = K_1 \times [0, \epsilon)$, $H_2 = K_2 \times [0, \epsilon)$ and $H_3 = \mf = Y \times \{ 0 \}$. The structure is specified as follows: a distinguished boundary defining function is given by $\rho x_3$ where $\rho$ is a distinguished boundary defining function for $Y$; the bases of the fibrations are given by $Z_1$ and $Z_2$ and $Z_3 = Y$; and the fibrations are given by
 \begin{equation}\begin{gathered}
 \phi_1 : K_1 \times [0, \epsilon) \to Z_1 = \psi_1 \circ \Pi_1 \\
 \phi_2 : K_2 \times  [0, \epsilon) \to Z_2 = \Pi_2, \\
 \phi_3 = \id
 \end{gathered}\end{equation}
 where $\Pi_i : K_i \times [0, \epsilon) \to K_i$ is projection onto the first factor. It is easily checked that this satisfies all conditions of a scattering-fibred structure on $X$.

 (We remark that we are ignoring the fact that $X$ is noncompact, contrary
   to the above definition; this is harmless since we will in practice only be concerned with compactly supported distributions on $X,$ supported in say $x_3 \leq \epsilon/2$.)  
   
A special case of this is of course $Y = \MMb$, the b-double space of a manifold with boundary $M$; we have discussed this space already in Section~\ref{fir}. In this case, $H_1 = \lb \cup \rb$, $H_2 = \bfc$ and $H_3 = \mf$.
Let us further consider the boundary fibration structure in this case. 
The fibrations are given by the identity on $\mf$, by the
  projection off the $h$ factor on $\bfc$, and by the projection to
  $\partial M$ on $\lb$ and $\rb$. 
  
  Consider a point on the codimension
  three face $\bfc \cap \rb \cap \mf$, which is naturally diffeomorphic to
  $(\partial M)^2$. Recall that the total boundary defining function $\rho$ for $\MMb$ is given by $\rho = ( x^{-2} + (x')^{-2})^{-1/2}$ where $x$ is a boundary defining function for $M$ lifted via the left projection, and $x'$ is the lift of the same boundary defining function via the right projection.   
  Local coordinates near this point are $x_1 = \rho/x$, $x_2 = x$, $x_3 = h$, $y_1 = y'$, $y_2 = y$. (Notice that $x_1 = x'/x(1 + (x'/x)^2)^{-1/2}$ is equivalent to $x'/x$ for $x'/x$ small.) Then the fibrations take the form 
  \begin{align*}
\phi_1: (x, h, y', y) &\mapsto y' \text{ on } \rb, \\ \phi_2:
(x_1, h, y', y) &\mapsto (x_1, y', y) \text{ on } \bfc, \\ \phi_3:
(x_1, x, y', y) &\mapsto (x_1, x, y', y) \text{ on } \mf.
\end{align*}
Moreover, the product of the three boundary defining functions 
  satisfies
$$ x_1 \cdot x \cdot h = \xx ,  $$
so these coordinates satisfy the conditions of Proposition~\ref{X-coords}. 
    \end{example}

\subsection{Scattering-fibred tangent and cotangent bundles}
We return briefly to the case of corners of arbitrary codimension. 

\begin{defn} The space $\CIsf(X)$ is the space of $\CI$ functions $f$ on $X$
which are constant on the fibres of $\Phi$.
\end{defn}

\jw{Replaced a bunch of $x$'s with $\xx$'s in what follows.}
It is not hard to check that changing the total boundary defining function
$\xx$ to $f\xx$, where $f \in \CIsf(X) > 0$, leads to the same
scattering-fibred structure. Hence the total boundary defining function is
distinguished up to multiplication by elements of $\CIsf(X)$.

\begin{defn}\ilabel{sfvf} The Lie algebra of scattering-fibred vector fields $\Vsf$ is defined by
\begin{equation}
V \in \Vsf(X) \text{ iff } V \in \Vb(X), V(\xx) = O( \xx^2 ) \text{ and } V(f) = O(\xx) \text{ for all } f \in \CIsf(X).
%\text{ is tangent to } \Phi.
\end{equation}
Here we recall that $\Vb(X)$ is the Lie algebra of smooth vector fields on
$X$ which are tangent to each boundary hypersurface. We remark that the condition $V(f) = O(\xx)$ is equivalent to $V | H_i$ being tangent to the fibres of $\phi_i$. 
\end{defn}

It is easy to check that this {\it is} a Lie algebra. For if $V, W \in
\Vsf(X)$ then $V(W\xx) = V(\xx^2 g)$ for some smooth $g$, and this is $O(\xx^2)$
since $V$ is a b-vector field. Thus $[V,W]\xx = VW\xx - WV\xx =
O(\xx^2)$. Similarly, if $V(f) = O(\xx)$ and $W(f) = O(\xx)$ then  $[V,W]f = O(\xx)$. It is equally clear that $\Vsf(X)$ is invariant
under multiplication by smooth functions on $X$, and thus can be localized
in any open set.

Using coordinates as in Proposition~\ref{X-coords}, it may be checked that the Lie algebra $\Vsf(X)$ is the $\CI(X)$-span of the vector fields 
\begin{equation}\begin{aligned}
& (x_1x_2x_3 \dots x_d) x_1 \partial_{x_1}, & (x_1 x_2 x_3 \dots x_d) \partial_{y_1}, \\ 
& (x_2 x_3 \dots x_d) (x_1 \partial_{x_1} - x_2 \partial_{x_2}), & (x_2 x_3 \dots x_d) \partial_{y_2}, \\ 
& (x_3 \dots x_d) (x_2 \partial_{x_2} - x_3 \partial_{x_3}), & (x_3 \dots x_d)\partial_{y_3}, \\
& \vdots & \vdots \\
& x_d (x_{d-1} \partial_{x_{d-1}} - x_d \partial_{x_d}), & x_d \partial_{y_d}
\end{aligned}
\end{equation}
(where we write $\partial_{y_i}$ for the $k_i$-tuple of vector fields $\partial_{y_i^j}$, $1 \leq j \leq k_i$, if $\dim y_i = k_i$).
Thus, in the codimension three case, any vector field
in $\Vsf(X)$ is a linear combination  of 
\begin{equation}\begin{aligned}
& (x_1x_2x_3) x_1 \partial_{x_1}, & (x_1 x_2 x_3) \partial_{y_1}, \\ & x_2
x_3 (x_1 \partial_{x_1} - x_2 \partial_{x_2}), & (x_2 x_3) \partial_{y_2},
\\ & x_3 (x_2 \partial_{x_2} - x_3 \partial_{x_3}), & x_3 \partial_{y_3}.
\end{aligned}
\ilabel{sf-vfs}
\end{equation}
Therefore, locally near any point in $X$, the vector fields in $\Vsf(X)$ are arbitrary
linear combinations (over $\CI(X)$) of $N = \dim X$ vector fields. It
follows that $\Vsf(X)$ is the space of sections of a vector bundle over
$X$.

\begin{example}\ilabel{rem:ex} At the corner $\bfc \cap \rb \cap \mf$ of the space $X$ from Section~\ref{fir}, we have  $x_1 = x'/x$, $x_2 = x$, $x_3 = h$, $y' = y_1$, $y = y_2$; in these coordinates, we have
\begin{equation}\begin{aligned}
h (x')^2 \partial_{x'} &= (x_1 x_2 x_3) x_1 \partial_{x_1} \\
h x^2 \partial_x &= x_2 x_3 (x_1 \partial_{x_1} - x_2 \partial_{x_2})  
 \\
 h (x \partial_x + x' \partial_{x'} - h \partial_h) &= x_3 (x_2 \partial_{x_2} - x_3 \partial_{x_3})
 \end{aligned} \qquad \begin{aligned}
  h x' \partial_{y'} &= (x_1 x_2 x_3) \partial_{y_1}  \\
  h x \partial_y &= x_2 x_3 \partial_{y_2} 
\end{aligned}\end{equation}
so the vector fields arising in the discussion of Section~\ref{2sfs} generate the scattering-fibred Lie algebra. 
\end{example}

\begin{defn} We define $\Tsf(X)$, the \emph{scattering-fibred tangent
    bundle over} $X$, to be the vector bundle of which $\Vsf(X)$ is the
  space of sections; explicitly, the fibre $\Tsf_p(X)$ at $p \in X$ is
  given by $\Vsf(X)/I_p \cdot \Vsf(X)$, where $I_p(\Vsf(X))$ is the set of
  vector fields of the form $f V$, where $f \in \CI(X)$ vanishes at $p$
  and $V \in \Vsf(X)$.  We define $\Tsfstar(X)$, the
  \emph{scattering-fibred cotangent bundle over} $X$, to be the dual vector
  bundle to $\Tsf(X)$.
  
  We define $\Diffsf(X)$ to be the ring of differential operators generated by $\Vsf(X)$ over $\CI(X)$. 
\end{defn}

The vector bundle $\Tsfstar X$ is spanned by one-forms of the form $d(f/\xx)$
where $f\in\CIsf(M)$.  To see the duality between scattering-fibred vector
fields and differentials $d(f/\xx)$ for $f \in \CIsf(X)$, first observe that
there is a pairing between scattering-fibred vector fields and such
differentials for each $p \in X$ given by
\begin{equation}
\ang{d\left(\frac{f}{\xx}\right), V}_p =
V\left(\frac{f}{\xx}\right)(p). \ilabel{duality}
\end{equation}
This is finite for every $p \in X$ since $V(f) = O(\xx)$ and $V(\xx) =
O(\xx^2)$. In the codimension three case, choosing $f$ equal to
\begin{equation}\begin{gathered}
y_1^j, \quad x_1 y_2^j, \quad x_1 x_2 y_3^j, \\ 1, \quad x_1, \quad x_1
x_2,
\end{gathered}\ilabel{sf-f}\end{equation}
in turn, and pairing with the vector fields in \eqref{sf-vfs} gives a
non-degenerate matrix. Thus, we can identify the dual space of ${\Vsf}
_p(X)$, the scattering-fibred cotangent bundle at $p$, $\Tsfstar _p(X)$, as
\begin{equation}
\Tsfstar _p(X) = \big\{ d \big( \frac{f}{\xx} \big) \mid f \in \CIsf(X)
\big\} / \sim _p
\end{equation}
where $\sim _p$ is the equivalence relation of yielding the same pairing
\eqref{duality} at the point $p$.

The dual basis to the vector fields \eqref{sf-vfs} is
\begin{equation}
d \big( \frac{1}{x_1 x_2 x_3} \big), \quad d \big(\frac1{x_2 x_3} \big),
\quad d\big(\frac1{x_3}\big), \quad \frac{dy_1}{x_1 x_2 x_3}, \quad
\frac{dy_2}{x_2x_3}, \quad \frac{dy_3}{x_3}.  \ilabel{sf-cvs}
\end{equation}
Here $dy_i$ is shorthand for a $k_i$-vector of 1-forms, if $y_i \in
\RR^{k_1}$.  Any element of $\Tsfstar X$ may therefore be written uniquely
as
\begin{equation}
\nu_1 d\big(\frac{1}{x_1 x_2 x_3}\big) + \nu_2 d\big(\frac1{x_2 x_3} \big)
+ \nu_3 d\big(\frac1{x_3}\big) + \mu_1 \cdot \frac{dy_1}{x_1 x_2 x_3} +
\mu_2 \cdot \frac{dy_2}{x_2x_3} + \mu_3 \cdot \frac{dy_3}{x_3}.
\ilabel{can}\end{equation} The function $\nu_1$, regarded as a linear form
on the fibres of $\Tsfstar X$, can be identified with the vector field
$(x_1 x_2 x_3) x_1 \partial_{x_1}$, and similarly for the other fibre
coordinates. The same expression can be viewed as the canonical one-form on
$\Tsfstar X$. Taking $d$ of \eqref{can} therefore gives the symplectic form
on $\Tsfstar X$.

There is an alternative basis which is sometimes more convenient; instead of \eqref{sf-cvs} we use the basis
\begin{equation}
d \big( \frac{1}{x_1 x_2 x_3} \big), \quad  \frac{dx_1}{x_1x_2 x_3} ,
\quad \frac{dx_2}{x_2x_3}, \quad \frac{dy_1}{x_1 x_2 x_3}, \quad
\frac{dy_2}{x_2x_3}, \quad \frac{dy_3}{x_3}.  \ilabel{sf-cvs2}
\end{equation}
Using this basis, we can write any $q \in \Tsfstar X$ locally in the form
\begin{equation}\ilabel{nubarmu} q = \nubar_1 d\big( \frac1{x_1 x_2 x_3} \big) + \nubar_2
\frac{dx_1}{x_1x_2x_3} + \nubar_3 \frac{dx_2}{x_2x_3} + \mu_1
\frac{dy_1}{x_1x_2x_3} + \mu_2 \frac{dy_2}{x_2x_3} + \mu_3
\frac{dy_3}{x_3}.
\end{equation}
These are related to the $\nu_i$ by
\begin{equation}\begin{gathered}
\nubar_1 = \nu_1 + x_1 \nu_2 + x_1 x_2 \nu_3 \\ \nubar_2 = \nu_2 + x_2
\nu_3 \\ \nubar_3 = \nu_3.
\ilabel{nubar}\end{gathered}\end{equation}
In particular, $\nubar_1 = \nu_1$ at $x_1 = 0$ and $\nubar_2 = \nu_2$ at
$x_2 = 0$.

\subsection{Induced bundles and fibrations}

There is a natural subbundle of $\Tsfstar_{H_i}(X)$,\break namely\footnote{The restriction of $\Tsfstar X$ to a subset $S \subset X$ will be denoted $\Tsfstar_S X$} equivalence classes of differentials $d(f/\xx)$
where $f \in \CIsf(X)$ {\it vanishes} at $H_i$.
%(this means it vanishes at$H_j$ for all $j < i$). 
Let us denote this subbundle $\Tsfstar(F_i, H_i)$;
the reason for this notation will become evident below.  Notice that any $f
\in \CIsf(X)$ has a representation
\begin{equation}
f = f_1(y_1) + x_1 f_2(x_1, y_1, y_2) + x_1 x_2 f_3(x_1, x_2, y_1, y_2,
y_3) + \dots + x_1 x_2 x_3 \dots x_d \tilde f
\end{equation}
where $f_i$ and $\tilde f$ are smooth.  Thus the $i$th subbundle
corresponds to $f$ with $f_1 = \dots = f_{i} = 0$, while the $f_j$, $j >
i$, are arbitrary.  A point in the quotient bundle
$\Tsfstar_{H_i}(X)/\Tsfstar(F_i, H_i)$ is therefore given by a differential
$d(f/\xx)$ where only $f_1, \dots f_{i}$ are relevant. Since these
functions are constant on the fibres of $H$, they may be regarded as
functions on $Z_i$. Hence this is the lift to $H_i$ of a bundle over $Z_i$,
which we shall denote $\Nsf Z_i$.  Therefore there is an induced
fibration given by the composition
$$ \phit_i: 
\Tsfstar[H_i] X \to \Tsfstar[H_i] X/\Tsfstar(F_i, H_i) \to \Nsf Z_i.
$$
In the coordinates above, the subbundle $\Tsfstar(F_i, H_i)$ is given by $x_i = 0, \nu_1 = \dots = \nu_i = 0, \mu_1 = \dots = \mu_i = 0 $, while $(x_1, \dots, x_{i-1}, y_1, \dots, y_i, \nu_1, \dots, \nu_i, \mu_1, \dots, \mu_i)$ furnish coordinates on $\Nsf Z_i$ in a natural way. 
%In fact we will mostly be interested in the restriction of this fibration to $\Tsfstar[H_i \cap \mf], which we shall denote $\phit_i$. This maps onto $\Nsf Z_i$ since the fibres of $H_i$ are transverse to $\mf$ for $i < d$. 

The subbundle $\Tsfstar(F_i, H_i)$ can be interpreted as follows. 
We observe that each fixed fibre
$F_i$ of $H_i$ has an induced scattering-fibred structure, since $F_i$
meets $H_{i+1} \dots H_d$
and the fibrations $\phi_j$ for $j > i$ restrict to fibrations from $F \cap
H_j$ to a face of $Z_j$. Moreover, a total boundary defining function for
$F$ is given by $\xx/(x_1 \dots x_i)$, where $x_k$ for $k \leq i$ is chosen
to be constant on the fibres of $H_j$ for $j > k$. Then the bundle
$\Tsfstar(F_i, H_i)$ restricted to a single fibre $F$ of $H_i$ is naturally isomorphic to the
scattering-fibred cotangent bundle of $F$, $\Tsfstar F$.

For concreteness consider the codimension three case. Recall that 
$Z_1$ is a manifold without boundary, while $Z_2$ has a
boundary which we denote $\partial_1 Z_2$, and $Z_3 = \mf$ has two boundary
hypersurfaces which we denote $\partial_1 Z_3$ (the intersection with
$H_1$) and $\partial_2 Z_3$ (the intersection with $H_2$). 
Moreover, there is an induced fibration $\phi_{12} : \partial_1 Z_2 \to Z_1$, as in \eqref{CD}. We claim that  the commutative diagram \eqref{CD} with $i=1, j=2$ `lifts' to a commutative diagram at the level of cotangent spaces 
\begin{equation}
\begin{diagram}
\node{\Tsfstar[H_1] X} \arrow{s,l}{\phit_1}
%\arrow{s,l}{c} \arrow{ese}
\node{\Tsfstar[H_1 \cap H_2] X} \arrow{w,t}{inc}\arrow{e,t}{inc} \arrow{s,l}{\phit_2} 
\arrow{sw,t}{\phit_1}
\node{\Tsfstar[H_2] X} \arrow{s,r}{\phit_2} \\ %\arrow{wsw}\\
\node{\Nsf Z_1} \node{\partial_1 \Nsf Z_2} \arrow{w,b}{\phit_{12}}\arrow{e,b}{inc}
\node{\Nsf Z_2}
\ilabel{CD2}\end{diagram} .
\end{equation}

In this diagram everything has been explained except the existence and properties of the map $\phit_{12}$. To define it, note that the subbundle  $\Tsfstar[H_1 \cap H_2](F_2, H_2)$ is a subbundle of $\Tsfstar[H_1 \cap H_2] (F_1, H_1)$. Therefore there is an induced map on the quotient bundles, which we denote$$
\phit_{12} :  \partial_1 \Nsf Z_2 \to \Nsf Z_1,
$$
making the diagram \eqref{CD2} commute. 

(We remark that there is also a diagram analogous to \eqref{CD2} for $(i,j) = (1,3)$ or $(2,3)$ as well. In these cases, the maps $\phit_3$ is the identity, but the map $\phit_{i3}$, $i = 1,2$ is still of interest, mapping from $\Tsfstar[H_i \cap H_3] X$ to $\Nsf Z_i$.)

We shall often be interested in the restriction of the fibrations $\phit_i$ to $\Tsfstar_{H_i \cap \mf} X \to \Nsf Z_i$; notice that this is still onto since the fibres of $H_i$ are transverse to $\mf$, $i  < d$. We shall abuse notation slightly and call the restriction $\phit_i$ also. Thus, restriction to $\mf$ gives the following variant of \eqref{CD2}:
\begin{equation}
\begin{diagram}
\node{\Tsfstar[H_1 \cap H_3] X} \arrow{s,l}{\phit_1}
\node{\Tsfstar[H_1 \cap H_2\cap H_3] X} \arrow{w,t}{inc}\arrow{e,t}{inc} \arrow{s,l}{\phit_2} 
\arrow{sw,t}{\phit_1}
\node{\Tsfstar[H_2\cap H_3] X} \arrow{s,r}{\phit_2} \\ 
\node{\Nsf Z_1} \node{\partial_1 \Nsf Z_2} \arrow{w,b}{\phit_{12}}\arrow{e,b}{inc}
\node{\Nsf Z_2}
\ilabel{CD3}\end{diagram} .
\end{equation}

\begin{remark} Each space in the diagram
  above has a simple form in terms of the coordinates $x_i, y_i, \nu_i,
  \mu_i$. For example, the top left space is $\{ x_1 = x_3 = 0 \}$, the top
  middle space is $\{ x_1 = x_2 = x_3 = 0 \}$, the top right space is $\{
  x_2 = x_3 = 0 \}$, while on the bottom row the left space is $\{ x_1 =
  x_2 = x_3 = 0, y_2 = y_3 = 0, \mu_2 = \mu_3 = 0, \nu_2 = \nu_3 = 0 \}$,
  the middle space is $\{ x_1 = x_2 = x_3 = 0, y_3 = 0, \mu_3 = 0, \nu_3 =
  0 \}$ and the right space is $\{ x_2 = x_3 = 0, y_3 = 0, \mu_3 = 0, \nu_3
  = 0 \}$. Moreover, all the maps are the obvious coordinate projections or
  inclusions.
\end{remark}

\subsection{Contact structures}
In the remainder of this paper we restrict attention to the codimension
three case.  We define a 1-form $\chi$ on $\Tsfstar[\mf] X$ by contracting
the symplectic form $\omega$ with $\xx x_3 \partial_{x_3}$ (where $x_3$ is
a boundary defining function for $\mf$) and restricting to $\mf$. This
yields a contact structure (i.e.\ the form $\chi$ is non-degenerate in the
sense that $\chi \wedge (d\chi)^{N-1} \neq 0$, $N = \dim X$) in the
interior of $\mf$. However, this contact structure degenerates at the
boundary of $\mf$. In local coordinates \eqref{can}, the contact structure
takes the form
\begin{equation}\ilabel{contact}
\chi = d\nu_1 + x_1 d\nu_2 + x_1 x_2 d\nu_3 - \mu_1 \cdot dy_1 - x_1 \mu_2
\cdot dy_2 - x_1 x_2 \mu_3 dy_3
\end{equation}
and this degeneration is evident. Indeed, at $\Tsfstar[H_i \cap \mf] X$, $\chi$ vanishes on the subbundle
$\Tsfstar[H_i \cap \mf](F_i, H_i)$. However, it is not difficult to see that $\chi | \Tsfstar[H_i \cap \mf] X$ is the lift of a one-form from $\Nsf Z_i$. This is most easily seen in local coordinates; at $x_1 = 0$, $\chi = 
d\nu_1 - \mu_1 \cdot dy_1$ is the lift of a one-form $\chi_1$ from $\Nsf Z_1$ since it is expressible in terms of the coordinates $y_1, \nu_1, \mu_1$ which are the lifts of functions on $\Nsf Z_1$. Similarly,   at $x_2 = 0$, $\chi = d\nu_1 + x_1 d\nu_2  - \mu_1 \cdot dy_1 - x_1 \mu_2 \cdot dy_2$ is the lift of a one-form $\chi_2$ from $\Nsf Z_2$. Moreover, $\chi_1$ is nondegenerate, i.e.\ is a contact form, on $\Nsf Z_1$, while $\chi_2$ is nondegenerate except at $\partial_1 \Nsf Z_2$. 

In the coordinates \eqref{nubarmu} the contact form takes the form
\begin{equation}\ilabel{chibar}
\chi = d\nubar_1 - \nubar_2 dx_1 - x_1 \nubar_3 dx_2 - \mu_1 \cdot dy_1 - x_1 \mu_2
\cdot dy_2 - x_1 x_2 \mu_3 dy_3.
\end{equation}
These coordinates are more convenient when analyzing Legendre distributions (see Section~\ref{sec:legendrian}).

The degeneration of $\chi$ at $\Tsfstar[\mf \cap H_i] X$ and of $\chi_2$
on $ \Nsf_{\partial_1 Z_2} Z_2$ is captured by contact structures on the
fibres of the maps $\phit_i$ and $\phit_{12}$. To define these we make the following definition. 

\begin{defn} Suppose that $M$ is a manifold, $S \subset M$ a hypersurface with boundary defining function $s$,  and $\alpha$ a one-form on $M$ that vanishes at $S$. Thus $\alpha = s\beta$ for some one-form\footnote{Note that vanishing at $S$ is a strictly stronger condition than vanishing when restricted to $S$; e.g. $ds$ does not vanish at $S$ although it vanishes when restricted to $S$.} 
$\beta$. We call $\beta$ the \emph{leading part} of $\alpha$ at $S$. It is well defined up to multiplication by a nonzero function. This remains true even if $\alpha$ itself is only well-defined up to multiplication by a nonzero function. 
\end{defn} 

Notice that $\chi = \phit_1^* \chi_1$ at $\partial_1 \Nsf Z_3 \equiv \Tsfstar[H_1 \cap \mf] X$, that $\chi = \phit_2^* \chi_2$ at 
$\partial_2 \Nsf Z_3 \equiv \Tsfstar[H_2 \cap \mf] X$, and that
$\chi_2 = \phit_{12}^* \chi_1$ at $\partial_1 \Nsf Z_2$. Using the definition we can define   $\chi_{13}$ to be the leading part of $\chi -\phit_1^* \chi_1 $ at $ \Tsfstar[H_1 \cap \mf] X \subset \Tsfstar[\mf] X$,  $\chi_{23}$ to be the leading part of $(\chi - \phit_2^* \chi_2)/x_1$ at $ \Tsfstar[H_2 \cap \mf] X \subset \Tsfstar[\mf] X$ and $\chi_{12}$  to be the leading part of $\chi_2 - \phit_{12}^* \chi_1$ at $\partial_1 \Nsf Z_2 \subset \Nsf Z_2$. Using the invariance property in the last part of the definition, we see that these one-forms are well-defined up to multiplication by nonzero functions. 
In local coordinates, we have
$$\begin{gathered}
\chi_{12} = d\nu_2 - \mu_2 \cdot dy_2 ,\\
\chi_{23} = d\nu_3 - \mu_3 \cdot dy_3, \\
\chi_{13} = d\nu_2 - \mu_2 \cdot dy_2 + x_2 \big( d\nu_3 -
\mu_3 \cdot dy_3 \big) .
\end{gathered}$$
Hence we have well-defined contact structures (i.e.\ $\chi_{12}$ and $\chi_{23}$ are nondegenerate)   
on the fibres of $\phit_{12}$ and $\phit_{2}$ in \eqref{CD3}, while $\chi_{13}$ is nondegenerate on the fibres of $\phit_1$ for $x_2 > 0$.

%%%%%%%%%%%%%%%%%%%%%%%%%%%%%%%%%%%
\part{Machinery}
\section{Legendrian submanifolds and distributions} In this section we define Legendre distributions on a scattering-fibred manifold $X$ with corners of codimension 3. These will be smooth functions in the interior of $X$ which oscillatory behaviour at the boundary. 
\ilabel{sec:legendrian}

\subsection{Legendre submanifolds}\ilabel{leg}

\begin{defn}\ilabel{leg-sm-def} A Legendre
submanifold is a submanifold $G$ of dimension $N$ of $\Tsfstar[\mf] X$ on
which the contact form $\chi$  vanishes, and such that $G$ is
transverse to each boundary $\Tsfstar[\mf \cap H_i] X $ of $\Tsfstar[\mf] X$. 
\end{defn}

\begin{example} Let $f \in \CIsf(X)$. Then the graph of $d(f/\xx)$, restricted to $\Tsfstar[\mf] X$,  is a Legendre submanifold. The condition that $f \in \CIsf(X)$, as opposed to $\CI(X)$, is essential; see Section~\ref{leg-param}. 
\end{example} 

As a consequence of this definition, $G$ is well-behaved with respect to  the fibrations  $\phit_i : \Tsfstar[\mf \cap H_i] X  \to \Nsfstar Z_i$. To ease notation, we write $\partial_i G$ for the boundary hypersurface of $G$ lying over $H_i$, and $\partial_{12} G$ for the corner lying over $H_1 \cap H_2$. 

\begin{proposition}\ilabel{prop:fibre} 
(i) The restriction of $\phit_i$ to $\partial_i G$ is locally a fibration 
$$
\phi_i^G : \partial_i G \to G_i
$$
to an immersed Legendre submanifold $G_i \subset \Nsf Z_i$, and the fibres of $\phi_i^G$ are Legendre submanifolds for the contact structure for the fibres of $\phit_i$, i.e.\ for the contact form $\chi_{i3}$. 

(ii) The manifold $G_2$ is a manifold with boundary $\partial_1 G_2$. 
The restriction of $\phit_{12}$ to  $\partial_{1} G_2$ is locally a fibration
$$
\phi_{12}^G : \partial_1 G_2 \to G_1
$$
 and the fibres of $\phi_{12}^G$ are Legendre submanifolds for the contact structure for the fibres of $\phit_{12}$, i.e.\ for the contact form $\chi_{12}$. The maps form a commutative diagram
\begin{equation}
\begin{diagram}
\node{\partial_1 G} \arrow{s,l}{\phi^G_1}
\node{\partial_{12} G} \arrow{w,t}{inc}\arrow{e,t}{inc} \arrow{s,l}{\phi^G_2} 
\arrow{sw,t}{\phi^G_1}
\node{\partial_2 G} \arrow{s,r}{\phi^G_2} \\ 
\node{G_1} \node{\partial_1 G_2} \arrow{w,b}{\phi^G_{12}}\arrow{e,b}{inc}
\node{G_2}
\ilabel{CD4}\end{diagram} \qquad .
\end{equation}
Notice that each object in \eqref{CD4} is an element of the corresponding space in \eqref{CD3}, and the maps are induced from those in \eqref{CD3}. 
\end{proposition}

\begin{proof} For conceptual ease we first prove this result in the codimension two case. Thus suppose that $Y$ is a scattering-fibred manifold with codimension two corners. Near the corner, there are local coordinates $(x_1, x_2, y_1, y_2)$, such that the fibration $\phi_1$ on $H_1 = \{ x_1 = 0 \}$ takes the form $(x_2, y_1, y_2) \mapsto y_1$, while the fibration on the main face $H_2 = \{ x_2 = 0 \}$ is the identity. The contact form on $\Tsf_{H_2} Y$  is $\chi = d\nu_1 + x_1 d\nu_2 - \mu_1 \cdot dy_1 - x_1 \mu_2 \cdot dy_2$. 
Let $k_1 = \dim y_1$ and $k_2 = \dim y_2$. This local model applies everywhere except near $\partial_{12} G$, which we treat later. 

In the proof we shall need the following consequence of the implicit
function theorem: if $V$ is a compact manifold, $W$ is a manifold and $f :
V \to W$ is a smooth map of constant rank, then $f(V)$ is an immersed
submanifold of $Y$ and $f : V \to f(V)$ is (locally) a fibration.  

By assumption, $G$ is transversal to $\{ x_1 = 0 \}$ and the restriction of
$\chi$ to $G$ vanishes.   Given $p \in \partial_1 G,$ let $T_p( \text{fibre}
)$ denote the tangent space to the fibre of $\phit_1.$  Now consider the
space
$$
T_p \partial_1 G \cap
T_p( \text{fibre} );
$$ we claim that $d\nu_2 + \mu_2 \cdot dy_2 = 0$ restricted to this space
vanishes. To prove this, let $V$ be any vector in $T_p \partial_1 G \cap
T_p( \text{fibre} )$, and let $W$ be a vector tangent to $G$ and transverse
to $\{ x_1 = 0 \}$. Then $d\chi(V, W) = 0$. But $d\chi = -d\mu_1 \wedge dy_1
+ dx_1 \wedge (d\nu_2 - \mu_2 \cdot dy_2)$ at $\partial_1 G$. Since the
fibres of $\phit_1$ are given by $y_1$, $\nu_1,$ $\mu_1$ constant, it
follows that $(d\mu_1 \wedge dy_1)(V, W)$ vanishes. Also, $dx_1(V)$
vanishes, but $dx_1(W)$ does not. This forces $(d\nu_2 - \mu_2 \wedge
dy_2)(V)$ to vanish, which proves that the restriction of $d\nu_2 - \mu_2
\cdot dy_2$ to $T_p \partial_1 G \cap T_p( \text{fibre} )$ vanishes. Taking
the differential, we see also that $d\mu_2 \wedge dy_2 = 0$ vanishes when
restricted to $T_p \partial_1 G \cap T_p( \text{fibre} )$.

Now recall that coordinates on the fibres of $\phit_1$ are $(y_2, \nu_2,
\mu_2)$.   Since $d\mu_2 \wedge
dy_2 = 0$ on $T_p \partial G \cap T_p (
\text{fibre} ),$ the dimension of the projection of this space to the span
of the variables $\pa_{y_2},\pa_{\mu_2}$ is at most $k_2;$ since we further
have $d\nu_2 - \mu_2 \cdot dy_2 = 0,$ we in fact have
\begin{equation}
\dim \big(T_p(\partial G) \cap  T_p (\text{fibre}) \big) \leq k_2
\ilabel{dim1}\end{equation}
for any $p \in \partial G$. 

On the other hand, we can look at the projection of $\partial G$ onto
$\Nsfstar Z_1$, via $\phit_1$. We show that the rank of this map,
restricted to $\partial G,$ is at most $k_1$. For if not, then let $k >
k_1$ be the maximal rank of this map, and $p \in \partial G$ a point where
this maximum is attained.  Then the rank is exactly $k$ in a neighbourhood
of $p$. Using the implicit function theorem as above we see that the image
of $\partial G$ is locally a submanifold of dimension $k > k_1$. However,
the form $d\nu_1 + \mu_1 \cdot dy_1$ is zero on this image since it
vanishes on $\partial G$. Therefore the dimension of the projection of
$\partial G$ is Legendre and can have dimension at most $k_1$, which
contradicts $k > k_1$. It follows that
\begin{equation}
\text{rank } \phit_1 |_{ \partial G} = \dim( T_p (\partial G) ) - \dim
\big(T_p(\partial G) \cap T_p (\text{fibre}) \big) \leq k_1.
\ilabel{dim2}\end{equation} On the other hand, $\dim \partial G = k_1 +
k_2$, so the sum of the LHSs of \eqref{dim1} and \eqref{dim2} is everywhere
$k_1 + k_2$. Consequently, the dimension of $T_p(\partial G) \cap T_p
(\text{fibre})$ is exactly $k_2$ and the rank of $\phit_1 |_{ \partial G}$
is exactly $k_2$, and hence $\phit_1 : \partial G \to \Nsfstar Z_1$ has
constant rank $k_1$. By the implicit function theorem, the image of
$\partial G$ in $\Nsfstar Z_1$ is an immersed submanifold, which the
reasoning above shows is Legendrian; the fibres of this map are Legendre
submanifolds with respect to the contact structure on the fibres of
$\phit_1$.

Now we treat the codimension three case. The codimension two argument applies locally everywhere except for a neighbourhood of the corner $\partial_{12} G$ where we have be more careful. We claim that the manifold $G_2$ is transverse (and in particular, regular) up to the boundary of $\Nsf Z_2$. To prove this, we note that the implicit function statement above remains true if $V$ and $W$ are manifolds with boundary, provided that $f$ pulls back a boundary defining function for $W$ to a boundary defining function for $V$. The argument above that $\dim 
\big(T_p(\partial G) \cap  T_p (\text{fibre}({\phit_2})) \big) \leq \dim y_3$ is valid uniformly to the corner, but the argument on the base of the fibration does not extend automatically to the corner because the contact form $\chi_2$ on $\Nsf Z_2$ degenerates there. Instead, we must further analyze the structure of $G$ at the corner $\partial_{12} G$. Arguing as above, we see that for $p \in \partial_{12} G$,
$$
\chi_{23} \text{ vanishes on } T_p (\partial_{12} G) \cap T_p (\fibre(\phit_2)),
$$
$$
\chi_{12} \text{ vanishes on } (\phit_2)^* T_p (\partial_{12} G) \cap T_p (\fibre(\phit_{12})),
$$
and
$$
\chi_{1} \text{ vanishes on } (\phit_1)^* T_p (\partial_{12} G) .
$$
The dimension counting argument then shows that $\dim T_p (\partial_{12} G) \cap T_p (\fibre(\phit_2)) = \dim y_3$, $\dim (\phit_2)^* T_p (\partial_{12} G) \cap T_p (\fibre(\phit_{12})) = \dim y_2$ and $\dim (\phit_1)^* T_p (\partial_{12} G) = \dim y_1$ are all constant. This establishes the constancy of the rank of $\phit_2 : \partial_2 G \to \Nsf Z_2$ uniformly to the boundary and thus the regularity of $G_2$, as well as showing that $\partial_1 G_2$ fibres over $G_1$ with Legendrian fibres. 
\end{proof}

\begin{remark}Notice that, because of our assumption that the fibration at the main face
$\mf$ is the identity, the scattering-fibred structure locally near the
interior of the main face is the same as the scattering structure: locally,
we have $\Vsf(X) = \Vsc(X)$ near the interior of the main
face. Consequently, the theory coincides with the theory of Legendre
distributions as defined by Melrose and Zworski in the interior of $\mf$.
\end{remark}

\subsection{Parametrization}\ilabel{leg-param}

Before considering the general case let us consider the special case of
Legendrians $G$ which are projectable, meaning that the projection from $G
\subset \Tsfstar[\mf] X \to \mf$ is a diffeomorphism. In this case, $G$ is
necessarily given by the graph of the differential of a function. We claim that it is
necessarily of the form $f/(x_1x_2x_3)$, where $f \in \CIsf(X)$. In fact,
consider the graph of $d(f/(x_1x_2x_3))$ for a general smooth
$f$. Expanding this in the basis \eqref{sf-cvs2}, we find that the
coordinates $\nubar_i$ and $\mu_i$ are given by
\begin{gather*}
\nubar_1 = f - x_3 \partial_{x_3} f, \quad \nubar_2 = \partial_{x_1} f -
\frac{x_3}{x_2} \partial_{x_3} f, \quad \nubar_3 = \frac1{x_1} \partial_{x_2}
f - \frac{x_3}{x_1 x_2} \partial_{x_3} f \\ \mu_1 = \partial_{y_1} f ,
\quad \mu_2 = \frac1{x_1} \partial_{y_2} f , \quad \mu_3 = \frac1{x_1 x_2}
\partial_{y_3} f
\end{gather*}
For this to be a smooth submanifold, it follows that $\partial_{x_2} f$ and
$\partial_{y_2} f$ are $O(x_1)$ and $\partial_{x_3} f$ and $\partial_{y_3}
f$ are $O(x_1 x_2)$. Thus $f$ is of the form $$f = f_1(y_1) + x_1 f_2(x_1,
y_1, y_2) + x_1 x_2 f_3(x_1, x_2, x_3, y_1, y_2, y_3),$$ which is to say
that $f \in \CIsf(X)$. %Moreover in this case, $
%$$G_1 = \{ \nu_1 = f_1(y_1),\ \mu_1 = d_{y_1} f_1 \}$$ and
%$$G_2 = \{ \nu_1 = f_1(y_1) + x_1 f_2 + x_1^2 \partial_{x_1} f_1,\ \mu_1 =
%d_{y_1} f_1,\ \mu_2 = d_{y_2} f_2 \}.$$

Now consider the general case.  We will use the notation $\vecx,$ $\vecnu,$
$\vecy,$ $\vecmu,$ $\vecv$ respectively to denote the sets of coordinates
$(x_1,x_2,x_3),$ $(\nu_1,\nu_2,\nu_3),$ $(y_1,y_2,y_3),$
$(\mu_1,\mu_2,\mu_3),$ $(v_1,v_2,v_3).$ A (local) non-degenerate
parametrization of $G$ near a point $q \in G \cap \Tsfstar[H_1 \cap H_2
\cap \mf] X$ given in these coordinates as $q = (\vecx=0, \vecy^*, \vecnu^*,\ 
\vecmu^*)$ is a smooth function $\psi(\vecx,\vecy,\vecv)$ such that $\psi$
has the form
\begin{equation}\ilabel{eq:fib-ph-9}
\phi(\vecx,\vecy,\vecv)=\psi_1(y_1, v_1)+ x_1\psi_2(x_1, y_1, y_2, v_1,
v_2) + x_1 x_2 \psi_3
\end{equation}
such that $\psi_1$, $\psi_2$, $\psi_3$ are defined on neighborhoods of
$(y_1^*, v_1^*)$, $(0, y_1^*, y_2^*, v_1^*, v_2^*)$ and $q' = (0,0,0,
y_1^*, y_2^*, y_3^*, v_1^*, v_2^*, v_3^*)$ respectively with
\begin{equation}
d\big( \frac{\psi}{\xx} \big)(q')=q, \quad  d_{\vecv} \psi(q') = 0,
\end{equation}
$\psi$ is non-degenerate in the sense that
\begin{equation}\ilabel{eq:corner-20}
d_{(y_1, v_1)}\frac{\partial \psi_1}{\partial {v_1^i}},\ d_{(y_2, v_2)}
\frac{\partial \psi_2}{\partial {v_2^j}}, \ d_{(y_3, v_3)} \frac{\partial
\psi_3}{\partial {v_3^k}}
\end{equation}
are independent at $(y_1^*,v_1^*)$, $(y_1^*, y_2^*, v_1^*, v_2^*)$ and $q'$
respectively, and locally near $q$, $G$ is given by
\begin{equation}
G=\{d\big( \frac{\psi}{\xx} \big) \mid (\vecx,\vecy,\vecv) \in C_\psi\}
\ilabel{Gleg}\end{equation}
where
\begin{equation}\ilabel{eq:corner-22}
C_\psi=\{ (\vecx, \vecy, \vecv) \mid d_{\vecv} \psi = 0 \}.
\end{equation}
Note that the non-degeneracy conditions imply that $C_{\psi}$ is a smooth
submanifold of codimension $k_1 + k_2 +k_3$ of $X \times \RR^{k_1 + k_2 +
k_3}$, and that in the interior of $\mf$, the parametrization is
non-degenerate in the sense of \cite{MZ}. 

\begin{remark} We also have
\begin{multline}
\text{$\psi_1$ is a non-degenerate parametrization of $G_1$ } \\ \text{ and
$\psi_1 + x_1 \psi_2$ is a non-degenerate parametrization of $G_2$.}
\ilabel{subpar-1}\end{multline}
In addition, for fixed $(y_1, v_1)$ with $d_{v_1} \psi_1 = 0$, the phase
function $$\psi_1(y_1, v_1) + x_1 \psi_2(x_1, y_1, y_2, v_1, v_2)$$
parametrizes the fibres of the map $\phit_{12}$, while for fixed $(x_1,
y_1, y_2, v_1, v_2)$ with $d_{v_1, v_2} \psi_2 = 0$, the phase function
$\psi$ parametrizes the fibres of the map $\phit_{23}$.
\end{remark}

\subsection{Existence of parametrizations}

\begin{proposition}\ilabel{prop:param} Let $G$ be a Legendre submanifold. Then for any point
$q \in \pa G$ there is a non-degenerate parametrization of $G$ in some
neighbourhood of $q$.
\end{proposition}

\begin{proof} It is only necessary to do this in the case of a point $q$
  lying over $H_1 \cap H_2 \cap \mf$, since the other cases have already
  been proven in \cite{hasvas1}. By definition of a Legendre submanifold, the
  boundary $\partial_2 G$ of $G$ at $\{ x_2 = 0 \}$ fibres, via the map
  $\phit_{23}$, over $G_2$ with fibres that are Legendre submanifolds of
  $\Tscstar[\partial F] F$, where $F$ denotes a fibre of
  $H_2$. Coordinates on $\Tscstar[\partial F] F$ are $(y_3, \nu_3, \mu_3)$
  and, as in \cite{MZ}, Proposition~5, we can find coordinates $y_3 = (y_3^\flat, y_3^\sharp)$ near
  $\phi(q)$ so that  $(y_3^\sharp,
  \mu_3^\flat)$ form coordinates on the fibres of $\partial_2 G \to
  G_2$. In turn, $\partial_1 G_2$ fibres over $G_1$ with fibres that are
  Legendrian with respect to the contact form $\chi_{12} = d\nu_2 - \mu_2 \cdot dy_2$;
  hence we can find coordinates $y_2 = (y_2^\flat,
  y_2^\sharp)$ near $\phi(q)$ so that $(y_2^\sharp, \mu_2^\flat)$ form coordinates on the
  fibres of $\partial_1 G_2 \to G_1$. Lastly, since $G_1$ is Legendrian, we
  can find coordinates $y_1 =
  (y_1^\flat, y_1^\sharp)$ on $Z_1$ near $\phit_{13}(q)$ so that $(y_1^\sharp, \mu_1^\flat)$ form coordinates
  on $G_1$ locally. Using the transversality of $G$ to $\{ x_1 = 0 \}$ and
  $\{ x_2 = 0 \}$ we see that
$$ \cZ = (x_1, x_2, y_1^\sharp, y_2^\sharp, y_3^\sharp, \mu_1^\flat,
\mu_2^\flat,\mu_3^\flat)$$ form coordinates on $G$ near $q$. Consequently,
we can write the other coordinates as functions of these coordinates when
restricted to $G$.

We now use the coordinates \eqref{sf-cvs2} on the scattering cotangent bundle. The reason is that, in terms of a phase function $\Phi$ parametrizing a Legendrian
$G$, the value of $\nubar_1$ on $G$ is given simply by $\Phi$. The contact form is given by
\begin{equation}
d\nubar_1 - \nubar_2 dx_1 - x_1 \nubar_3 dx_2 - \mu_1 \cdot dy_1 - x_1
\mu_2 \cdot dy_2 - x_1 x_2 \mu_3 \cdot dy_3.
\ilabel{contactbar}\end{equation}

Writing $\nubar_i$, $y_i^\flat$ and $\mu_i^\sharp$ in terms of the
coordinates $\cZ$ on $G$, we have
\begin{equation}\begin{aligned}
\nubar_1 &= N_1(\cZ) \\ \nubar_2 &= N_2(\cZ) \\ \nubar_3 &= N_3(\cZ) \\
y_i^\flat &= Y_i^\flat(\cZ), i = 1 \dots 3 \\ \mu_i^\sharp &=
M_i^\sharp(\cZ), i = 1 \dots 3
\end{aligned} \qquad  \text{ on } G.
\end{equation}
Since $G$ is Legendrian, we have
\begin{equation}\begin{gathered}
dN_1 - N_2 dx_1 - x_1 N_3 dx_2 - \mu_1^\flat \cdot dY_1^\flat - M_1^\sharp
\cdot dy_1^\sharp \\ - x_1 \big( \mu_2^\flat \cdot dY_2^\flat - M_2^\sharp
\cdot dy_2^\sharp \big) - x_1 x_2 \big( \mu_3^\flat \cdot dY_3^\flat -
M_3^\sharp \cdot dy_3^\sharp \big) = 0.
\end{gathered}\ilabel{contact0}\end{equation}
We claim that the function
\begin{equation}
\Phi = N_1 + (y_1^\flat - Y_1^\flat) \cdot \mu_1^\flat + x_1 \big(
(y_2^\flat - Y_2^\flat) \cdot \mu_2^\flat \big) + x_1 x_2 \big( (y_3^\flat
- Y_3^\flat) \cdot \mu_3^\flat \big)
\ilabel{Phidefn}\end{equation}
is a local parametrization of $G$. To avoid confusion, let us write $v_i$
instead of $\mu_i^\flat$ for the corresponding arguments of $\Phi$.

First, observe that $N_1$ has the form $N_1 = N_{1,1}(y_1, v_1) + O(x_1)$
since at $x_1 = 0$, $\partial_1 G$ fibres over $G_1$ where the value of
$\nubar_1 = \nu_1$ is determined by $(y_1^\sharp, v_1)$ since these are
coordinates on $G_1$. Similarly, $N_1$ is a function of $(x_2, y_1, y_2,
v_1, v_2)$ plus $O(x_1 x_2)$, $Y_1^\flat$ is a function of $(y_1, v_1)$
plus $O(x_1)$, etc. It follows that $\Phi$ has the form
\eqref{eq:fib-ph-9}.

Second, suppose that $d_{v_1} \Phi = 0$. This means that
\begin{equation}\begin{gathered}
d_{v_1} N_1 + y_1^\flat - Y_1^\flat - d_{v_1} Y_1^\flat \cdot v_1 \\ - x_1
\big( d_{v_1} Y_2^\flat \cdot \mu_2^\flat \big) - x_1 x_2 \big( d_{v_1}
Y_3^\flat \cdot \mu_3^\flat \big) = 0.
\end{gathered}\end{equation}
Using the $dv_1$ component of \eqref{contact0}, this is the same thing as
saying that $y_1^\flat = Y_1^\flat$. In a similar way, the conditions that
$d_{v_i} \Phi = 0$ imply that $y_i^\flat = Y_i^\flat$, $i = 2, 3$. This
also shows the non-degeneracy condition, since $d(\partial_{v_j^i} \Phi) =
dy_j^i$ at $q$ which are manifestly linearly independent differentials.

To see that the set
$$ G' = \{ d \big( \frac{\Phi}{x_1 x_2 x_3} \big) \mid d_{v_1, v_2, v_3}
\Phi = 0 \}
$$ coincides with $G$ locally near $q$, first consider the value of
$\mu_1^\flat$; this is given by $d_{y_1^\flat} \Phi = v_1$. So we can
re-identify $\mu_1^\flat$ with $v_1$. Similarly we can re-identify
$\mu_2^\flat$ with $v_2$ and $\mu_3^\flat$ with $v_3$.

Next consider the value of $\nubar_1$ on $G'$. It is given by the value of
$\Phi$, that is, by \eqref{Phidefn}. This simplifies to $N_1$ when $d_{v_i}
\Phi = 0$, since we have $y_i^\flat = Y_i^\flat$ when $d_{v_i} \Phi =
0$. Now consider the value of $\nubar_2$. This is given by $d_{x_1} \Phi$
which is equal to
$$ d_{x_1} N_1 - d_{x_1} Y_1^\flat \cdot \mu_1^\flat - x_1 d_{x_1}
Y_2^\flat \cdot \mu_2^\flat - x_1 x_2 d_{x_1} Y_3^\flat \cdot \mu_3^\flat
$$ (again using $y_i^\flat = Y_i^\flat$ when $d_{v_i} \Phi = 0$). Since the
$dx_1$ component of \eqref{contactbar} vanishes, this is equal to $N_2$. So
$\nubar_2 = N_2$ on $G'$. In a similar way we deduce that $\nubar_3 = N_3$,
and $\mu_i^\sharp = M_i^\sharp$ on $G'$. It follows that $G'$ coincides
with $G$.
\end{proof}

\subsection{Equivalence of phase functions}

In this section we shall give a necessary and sufficient condition for equivalence of two phase functions parametrizing a given Legendrian. 
This is the key step in showing, in the following subsection, that the class of Legendre distributions does not depend on the choice of phase function, which is crucial for deducing that the class of Legendre distributions has a useful symbol calculus. 

Two phase functions $\phi$, $\tilde{\phi}$ are said to be \emph{equivalent}
if they have the same number of phase variables of each type $v_1,v_2,v_3$
and there exist maps
$$ V_1(\vecx,\vecy,\vecv),\ V_2(\vecx,\vecy,\vecv),\ V_3(\vecx,\vecy,\vecv)
$$ such that
$$ \tilde{\phi}(\vecx,\vecy,V_1,V_2,V_3) = \phi.
$$

\begin{proposition}\ilabel{prop:equiv1}
The phase functions $\phi=\psi_1+x_1 \psi_2 + x_1x_2\psi_3$ and
$\tilde\phi=\tilde\psi_1+x_1 \tilde\psi_2 + x_1x_2\tilde\psi_3$ are locally
equivalent iff
\begin{enumerate}
\item
They parametrize the same Legendrian,
\item
They have the same number of phase variables of the form $v_1,$ $v_2,$ and
$v_3$ separately,
\item
\begin{align*}
\sgn d_{v_1}^2 \psi_1 &= \sgn d_{v_1}^2 \tilde\psi_1,\\ \sgn d_{v_2}^2
\psi_2 &= \sgn d_{v_2}^2 \tilde\psi_2, \\ \sgn d_{v_3}^2 \psi_3 &= \sgn
d_{v_3}^2 \tilde\psi_3.
\end{align*}
\end{enumerate}
\end{proposition}
\begin{proof}
The proof follows \cite{Ho:FIO1}, Theorem~3.1.6 quite closely (and
Lemma~4.5 of \cite{hasvas1} even more so), hence we will be brief.  To
begin, we let $C$ and $\tilde C$ denote the respective sets where
$d_{\vecv} \psi=0,$ $d_{\vecv} \tilde\psi=0,$ near a given point in the
codimension-three corner.

We begin by noting that when we restrict to the face $H_1$ we have a phase
function $\psi = \psi_1(y_1,v_1)$ parametrizing $G_1.$ Hence by the usual
argument for equivalence of phase functions (\cite{Ho:FIO1}, as extended to
Legendrians in \cite{MZ}), there exists a fiber diffeomorphism $\tilde v_1
= V_1(y_1,v_1)$ such that $\psi_1(y_1,\tilde v_1) = \tilde
\psi_1(y_1,v_1).$ Furthermore on the face $H_2,$ equivalence of phase
functions is guaranteed by \cite{hasvas1}.  Hence we need only extend from
$H_1$ and $H_2$ to obtain equivalence on $H_3$ as well.

The manifolds $C_\psi$ and $C_{\tilde\psi}$ are diffeomorphic, via their
common fiber-preserving diffeomorphism with the Legendrian they
parametrize.  As they are smooth manifolds, we may extend this
diffeomorphism to a fiber-preserving diffeomorphism $F$ of an open
neighborhood of $C_\psi$ with an open neighborhood of $C_{\tilde\psi}.$
Then the phase function $\bar\psi:= F^*(\tilde\psi)$ has the property that
$C_{\bar\psi}=C_\psi=:C,$ and $\psi=\bar\psi$ to second order along $C.$
Therefore we have reduced by this initial change of variables to the case
in which we may take $\psi,\tilde\psi$ equal to second order along $C.$

We now improve this result to exact equivalence of $\psi$ and $\tilde \psi$
on $H_3,$ under the assumption that the functions agree to second order on
$C.$ As we have equivalence on $H_1,H_2$ we may write
$$ \psi = \tilde \psi_1+x_1 \tilde \psi_2 +x_1x_2  \psi_3.
$$ We may expand in a Taylor series on $H_3$:
$$ \tilde \psi_3 - \psi_3 = \frac 12 (\nabla'_{\vecv}\psi)^t B
(\nabla'_{\vecv}\psi)
$$ for some matrix $B=B(\vecx,\vecy,\vecv),$ where we define $\nabla' \psi
= (\pa_{v_1} \psi, \pa_{v_2}(\psi_2+x_2 \psi_3), \pa_{v_3} \psi_3).$
Observe as in \cite{Ho:FIO1} that the non-degeneracy assumptions on $\psi_3,
\tilde \psi_3$ means precisely that $\det (I+B_{33}\pa_{v_3
v_3}^2\psi_3)\neq 0$ where $B_{33}$ is the $(3,3)$ block of the matrix $B$.  We now expand
$$ \psi(\vecx,\vecy,\tilde{\vecv}) - \psi(\vecx,\vecy,{\vecv}) =
(\tilde{\vecv}-\vecv ) \cdot \pa_{\vecv}\psi + O((\tilde{\vecv}-\vecv )^2).
$$ Set
$$ (\tilde v_1, \tilde v_2,\tilde v_3) -(v_1,v_2,v_3) = (x_1 x_2 w_1, x_2
w_2, w_3)\cdot \nabla'_{\vecv}\psi
$$ where $w_i=w_i(\vecx,\vecy,\vecv)$ is a matrix for $i=1,2,3.$ We thus
have
$$ \psi(\vecx,\vecy,\tilde{\vecv}) - \psi(\vecx,\vecy,\vecv) = x_1 x_2
(\nabla'_{\vecv}\psi)^t (w+O(w^2)) (\nabla'_{\vecv}\psi).
$$ We want
\begin{align*} \psi(\vecx,\vecy,\tilde{\vecv}) - \psi(\vecx,\vecy,{\vecv}) &=
\tilde\psi(\vecx,\vecy,\vecv) -\psi(\vecx,\vecy,\vecv)\\ &=x_1 x_2
(\tilde\psi_3(\vecx,\vecy,\vecv) -\psi_3(\vecx,\vecy,\vecv)).
\end{align*}
We thus need to solve
$$ x_1 x_2 (\nabla'_{\vecv}\psi)^t (w+O(w^2)) (\nabla'_{\vecv}\psi)=
\frac{x_1 x_2}2 (\nabla'_{\vecv}\psi)^t B (\nabla'_{\vecv}\psi)
$$ for $w.$ This can be accomplished for $B$ small, i.e.\ for $\psi_3$ and
$\tilde \psi_3$ close, by the inverse function theorem; Lemma~3.1.7 of
\cite{Ho:FIO1} enables us to extend to the case of arbitrary $\psi_3,$
$\tilde\psi_3$ using the hypotheses on the signatures of $\pa_{v_3 v_3}^2
\psi_3$ and $\pa_{v_3 v_3}^2 \tilde\psi_3.$
\end{proof}

\subsection{Legendrian distributions}

Let $m,r_1,r_2$ be real numbers, let $N = \dim X$, let $G \subset \Tsfstar[\mf] X$ be a Legendre submanifold, and let $\nu$ be a smooth
nonvanishing scattering-fibred half-density. The set of (half-density) Legendre distributions of
order $(m; r_1, r_2)$ associated to $G$, denoted $I^{m, r_1, r_2}(X, G; \sfOh)$, is the set of half-density distributions that can be written in 
the form $u_1 + u_2 + (u_3 + u_4 + u_5) \nu$, such that  
\begin{itemize}
\item $u_1$ is a Legendre distribution of order $(m; r_1)$ associated to $G$ and
supported away from $H_2$, 
\item $u_2$ is a Legendre distribution of order $(m; r_2)$ associated to $G$ and supported away from $H_1$ (both of these are defined in \cite{hasvas1}), 
\item $u_3$ is given by an finite sum of
local expressions of the form
\begin{equation}\begin{gathered}\ilabel{eq:dist-1}
\int \int \int e^{i\psi(x_1, x_2, \vecy,\vecv)/\xx} a(\vecx,\vecy,\vecv)\\
x_3^{m-(k_1+k_2+k_3)/2+N/4} x_2^{r_2 - (k_1 + k_2)/2 - f_2/2 + N/4}
x_1^{r_1-k_1/2-f_1/2 + N/4}\,dv_1 \, dv_2 \, dv_3,
\end{gathered}\end{equation}
with $v_i \in \RR^{k_i},$ $a$ smooth and compactly supported,  $f_i$ the dimension of the
fibres of $H_i$ and $\psi = \psi_1 + x_1 \psi_2 + x_1 x_2 \psi_3$ a phase function locally parametrizing $G$ near a corner point $q \in \partial_{12} G$, as in Section~\ref{leg-param}, 

\item $u_4$ is  given by  a finite sum of terms of the form
\begin{equation}\begin{gathered}
 \int \int e^{i(\psi_1 + x_1 \psi_2)/\xx} b(x_1,
 x_2, x_3,  y_1, y_2, y_3, v_1, v_2) \\ x_2^{r_2 - (k_1 + k_2)/2 - f_2/2 + N/4}
 x_1^{r_1-k_1/2-f_1/2 + N/4}\,dv_1 \, dv_2
\end{gathered}\ilabel{eq:dist-2}
\end{equation}
with $\psi_1, \psi_2$ and $f_i$ as above, $b$ smooth with support compact and $O(x_3^\infty)$ at $\mf$, and

\item $u_5 \in \CIdot (X)$. (We use the notation $\CIdot(X)$ for $\xx^\infty \CI(X)$.)
\end{itemize}

\begin{remark} The convention regarding orders is as follows: the order
  increases as the distribution gets `better', i.e.\ vanishes more rapidly,
  and it is `zeroed' so that $N/4$ is critical for $L^2$-membership,
  i.e.\ for a distribution with positive symbol, $u$ is in $L^2$ iff all the
  orders are more than $N/4$. This somewhat peculiar choice is to conform
  to the order convention for pseudodifferential operators (apart from the
  change of sign) on a manifold of dimension $n$, whose kernels are in
  $L^2$ provided the order is less than $-n/2 = -N/4$, where  $N = 2n$ is
  the dimension of the space on which the kernel is defined. In any case,
  the order convention agrees with that of \cite{MZ}, \cite{hasvas1} and
  \cite{hasvas2}.  
\end{remark}

\begin{proposition}\ilabel{prop:any} Let $u \in I^{m,r_1,r_2}(X, G; \sfOh)$ be a Legendre distribution, and let $\psi$ be any local parametrization of some subset $U \subset G$. After localization to $U$, the $u$ may be expressed as an oscillatory integral with respect to $\psi$, modulo $\CIdot(X)$. 
\end{proposition}

\begin{proof}
We give a brief sketch of this proof, which follows standard lines.

By definition, $u$ can be written with respect to \emph{some} phase function parametrizing $G$, say $\psi'$. 

One can modify any phase function (without changing the Legendrian parametrized) by adding a nondegenerate quadratic form
$Q_1(w_1)+ x_1 Q_2(w_2) + x_1 x_2 Q_3 (w_3)$ in extra variables $w_i \in \RR^{l_i}$. This does not change, modulo $O(\xx^\infty)$, the distributions that can be written with respect to the phase function since the extra oscillatory factor only contributes a factor 
$$
c x_1^{l_1/2} x_2^{(l_1 + l_2)/2} x_3^{(l_1 + l_2 + l_3)/2}
$$
which is just an adjustment of the orders. However it allows us to change the number of phase variables of each type, and the corresponding signature. By modifying both $\psi$ and $\psi'$ in this way we may arrange that they satisfy the conditions of Proposition~\ref{prop:equiv1}. (This requires some mod 2 compatibility conditions between $\dim v_i$ and the signature of $d^2_{v_i v_i} \psi_i$ but these are automatically satisfied; see Theorem 3.1.4 of \cite{Ho:FIO1}.) One can then use the change of variables given by Proposition~\ref{prop:equiv1} to write $u$ in terms of the modified phase function $\psi$, and therefore in terms of $\psi$ itself. 
\end{proof}

\subsection{Symbol calculus}

The previous proposition implies that there is a symbol calculus for Legendre distributions. Since this follows standard lines, we omit the proof.
 
Let $X$ be a scattering-fibred manifold with codimension $3$ corners, let
$N = \dim X$, and let $G$ be a Legendre submanifold. Let $\xx$ denote the  distinguished total
boundary defining function for $X$, and $x_1,x_2,x_3$ be the set of
boundary defining functions for each $H_i \in M_1(X) \setminus \{ \mf
\}$. The Maslov bundle $M$ and the $E$-bundle are defined via the scattering
structure over the interior of $G$ and extend to smooth bundles over the
whole of $G$ (that is, they are smooth up to each boundary of $G$ at $\Tsf
_{H_i \cap \mf} X$); see \cite{hasvas2}. We define $N^*_{\mf} \partial X$ to be the bundle over $\mf$ given by differentials $df$ of smooth functions $f$ on $X$ vanishing at each boundary hypersurface. It is a line bundle with nonzero section $d\xx$. 

We define the symbol bundle $S^{[m]}(G)$ of order $m$ over $G$ to be the bundle 
\begin{equation}
S^{[m]}(G) = M(G) \otimes E \otimes \big| N^*_{\mf} \partial X \big|^{m-N/4},
\ilabel{symbolbundle}\end{equation} following \cite{hasvas2}.

\begin{proposition}\ilabel{prop:symb-calc} The symbol map for Legendre distributions, defined in the interior of $G$ \cite{MZ}, extends by continuity to 
give an exact sequence
$$ 0 \to I^{m+1, r_1, r_2}(X, G; \sfOh) \to I^{m, r_1, r_2}(X, G; \sfOh) \to
x_1^{r_1-m} x_2^{r_2 - m} \CI(G, \Omega^\half_b \otimes S^{[m]}(G)) \to 0.
$$ 
If $P \in \Diffsf(X; \sfOh)$ has principal symbol $p$ and $u \in I^{m,r_1, r_2}(X, G; \sfOh)$,
then $Pu \in I^{m,r_1,r_2}(X, G; \sfOh)$ and
$$ \sigma^m(Pu) = \big( p \rest G \big) \sigma^m(u).
$$ Thus, if $p$ vanishes on $G$, then $Pu \in I^{m+1,r_1, r_2}(X, G; \sfOh)$. The symbol of order $m+1$ of $Pu$ in this case is given
by\jw{Changed $\tau$ to $\nu_1$ here.  Correct?}
\begin{equation}
\Big( -i \mathcal{L}_{\scH_p} -i \big(\half + m - \frac{N}{4} \big) \frac{ \pa
  p}{\pa \nu_1} + p_{\sub} \Big) \sigma^m(u) \otimes |d\xx|,
  \ilabel{transport}\end{equation} where $\scH_p$ is the scattering Hamilton vector field
  of $p$ (that is, the Hamilton vector field multiplied by $\xx^{-1}$ and restricted to $G$), $\nu_1$ is the coordinate in the coordinate system \eqref{can}, and
  $p_{\sub}$ is the subprincipal symbol of $P$.
\end{proposition}

\begin{remark} The subprincipal symbol of a differential operator has the following properties: (i) for a multiplication operator $f$,  it is the $O(\xx)$ part of the Taylor series of $f$ at $x_3 = 0$. (ii) 
The subprincipal symbol of  $i(V - V^*)$,  where $V$ is a real vector field,  is zero. (iii) The subprincipal symbol of  the composition of two differential operators $P$ and $Q$ is $\sigma(P) \sigma_{\sub}(Q) + \sigma(Q) \sigma_{\sub}(P) -i/2 \{ \sigma(P), \sigma(Q)\}$. These properties in fact  uniquely determine the subprincipal symbol for any differential operator. 
\end{remark}

\begin{example} A very simple example may help to illustrate the symbol calculus. Let $P$ be the differential operator $x_1 x_2 x_3 (x_1 D_{x_1})$, $D = -i \partial$,  and let $u$ be the Legendre distribution
$$
u = x_3^{m + N/4} x_2^{r_2 - f_2/2 + N/4} x_1^{r_1 - f_1/2 + N/4} \Big| \frac{dx_1 dx_2 dx_3 dy_1 dy_2 dy_3}{x_3^{N+1} x_2^{N+1-f_2} x_1^{N+1-f_1}} \Big|^{1/2},
$$
a distribution of order $(m, r_1, r_2)$ associated to the zero section (which is a Legendrian submanifold). We assume that the half-density factor above, which is a smooth nonvanishing scattering-fibred half-density, is covariant constant. Hence $Pu = -i(r_1 - f_1/2 + N/4)x_1 x_2 x_3 u$. 

In terms of the symbol calculus, the symbol of $P$ is $\nu_1$ which vanishes on the Legendrian, so Proposition~\ref{prop:symb-calc} tells us that the result is a Legendre distribution is of order $(m+1, r_1, r_2)$ and the principal symbol is given by \eqref{transport}.

The symbol of $u$ at $x_3 = 0$ is the half-density (where for convenience we write $u$ as a b-half-density on the Legendrian)
$$
\sigma^m(u) = x_2^{r_2 - m} x_1^{r_1 - m} \Big| \frac{dx_1 dx_2 dx_3 dy_1 dy_2 dy_3}{x_3 x_2 x_1} \Big|^{1/2} \otimes |d(x_1x_2x_3)|^{m-N/4}.
$$
The scattering Hamilton vector field of $P$ is $x_1 \partial_{x_1}$. The subprincipal symbol of $P$ is $-i(N-1-f_1)$, which is easily obtained from the fact that $P + P^*$ has vanishing subprincipal symbol. Finally $\partial p/\partial \nu_1 = 1$.  Thus, noting that $\mathcal{L}_{\scH_p}$ leaves the b-half-density $dx_1/x_1$ invariant,  \eqref{transport} says that
$$
\sigma^{m+1}(Pu) =  \Big( -i(r_1 - m) - i(\half + m - N/4) + \frac{-i}{2} (N - 1 - f_1) \Big) \sigma^m(u) =  -i(r_1 -f_1/2 + N/4) \sigma^m(u)
$$
in agreement with the direct calculation. 
\end{example}

\subsection{Residual space}\ilabel{resid} The residual space for the spaces of Legendre distributions\break $I^{m,r_1, r_2}(X, G; \sfOh)$ is, by definition, the intersection of these spaces over all $m \in \RR$, and is denoted $I^{\infty,r_1, r_2}(X, G; \sfOh)$. Let us consider the special case that $X = Y \times [0, \epsilon)$ as in Example~\ref{23}.
In that case, for a fixed $x_3 > 0$ an element of $I^{m,r_1, r_2}(X, G; \sfOh)$ is (after division by $|dx_3|^{1/2}$) a Legendre distribution on $Y$ belonging to $I^{r_1 - 1/4, r_2 - 1/4}(Y, x_3^{-1} G_2, \sfOh)$, in particular  associated to the Legendre submanifold $x_3^{-1} G_2$, where $G_2 = \partial_2 G$ is the boundary of $G$ over $H_2$ and the factor $x_3^{-1}$ scales the cotangent variables (this follows immediately from \eqref{eq:dist-2}). 
We may regard the spaces $I^{r_1 - 1/4, r_2 - 1/4}(Y, x_3^{-1} G_2, \sfOh)$ as forming a smooth bundle over $(0, \epsilon)_{x_3}$. 
The residual space $I^{\infty,r_1, r_2}(X, G; \sfOh)$ can then be described as a smooth, $O(x_3^\infty)$ section of this bundle on $[0, \epsilon)$. We write this (with a minor abuse of notation) as
$$
I^{\infty,r_1, r_2}(X, G; \sfOh) \equiv x_3^\infty C^\infty \big( [0, \epsilon]; I^{r_1 - 1/4, r_2 - 1/4}(Y, x_3^{-1}G_2; \sfOh) \big) \otimes |dx_3|^{1/2}.
$$
We remark that the rather irritating drop of $1/4$ in the orders, when regarding elements of $I^{\infty,r_1, r_2}(X, G; \sfOh)$ as distributions on $Y$ parametrized by $x_3$, follows from the order convention where a 
Legendre distribution is order $N/4$ if it is borderline $L^2$. In terms of \eqref{eq:dist-1} and \eqref{eq:dist-2} it can be seen since $f_i$ and $N$ both decrease by $1$ when we fix a value of $x_3 > 0$.

%%%%%%%%%%%%%%%%%%%%%%%%%%%%%%%%%%%

\section{Intersecting Legendre distributions}\ilabel{section:ild}
For a manifold with boundary, $M$, intersecting Legendre distributions 
were defined in \cite{hasvas1} as the analogue of the intersecting Lagrangian distributions of \cite{MU}. They are related to a pair of Legendre submanifolds in $\Tscstar[\partial M] M$ that intersect cleanly in codimension 1. Here we define the analogue for a scattering-fibred manifold with codimension two corners.

\subsection{Intersecting Legendre submanifolds}\ilabel{sec4.1}
Let $X$ be a scattering-fibred manifold with codimension two corners. 
By Proposition~\ref{X-coords}, locally near the corner, there are local coordinates $(x_1, x_2, y_1, y_2)$ with
respect to which the main face is given by $x_2 = 0$, the boundary
hypersurface $H_1$ is given by $x_1 = 0$ and the fibration at $H$ is given by
$(x_2, y_1, y_2) \mapsto y_1$. We define a pair of intersecting Legendre
submanifolds , $\Lt = (L, \Lambda)$, in $\Tsfstar[\mf] X$, to be a pair consisting of a Legendre submanifold $L$ in the sense of Definition~\ref{leg-sm-def}, thus a manifold with boundary meeting $\Tsfstar[\mf \cap H_1] X$ transversally, together with a submanifold $\Lambda$ with codimension two corners of $\Tsfstar[\mf] X$ which is Legendre, transversal to $\Tsfstar[\mf \cap H_1] X$, and satisfying the following:

\begin{itemize}
\item $\Lambda$ has two boundary hypersurfaces, $\partial_1 \Lambda = \Lambda \cap \Tsfstar[\mf \cap H_1] X$, and $\partial_L \Lambda = L \cap \Lambda$;

\item the intersection $L \cap \Lambda$ is clean;

\item the images $L_1 = \phit_1(\partial_1 L)$ and $\Lambda_1= \phit_1 (\partial_1 \Lambda)$ (which are Legendre in $\Nsf Z_1$ by Proposition~\ref{prop:fibre}) form an intersecting pair of Legendre submanifolds in $\Nsf Z_1$. 
\end{itemize}

\subsection{Parametrization}

A local parametrization of $(L, \Lambda)$ near $q \in L \cap \Lambda \cap
\Tsfstar _{\mf \cap H_1} X$ is a function of the form
\begin{multline}
\Phi(x_1, y_1, y_2, v_1, v_2, s) = \phi_{00}(y_1, v_1) + s \phi_{10}(y_1,
v_1, s) \\ + x_1\phi_{01}(x_1, y_1, y_2, v_1, v_2) + x_1 s \phi_{11}(x_1,
y_1, y_2, v_1, v_2, s),
\ilabel{intleg-form}\end{multline}
defined in a neighbourhood of $q' = (0,y_1^*,y_2^*, v_1^*, v_2^*,0)$ in
$\mf \times \RR^{k_1 + k_2} \times [0, \infty)$ such that $d_{v_1, v_2, s}
\Phi = 0$ at $q'$, $q = (0, y_1, y_2, d (\Phi/x_1x_2)(q'))$, $\Phi$
satisfies the non-degeneracy hypothesis
$$ ds, \ d\phi_{10}, \ d\Big( \frac{ \pa \phi_{00}}{\pa v_1^j} \Big), \
d\Big( \frac{ \pa \phi_{01}}{\pa v_2^k} \Big) \text{ are linearly independent
at } q',
$$ and near $q$,
$$ \begin{gathered} L = \{ \big(x_1, y_1, y_2, d \left( \frac{\Phi}{\xx}
\right) \big) \mid s=0, d_{v_1, v_2} (\Phi)
%_{00} + x_1 \phi_{01}) 
= 0 \}, \\ \Lambda = \{ \big(x_1, y_1, y_2, d \left( \frac{\Phi}{\xx} \right)
\big) \mid s\geq 0, d_s \Phi = 0, \ d_{v_1, v_2} \Phi = 0 \}.
\end{gathered}
$$

\subsection{Existence of parametrizations}

For simplicity we shall prove existence of para\-metrizations only in a special case, which nevertheless suffices for our application. We shall assume that $L$ is a `conormal bundle' of a submanifold $N \subset \mf$ that meets the boundary $x_1 = 0$ transversally. We shall further assume that the projection from $L \cap \Lambda$ to $N$ everywhere has maximal rank. We need only prove existence of a parametrization locally near a point $q \in L \cap \Lambda \cap
\Tsfstar _{\mf \cap H} X$ as above, since existence near other points has been shown in \cite{hasvas1} or in the previous section. 

By Proposition~\ref{prop:fibre}, the boundary of $N$ necessarily fibres over a submanifold $N_1 \subset Z_1$. Choose coordinates $y_1 = (y_1', y_1^\sharp)$ on $Z_1$ so that $N_1 = \{ y_1' = 0 \}$ locally. We can then find a splitting 
 $y_2 = (y_2^\flat, y_2^\sharp)$ with respect to which $N$ locally takes the form $\{ y_1' = 0, y_2^\flat = 0 \}$. Our assumption on $L$ reads as follows in local coordinates: $$L = \{ y_1' = 0, y_2^\flat = 0, \mu_1^\sharp = 0, \mu_2^\sharp = 0, \nu_1 = 0, \nu_2 = 0 \}.$$

Let us first parametrize the intersecting pair of Legendrians $(L_1, \Lambda_1)$. We first claim that one can split (after a suitable linear change of $y_1$ variables) $y_1'$ as $y_1' = (y_1^\flat, y_1^\natural)$, where $\dim y_1^\natural = 1$, in such a way that 
$(y_1^\natural, y_1^\sharp, \mu_1^\flat)$ form coordinates locally on $\Lambda_1$. In fact, we have local coordinates $(y_1^\sharp, \mu_1')$ on $L_1$. The second assumption above has the consequence that local coordinates on $L_1 \cap \Lambda_1$ are furnished by $y_1^\sharp$ and all but one of the $\mu_1'$ variables; after making a linear change of variables, we may split $\mu' = (\mu_1^\flat, \mu_1^\natural)$ dual to the splitting of the $y_1'$ variables so that $y_1^\sharp$ and $\mu_1^\flat$ are coordinates on $L_1 \cap \Lambda_1$. It then follows from the condition that $\Lambda_1$ is Legendre with respect to the contact structure $d\nu_1 + \mu_1 \cdot dy_1$ that $(y_1^\natural, y_1^\sharp, \mu_1^\flat)$ furnish coordinates on $\Lambda_1$ locally. Thus we can write the other variables $y_1^\flat, \mu_1^\natural, \mu_1^\sharp, \nu_1$, restricted to $\Lambda_1$, uniquely as smooth functions of $(y_1^\natural, y_1^\sharp, \mu_1^\flat)$. In particular we have
$$
\nu_1 = N_{1,1}(y_1^\natural, y_1^\sharp, \mu_1^\flat), \quad y_1^\flat = Y_{1,1}^\flat(y_1^\natural, y_1^\sharp, \mu_1^\flat),
$$
and each of these functions is $O(y_1^\natural)$ since they vanish at $L_1 \cap \Lambda_1$ which is  $\Lambda_1 \cap \{ y_1^\natural = 0 \}$. 

Then a local parametrization of $(L_1, \Lambda_1)$ is given by
$$
(y_1^\natural - s) v_1^\natural + (y_1^\flat - Y_{1,1}^\flat(s, y_1^\sharp, v_1^\flat)) \cdot v_1^\flat + N_{1,1}(s, y_1^\sharp, v_1^\flat);
$$
the reasoning is the same as in the proof of Proposition~\ref{prop:param}. 

We now parametrize $(L, \Lambda)$ in a neighbourhood of a point on $L \cap \Lambda \cap \{ x_1 = 0 \}$. In this case 
%there is some splitting $y_2 = (y_2^\flat, y_2^\sharp)$ of the $y_2$ coordinates such that  
$(x_1, y_1^\sharp, y_2^\sharp, \mu_1^\flat, \mu_1^\natural, \mu_2^\flat)$ furnish local coordinates on $L$ and $(x_1, y_1^\natural, y_1^\sharp, y_2^\sharp, \mu_1^\flat, \mu_2^\flat)$ furnish local coordinates on $\Lambda$. As before we write 
\begin{gather*}
\nu_1 = N_1(x_1, y_1^\natural, y_1^\sharp, y_2^\sharp, \mu_1^\flat, \mu_2^\flat), \quad y_1^\flat = Y_1^\flat(x_1, y_1^\natural, y_1^\sharp, y_2^\sharp, \mu_1^\flat, \mu_2^\sharp), \\
Y_2^\flat(x_1, y_1^\natural, y_1^\sharp, y_2^\sharp, \mu_1^\flat, \mu_2^\sharp).
\end{gather*}
Due to the conditions on $L$ and $\Lambda$ at $x_1 = 0$ we have
$N_1 = N_{1,1} + x_1N_{1,2}$, $Y_1^\flat = Y_{1,1}^\flat + x_1 Y_{1,2}^\flat$ and $Y_2^\flat =  x_1 Y_{2,2}^\flat$ for some smooth functions $N_{1,2}, Y_{1,2}^\flat$ and $Y_{2,2}^\flat$. Then the function
\begin{gather*}
(y_1^\natural - s) v_1^\natural + (y_1^\flat - Y_{1}^\flat(x_1,  y_1^\sharp, y_2^\sharp, v_1^\flat, s, v_2^\flat)) \cdot v_1^\flat + N_{1}(x_1, y_1^\sharp, y_2^\sharp, v_1^\flat, s, v_2^\flat) \\
= (y_1^\natural - s) v_1^\natural + (y_1^\flat - Y_{1,1}^\flat(s, y_1^\sharp, v_1^\flat)) \cdot v_1^\flat + N_{1,1}(s, y_1^\sharp, v_1^\flat) + O(x_1) \\ 
= y_1^\natural v_1^\natural + y_1^\flat \cdot v_1^\flat + y_2^\flat \cdot v_2^\flat + O(s)
\end{gather*}
has the form \eqref{intleg-form} and parametrizes $(L, \Lambda)$.

\subsection{Equivalence of phase functions}

Two phase functions $\Phi$, $\tilde{\Phi}$ are said to be \emph{equivalent}
if they have the same number of phase variables of each type $v_1,v_2$ and
there exist maps
$$ V_1(x_1,\vecy,\vecv,s),\ V_2(x_1,\vecy,\vecv,s),\ S(x_1,\vecy,\vecv,s)
$$ such that
$$ \tilde{\Phi}(\vecx,\vecy,V_1,V_2,S) = \Phi.
$$

\begin{proposition}
The phase functions $\Phi=\phi_{00}+s \phi_{10} +
x_1\phi_{01}+x_1s\phi_{11}$ and $\tilde\Phi=\tilde\phi_{00}+s
\tilde\phi_{10} + x_1\tilde\phi_{01}+x_1s\tilde\phi_{11}$ are locally
equivalent iff
\begin{enumerate}
\item
They parametrize the same Legendrians,
\item
They have the same number of phase variables of the form $v_1,$ $v_2$
separately,
\item
\begin{align*}
\sgn d_{v_1}^2 (\phi_{00}+s\phi_{10}) &= \sgn d_{v_1}^2
(\tilde\phi_{00}+s\tilde\phi_{10}),\\ \sgn d_{v_2}^2 (\phi_{01}+s\phi_{11})
&= \sgn d_{v_2}^2 (\tilde\phi_{01}+s\tilde\phi_{11}),
\end{align*}
\end{enumerate}
\end{proposition}
\begin{proof}
Using the equivalence of phase functions in the codimension one case from
\cite{hasvas1} to solve the problem at $x_1=0,$ and using
Proposition~\ref{prop:equiv1} to solve at $s=0,$ we may assume that we have
reduced to the case
$$ \tilde\Phi=\phi_{00}+s \phi_{10} + x_1 \phi_{01}+x_1s\tilde\phi_{11}.
$$ As before, we may further reduce by an initial change of variables to
the case in which we may take $\Phi,\tilde\Phi$ equal to second order along
$C = \{d_s\Phi = d_{\vecv}\Phi=0\}.$

As the two functions agree to second order on $C,$ we may expand in a
Taylor series
$$ \tilde\phi_{11} - \phi_{11} = \frac 12 (\nabla'_{\vecv,s}\Phi)^t B
(\nabla'_{\vecv,s}\Phi)
$$ where we define $\nabla'_{\vecv, s} \Phi = (\pa_{v_1} \Phi,
\pa_{v_2}(\phi_{01}+ s\phi_{11}), \pa_s\Phi).$ We further expand
$$ \Phi(x_1,\vecy,\tilde{\vecv},\tilde s) - \Phi(x_1,\vecy,\tilde{\vecv},s)
= (\tilde{\vecv}-\vecv ) \cdot \pa_{\vecv}\Phi+(\tilde{s}-s ) \cdot
\pa_{s}\Phi + O((\tilde{\vecv}-\vecv )^2+(\tilde s-s)^2).
$$ Set
$$ (\tilde v_1, \tilde v_2,\tilde s) -(v_1,v_2, s) = (x_1 w_1, w_2, x_1
w_3)\cdot \nabla'_{\vecv}\Phi
$$ for $w_i=w_i(x_1,\vecy,\vecv,s).$ We thus have
$$ \Phi(x_1,\vecy,\tilde{\vecv},\tilde{s}) - \Phi(x_1,\vecy,\vecv,s) = x_1
(\nabla'_{\vecv,s}\Phi)^t (w+O(w^2)) (\nabla'_{\vecv,s}\Phi).
$$ We want
\begin{align*} \Phi(x_1,\vecy,\tilde{\vecv},\tilde{s}) - \Phi(x_1,\vecy,{\vecv},s) &=
\tilde\Phi(x_1,\vecy,\vecv,s) -\Phi(x_1,\vecy,s)\\ &=x_1 s
(\tilde\phi_{11}(x_1,\vecy,\vecv,s)-\phi_{11}(x_1,\vecy,\vecv,s))
\end{align*}
We thus need to solve
$$ x_1 (\nabla'\Phi)^t (w+O(w^2)) (\nabla'\Phi)= \frac{x_1 s}2
(\nabla'\Phi)^t B (\nabla'\Phi)
$$ for $w.$ This can always be accomplished for $s$ small by the inverse
function theorem.
\end{proof}

\subsection{Intersecting Legendre distributions}\label{ilds}

Let $\nu$ be a smooth scattering-fibred half-density. The set of Legendre distributions of order $(m, r)$ associated to $\Lt$, denoted\break $I^{m,r}(X, \Lt; \sfOh)$,  is the set of 
half-density distributions of the form $u = u_1 + u_2 + u_3 + (u_4 + u_5 + u_6) \nu$, where  
\begin{itemize}

\item 
$u_1 \in I^{m,r}(X, \Lambda; \sfOh)$ with the microsupport of $u_1$ disjoint from $\partial \Lambda$, 

\item $u_2 \in I^{m+1/2,r+1/2}(X, L; \sfOh)$,  

\item $u_3$ has support disjoint from $H_1$ and is an intersecting Legendre distribution of order $(m,r)$ associated to $(L, \Lambda)$ as defined in \cite{hasvas1}, 

\item $u_4$ is a finite sum of terms, each 
supported near $\mf = \{ x_2 = 0 \}$, with an
expression
\begin{equation}\begin{gathered}
 x_1^{j_1} x_2^{j_2} \int_0^\infty  \int \int
e^{i\Phi(x_1,y_1,y_2,v_1,v_2,s)/x_1x_2} a(x_1,x_2,y_1,y_2,v_1,v_2,s) \,
dv_1 \, dv_2  \, ds, \\ j_1 = r-\frac{k_1 + 1}{2} + \frac{N}{4} -\frac{f}{2}, \ j_2 =
m- \frac{k_1 + k_2+1}{2} +\frac{N}{4}
\end{gathered}\ilabel{ild}\end{equation}
where $v_i \in \RR^{k_i}$, $a$ is smooth and compactly supported, $f$ is the dimension of the fibres of $H_1$, and $\Phi = \phi_{00} + s \phi_{01} + x_1 \phi_{10} + x_1 s \phi_{11}$ locally parametrizes $(L, \Lambda)$ near a point $q \in L \cap \Lambda \cap \Tsfstar[\mf \cap H_1] X$, as in \eqref{intleg-form},

\item $u_5$ is a finite sum of terms of the form
\begin{equation}
 x_1^{ r-\frac{k_1 + 1}{2} + \frac{N}{4} -\frac{f}{2}} \int_0^\infty  \int 
e^{i(\phi_{00} + s \phi_{01})/x_1x_2} b(x_1,y_1,x_2, y_2,v_1,s) \,
dv_1   \, ds, 
\end{equation}
where $\phi_{00}, \phi_{01}$, $f$ and $v_i$ are as above, and $b$ is smooth and $O(x_2^\infty)$ at $\mf$, and

\item $u_6 \in \CIdot(X)$.
\end{itemize}

 As in Section~\ref{sec:legendrian}, $u_3$ can be written with respect to any local
parametrization, up to an error in $\CIdot(X).$   This follows from the equivalence result above and the argument in Proposition~\ref{prop:any}. 
%The set of such half-densities is denoted $I^{m,r}(X, \Lt; \sfOh)$.

\subsection{Symbol calculus}\label{ildsc}
The geometry of intersecting Legendre distributions is such that the symbol on $L$ has a $1/\rho_1$ singularity at $\Lambda$, where $\rho_1$ is a boundary defining function for $\partial \Lambda \subset L$, while the symbol on $\Lambda$ is smooth up to the boundary at $L \cap \Lambda$. This allows one to symbolically solve away error terms at $L$ in the equation $Pu = f$ where $f$ is Legendrian on $L$, and the principal symbol of $P$ vanishes simply at $\Lambda$; what happens is that the singularities of the solution $u$  propagate from $L \cap \Lambda$ along $\Lambda$. 
The formal symbol calculus for intersecting Legendre distributions on $X$ follows
readily from the codimension one case; we follow the description from \cite{hasvas2} closely. 

Let $\Lt = (L, \Lambda)$ be a pair of intersecting Legendre
submanifolds as in Section~\ref{sec4.1}. We consider  $u \in I^{m, r}(X, \Lt; \sfOh)$. 
The symbol of $u$ takes values in a bundle over $L \cup \Lambda$. 
To define this bundle, let
$\rho_1$ be a boundary defining function for $\pa \Lambda$ as a submanifold of
$L$, and $\rho_0$ be a boundary defining function for $\pa \Lambda$ as a
submanifold of $\Lambda$. Note that the symbol on $L$
is defined by continuity from distributions in $I^{m+1/2, r+1/2}(X, L; \sfOh)$ microsupported away from $\Lambda$,
and takes values in
\begin{equation}
x_1^{r-m} \rho_1^{-1} \CI(\Omega_b^{1/2}(L)\otimes S^{[m+1/2]}(L)) = x_1^{r-m} \rho_1^{-1/2}
\CI(\Omega_b^{1/2}(L \setminus \pa \Lambda) \otimes S^{[m+1/2]}(L)),
\ilabel{int-Leg-dist}\end{equation} 
while the symbol on $\Lambda$, defined by
continuity from distributions in $I^{m,r}(X, \Lambda; \sfOh)$  microsupported away from $\pa \Lambda$, takes values in
$$ 
x_1^{r-m} \rho_0^{1/2}\CI(\Omega_b^{1/2}(\Lambda) \otimes S^{[m]}(\Lambda)).
$$ 
Melrose and Uhlmann showed that the Maslov factors were canonically isomorphic on $L \cap \Lambda$, so $S^{[m+1/2]}(L)$ is naturally isomorphic to
$S^{[m]}(\Lambda) \otimes |N^*_{\mf} \partial X|^{1/2}$ over $L \cap \Lambda$. Canonical restriction of the
half-density factors to $L \cap \Lambda$ gives terms in $\CI(\Omega^\half(L
\cap \Lambda) \otimes S^{[m]}(\Lambda) \otimes |N^*_{L}\pa \Lambda|^{-1/2} \otimes
|N^* \pa X|^{1/2}$ and $\CI(\Omega^\half(L \cap \Lambda) \otimes S^{[m]}(\Lambda)
\otimes |N^*_{\Lambda}\pa \Lambda|^{1/2}$ respectively. In fact $|N^*_{L}\pa \Lambda|
\otimes |N^*_{\Lambda}\pa \Lambda| \otimes |N_{\mf}^* \pa X|^{-1}$ is canonically
trivial; an explicit trivialization is given by
\begin{equation}
(d\rho_0, d\rho_1, (x_1x_2)^{-1}) \mapsto (x_1x_2)^{-1} \omega(V_{\rho_0}, V_{\rho_1})
\restriction L \cap \Lambda, \ilabel{triv}\end{equation} where $V_{\rho_i}$
are the Hamilton vector fields of the functions $\rho_i$, and $\omega$ is the standard symplectic form. Thus the two
bundles are naturally isomorphic over the intersection.

We define the bundle
$S^{[m]}(\Lt)$ over $\Lt = L \cup \Lambda$ to be that bundle such that smooth sections of
$\Omega^{1/2}_b(\Lt) \otimes S^{[m]}(\Lt)$ are precisely those pairs
$(a,b)$ of sections of $\rho_1^{-1} \CI(\Omega^{1/2}(L)\otimes
S^{[m+1/2]}(L))$ and $\rho_0^{1/2}\CI(\Omega^{1/2}_b(\Lambda) \otimes
S^{[m]}(\Lambda))$ such that\jw{Powers of $\rho_1$ don't match.}
\begin{equation}
\rho_1^{1/2} b = e^{i\pi/4} (2\pi)^{1/4} \rho_0^{-1/2} a \text{ at }L
\cap \Lambda \ilabel{compat}\end{equation} under the above  identification of
bundles (cf. equation (3.7) of \cite{hasvas2}). The symbol maps of order $m$
on $\Lambda$ and $m+1/2$ on $L$ then extend in a natural way to a symbol map
of order $m$ on $\Lt$ taking values in $\Omega_b^{1/2}(\Lt) \otimes
S^{[m]}(\Lt)$.

\begin{proposition} The symbol map on $\Lt$ yields an
exact sequence
\begin{equation}
%\begin{CD}
0 \to I^{m+1, r}(X, \Lt; \sfOh) \to I^{m, r}(X, \Lt; \sfOh) \to
x_1^{r-m}\CI(\Lt, \Omega^\half_b \otimes S^{[m]}) \to 0.
%\end{CD}
\ilabel{ex-int-1}\end{equation} Moreover, if we consider just the symbol
map to $\Lambda$, there is an exact sequence
\begin{multline}
0 \to I^{m+1, r}(X, \Lt; \sfOh) + I^{m+\half, r}(X, L; \sfOh)
\to I^{m, r}(X, \Lt; \sfOh) \\ \to x_1^{r-m}\CI(\Lambda,
\Omega^\half \otimes S^{[m]}) \to 0.  \ilabel{ex-int-2}\end{multline} 

If $P
\in \Diffsf(X; \sfOh)$ has principal symbol $p$ and $u \in I^{m, r}(X, \Lt; \sfOh)$, then $Pu
\in I^{m, r}(X, \Lt; \sfOh)$ and
$$ \sigma^m(Pu) = \big( p \rest \Lt \big) \sigma^m(u).
$$ Thus, if $p$ vanishes on $\Lambda$, then $Pu$ is an element of
$I^{m+1, r}(X, \Lt; \sfOh)$ $+ I^{m, r}(X, L; \sfOh)$ by \eqref{ex-int-2}. The
symbol of order $m+1$ of $Pu$ on $\Lambda$ in this case is given by
\eqref{transport}.
\end{proposition}

\subsection{Residual space}\ilabel{int-resid} The residual space for the spaces of intersecting Legendre distributions $I^{m,r}(X, \Lt; \sfOh)$ is
$$
I^{\infty, r}(X, \Lt; \sfOh) = \cap_m I^{m,r}(X, \Lt; \sfOh).
$$
If $X = Y \times [0, \epsilon]_{x_2}$ where $Y$ is a manifold with boundary, then the residual space may be identified with
$$
x_2^\infty C^\infty \big( [0, \epsilon]; I^{r - 1/4}(X, (x_2^{-1}L_1, x_2^{-1}\Lambda_1); \sfOh) \big) \otimes |dx_2|^{1/2}.
$$

%%%%%%%%%%%%%%%%%%%%%%%%%%%%%%%%%%

\section{Legendrian distributions with conic points}\ilabel{Lconicpts}

Here we shall define a more singular situation in which the Legendrian $G
\subset \Tsfstar[\mf] X$ has conic singularities. We first give a precise
description of `having conic singularities'. We recall the notion of real
blowup. Suppose that $X$ is a compact manifold with corners and $S \subset
X$ a compact product-type submanifold\footnote{All the submanifolds
considered in this paper are product-type submanifolds; from here on we
refer to them simply as submanifolds for brevity.}, which means that locally
near any point $s$ of $S$, there are local coordinates $x_1, \dots, x_j, y
= (y_1, \dots, y_k)$, $x_i \in [0, \epsilon)$, $y \in B(0, \epsilon)
\subset \RR^k$, with $s$ corresponding to the origin of coordinates, such
that $S$ is given locally by the vanishing of some subset of these
coordinates. Then by $[X;S]$ we denote the blow-up of $X$ around $S$. As a
set this is the union of $X \setminus S$ with the inward pointing spherical
normal bundle at $S$, which we denote $\tilde S.$  $[X;S]$ carries a natural
differentiable structure making it a compact manifold with corners, such
that $\tilde S$ is one of its boundary hypersurfaces.

\begin{defn} Let $X$ be a manifold with corners and $S \subset X$ a
  submanifold, and $G \subset X$ a closed set which is a submanifold locally near
  every point of $G \setminus S$. We say that $G$ has conic singularities
  at $S$ if the lift of $G$ to $[X;S]$, i.e.\ the closure of $G \setminus S$ in $[X; S]$, is  a smooth product-type submanifold
  $\hat G$ which is transverse to $\tilde S$.
\end{defn}

%\subsection{Review} 
Legendre submanifolds with conic singularities have
  been defined already in two different settings in \cite{MZ} and
\cite{hasvas1}, and we review these definitions for the convenience of the
reader.

The original setting of Melrose-Zworski was that of a Legendre
submanifold $G \subset \Tscstar[\partial X] X$ in the boundary of the
scattering cotangent bundle of a manifold $X$ with boundary, which has
conic singularities at a submanifold $\Jsharp$ which is the span of a
smooth projectable Legendrian $\Gsharp$. Projectability means that the
restriction of the projection $\pi: \Tscstar[\partial X] X \to \partial X$
to $\Gsharp$ is a diffeomorphism, or in other words $\Gsharp$ is a graph
over $\partial X$; then $\Jsharp$, which is obtained by replacing each
point of $\Gsharp$ by the ray in $\Tscstar[\partial X] X$ through this
point, is a submanifold with dimension equal to $\dim X$ (one greater than
$\dim \Gsharp$) . By choosing coordinates judiciously we may arrange that,
in local coordinates $(y, \nu, \mu)$ on $\Tscstar[\partial X] X$ given by
writing scattering covectors as
$$ \nu d\big( \frac1{x} \big) + \mu \cdot \frac{dy}{x},
$$ we have $\Gsharp = \{ \nu = 1, \mu = 0 \}$, and $\Jsharp = \{ \mu = 0
\}$. We say that $(G, \Gsharp)$ are a pair of Legendre submanifolds with
conic points, or a Legendrian conic pair for short, if $G$ has conic
singularities at $\Jsharp$, i.e.\ $G$ lifts to $[\Tscstar[\partial X] X;
\Jsharp]$ to a smooth submanifold $\Ghat$ transverse to $\tilde \Jsharp$.

We recall what it means to locally parametrize $(G,
\Gsharp)$. Transversality of $\Ghat$ to the span of $G^\sharp$ at $q \in
\Ghat \cap \tilde \Jsharp$ means that $d|\mu| \neq 0$ at $q$; we may assume
(after making a linear change of coordinates in $y$) that the first
component $\mu^1$ of $\mu$ is a local boundary defining function for the
blowup of the span of $G^\sharp$ near $q$.  Assuming this, a local
parametrization of $(G, G^\sharp)$ near $q$ is given by a phase function of
the form
$$ 1 + s \psi(y, s, v), \quad s \geq 0, v \in \RR^k
$$ defined in a neighbourhood of $(y^*, 0, v^*)$, satisfying the
non-degeneracy condition
\begin{equation}
d_{y^1} \psi \text{ and } d_{y, v} \big( \frac{ \partial \psi}{\partial
v^i} \big) \text{ are linearly independent at } (y^*, 0, v^*),
\end{equation}
such that $\Ghat$ is given by
\begin{equation}
\Ghat = \Big\{ d \big( \frac{ 1 + s \psi(y, s, v)}{x} \big) \mid d_{s,v}
\psi = 0 \Big\}.
\ilabel{ghat1}\end{equation}
 Furthermore we require that $d_v \psi(y^*, 0, v^*) = 0$ and that the point
on $\Ghat$ corresponding to $(y^*, 0, v^*)$ is $q$.  To be precise, the
meaning of \eqref{ghat1} is that when the set on the RHS is \emph{lifted}
to the space $[\Tscstar[\partial X] X; \Jsharp]$ obtained by blowup of
$J^\sharp$ it coincides with $\Ghat$. We remark that the correspondence in
\eqref{ghat1} lifts to a diffeomorphism from $\{ (y, s, v) \mid d_{s,v}
\psi = 0 \}$ to $\Ghat$, so the blowup is implicit in the parametrization
$\psi$.

\

Next we recall the definition of Legendre conic pairs in the case of a
manifold $X$ with fibred boundary and codimension 2 corners. Let
$\Gsharp$ be a smooth projectable
%\footnote{Here, \emph{projectable} means that that the projection  $\Tsfstar[\mf] X \to \mf$ restricts to $\Gsharp$ to be a diffeomorphism, or equivalently that $\Gsharp$ is given globally by $d(f/\xx)$ where $f \in \CIsf(X)$}   
Legendrian submanifold of
$\Tsfstar[\mf] X$, and $G$ be a Legendrian submanifold of
$\Tsfstar[\mf] X$ which is smooth away from $\Gsharp$ and which has
conic singularities at $\Jsharp \subset \Tsfstar[\mf] X$, where $\Jsharp$ is the span of $\Gsharp$. Let $\hat G$ denote the lift of $G$ to $[\Tsfstar_{\mf} X; \Jsharp]$; we assume that it is transverse to both boundary hypersurfaces of $[\Tsfstar_{\mf} X; \Jsharp]$ (that is, transverse to both the lift of $\Tsfstar_{\mf \cap H_1} X$ and the lift of $\Jsharp$). Let $\partial_1 \hat G$ and $\partial_1 \Gsharp$ denote the boundary hypersurface of $\hat G$, resp. $\Gsharp$, at (the lift of) $\Tsfstar_{\mf \cap H_1} X$. We say that
$(G, G^\sharp)$ form a conic Legendrian pair if  $\partial_1 \Ghat$ and $\partial_1 \Gsharp$ fibre over the  same Legendrian submanifold $G_1 \subset \Nsfstar
Z_1$ as base.

\begin{remark} This implies that the fibres of $\partial_1 \Ghat \to G_1$ and the fibres of $\partial_1 \Gsharp \to G_1$ form an intersecting pair of Legendre
submanifolds in $\Tscstar[\partial F] F$ for each fibre $F \subset H_1$. The reasoning is analogous to that in Proposition~\ref{prop:fibre}. 
\end{remark}

This differs from the structure above only over the codimension two corner
of $X$, so we shall consider a point of $\Tsfstar[H_1 \cap \mf] X$ lying
over the codimension two corner.  We shall use coordinates $(x_1, x_2, y_1,
y_2)$ as in Section~\ref{leg}, and associated dual coordinates $(\nu_1,
\nu_2, \mu_1, \mu_2)$ defined by writing scattering-fibred covectors in the
form
$$ \nu_1 d \big( \frac1{x_1 x_2} \big) + \nu_2 d \big( \frac1{x_2} \big) +
\mu_1 \cdot \frac{dy_1}{x_1 x_2} + \mu_2 \cdot \frac{dy_2}{x_2}.
$$ For definiteness we shall assume that $\Gsharp_2$ is the submanifold $\{
\nu_1 = 1, \nu_2 = 1, \mu_1 = 0, \mu_2 = 0 \}$ which is parametrized by the
function $1 + x_1$.  This is the form of $G_2^\sharp$ that turns up in our
application (and in any case, it can always be arranged by a change of
coordinates). Then the span of $G^\sharp_2$ is given by
\begin{equation}
J_2^\sharp = \{ x_2 = 0, \nu_1 = \nu_2, \mu_1 = 0, \mu_2 = 0 \}.
\end{equation}
The corresponding Legendrian in $\Nsfstar Z_1$ is

$$ \Gsharp_1 = \{ \nu_1 = 1, \mu_1 = 0 \}.
$$

The condition of being a conic Legendrian pair means that at $\{ x_1 = x_2
= 0\}$, if we set $\nu_1 = 1, \mu_1 = 0$ and fix $y_1$, then we have
remaining coordinates $(y_2, \nu_2, \mu_2)$ and these are local coordinates
on the fibre $\Tscstar[\partial F] F$ which is a contact manifold with
contact form $d\nu_2 + \mu_2 \cdot dy_2$; we are then asking that the
restriction of $G_2$ to this fibre have a conic singularity at (and
therefore becomes smooth after blowup of) $\{ \mu_2 = 0 \}$. In particular
$d|\mu_2| \neq 0$ on $\Ghat_2$ at its intersection with $\Jsharp_2$.

We next recall the form of a parametrization of $(G_2, \Gsharp_2)$ near a
point $q \in \Ghat_2$ on the codimension two corner of $\Ghat_2$, i.e.\
lying above $x_1 = 0$ and on $\tilde J_2^\sharp$.  Assume that coordinates
have been chosen so that $dy_2^1 \neq 0$ at $q$.  A local parametrization
of $(G_2, G_2^\sharp)$ near $q$ is given by a phase function of the form
$$ 1 + x_1 + sx_1 \psi(x_1, y_1, y_2, s, v), \quad s \geq 0, v \in \RR^k
$$ defined in a neighbourhood of $(0, y_1^*, y_2^*, 0, v^*)$, satisfying
the non-degeneracy condition
\begin{equation}
d_{y_2^1} \psi \text{ and } d_{y_2, v} \big( \frac{ \partial \psi}{\partial
v^i} \big) \text{ are linearly independent at } (0, y_1^*, y_2^*, 0, v^*),
\end{equation}
such that $\Ghat_2$ is given by
\begin{equation}
\Ghat_2 = \Big\{ d \big( \frac{ 1 + x_1 + sx_1 \psi(x_1, y_1, y_2, s,
v)}{x_1x_2} \big) \mid d_{s,v} \psi = 0 \Big\}.
\ilabel{ghat2}\end{equation}
 Furthermore we require that $d_v \psi(0, y_1^*, y_2^*, 0, v^*) = 0$ and
that the point on $\Ghat_2$ corresponding to $(0, y_1^*, y_2^* ,0, v^*)$ is
$q$.  The precise meaning of \eqref{ghat2} is that when the set in
\eqref{ghat2} is lifted to the space obtained by blowup of $J^\sharp_2$ it
coincides with $\Ghat_2$.

\begin{remark} As in the case above, the correspondence in \eqref{ghat2}
  lifts to a diffeomorphism from $\{ (x_1, y_1, y_2, s, v) \mid d_{s,v}
  \psi = 0 \}$ to $\Ghat$, so the blowup is implicit in the parametrization
  $\psi$. Also, if we fix a value of $\overline y_1$, or equivalently fix a
  point in the base $\Gsharp_1$ of the fibration $\phit_{12} |_G$, then the
  function $\psi(0, \overline y_1, y_2, s, v)$ parametrizes the fibre
  (which is a  Legendrian conic pair in $\Tscstar[\partial F] F$).
\end{remark}

\subsection{Legendre submanifolds with conic points}
We now define Legendre submanifolds with conic points in two new
situations, although both are closely analogous to the ones reviewed above.

\subsubsection{Codimension two corners}

Suppose that $X$ is a scattering-fibred manifold with
corners of codimension 2. Let $x_2$ be a boundary defining function for the
main face $\mf$ and $x_1$ a boundary defining function for the fibred face
$H_1$. Let $G_1^\sharp$ be a projectable Legendrian in $\Nsf Z_1$, and let $J$ be the lift of the span of $G_1^\sharp$ to $\Tsfstar[H_1 \cap \mf] X$ via the fibration
$\phit_{12}$. Let $\Ghat$ be the lift of $G$ to $[\Tsfstar[\mf] X, J]$. We shall say that $(G, G_1^\sharp)$ form a conic Legendrian pair of submanifolds if  $\Ghat$ has conic singularities at $J$, i.e.\ is transverse to both boundary hypersurfaces of $[\Tsfstar[\mf] X, J]$ (that is, transverse to both the lift of $\Tsfstar[\mf \cap H_1] X$ and to the lift $\tilde J$ of $J$). 

Let $\partial_1 \Ghat$ and $\partial_\sharp \Ghat$ denote the boundary hypersurfaces of $\Ghat$. Also, let $G_1$ denote the projection of $G \cap \{ x_1 = 0 \}$ to $\Nsf Z_1$ via
$\phit_{12}$. It follows from the definition that $G_1$ has conic singularities at $G_1^\sharp$; let $\hat G_1$ be the lift of $G$ to $[\Nsf Z_1; J_1]$ where $J_1$ is the span of $G_1^\sharp$. Then, as a consequence of $(G, G_1^\sharp)$ being a conic Legendrian pair,  the fibres of the map $\partial_1 \Ghat \to \hat G_1$ are Legendrian, while $\partial_\sharp \Ghat$ is itself Legendrian with respect to a natural contact structure on the lift of $G_1^\sharp$ to $\tilde J$ defined by the leading part of $\chi$.  \ah{statements in this para need to be checked}

\subsubsection{Codimension three corners}
Now let us assume that $X$ is a scattering-fibred manifold with
corners of codimension 3, and consider a Legendrian submanifold $G \subset
\Tsfstar[\mf] X$ which is singular at the boundary. We use the notation
$H_1, H_2, H_3$ for boundary hypersurfaces of $X$ and $x_1, x_2, x_3$ for
boundary defining functions as in Section~\ref{leg}.  Let $G_1 = \phi_{13}
(G \cap \{ x_1 = 0 \})$ and $G_2 = \phi_{23} (G \cap \{ x_2 = 0 \})$.  Here
we could consider the cases where either $G_1$ or $G_2$ have conic
singularities at some Legendrian $G_1^\sharp$ or $G_2^\sharp$; however, we
shall only consider the case where $G_2$ has conic singularities since that
is the case that occurs in our applications.  Thus, we consider a case
where $G_1$ is smooth, but $G_2$ has conic singularities, and indeed that
there is a projectable smooth Legendrian $\Gsharp_2 \subset \Nsf Z_2$ such
that $(G_2, G_2^\sharp)$ form a Legendrian conic pair. Thus, if $J_2$ is the span of $G_2^\sharp$ in $\Nsf Z_2$, then $G_2$ lifts to a smooth manifold $\hat G_2$ in $[\Nsf Z_2; J_2]$ that is transversal to $\tilde J_2$.  
Let $J$ denote the preimage of $J_2$ inside  $\Tsfstar[H_2 \cap \mf] X$ via $\phit_{23} : \Tsfstar[\mf \cap H_2] X \to \Nsfstar Z_2$. 
We shall say that $(G, G_2^\sharp)$ form a conic Legendrian
pair if $G$ has conic singularities at $J$, i.e.\ the lift $\Ghat$ of
$G$ to $[\Tsfstar[\mf] X; J]$ is smooth and transverse to $\tilde J$ as
well as to the lifts of $\Tsfstar[\mf \cap H_1] X$ and $\Tsfstar[\mf \cap
H_2] X$.

The manifold $\Ghat$ is a manifold with corners of codimension three. The boundary at $\Tsfstar[\mf \cap H_1] X$ (more precisely, at the lift of this to $[\Tsfstar[\mf] X; J]$) is denoted $\partial_1 \Ghat$, the boundary at the lift of $\Tsfstar[\mf \cap H_2] X$ is denoted $\partial_2 \Ghat$ and the boundary at $\tilde J$ is denoted $\partial_\sharp \Ghat$. 
It follows from the definition that $\partial_1 \Ghat$ fibres over $G_1$ with Legendrian fibres relative to $\chi_{13}$, that $\partial_2 \Ghat$ fibres over $\hat G_2$ via a map $\phi_{23}^G$ induced from $\phit_{23}$, with fibres that are Legendrian for the contact structure $\chi_{23}$, and  $\partial_\sharp \Ghat$ is Legendrian for the contact structure on
the lift of $G_2^\sharp$ to $\tilde J$ given by the leading part of $\chi$. \ah{ditto}

\subsection{Parametrization}

\subsubsection{Codimension two corners}\ilabel{conic-param-2}
In this situation, the lifted submanifold $\Ghat$ is a manifold with
corners of codimension two. The two boundary hypersurfaces of $\Ghat$ are denoted $\partial_1 \Ghat$ (at $x_1 = 0$ and
away from $\tilde J$) and  $\partial_\sharp \Ghat$ (at $\{ x_1 = 0 \} \cap \tilde J$). Locally near a point on the interior of $\partial_1 \Ghat$ the situation is as for a smooth Legendrian distribution, so consider a
point $q$ on $\partial_\sharp \Ghat$. We need to
distinguish two cases: the first is that $q$ is on the codimension two
corner $\partial_1 \Ghat \cap \partial_\sharp \Ghat$, and the second is that $q$ is on the interior of $\partial_\sharp \Ghat$. 

To make things concrete we shall assume that coordinates have been chosen
so that $\Gsharp_1$ is the Legendrian $\{ \nu_1 = 1, \mu_1 = 0 \}$, and
that $\mu_1^1$ is a local boundary defining function for $\tilde J$.  Then
a local parametrization of $(G, G_1^\sharp)$ near $q$ is a phase function
of the form
\begin{equation}
\psi(s, x_1, y_1, y_2, v_1, v_2) = 1 + s \psi_1(y_1, s, v_1) + x_1
\psi_2(s, \frac{x_1}{s}, y_1, y_2, v_1, v_2), \quad s\geq 0, v_i \in
\RR^{k_i}, \ilabel{psi1psi2}\end{equation} defined in a neighbourhood of
$(0, 0, y_1^*, y_2^*, v_1^*, v_2^*)$, where $\psi_1$ and $\psi_2$ are
smooth, satisfying the non-degeneracy condition
\begin{equation}
d_{y_1^1} \psi_1, d_{y_1, v_1} \big( \frac{ \partial \psi_1}{\partial
v_1^i} \big) \text{ and } d_{y_2, v_2} \big( \frac{ \partial
\psi_2}{\partial v_2^j} \big) \text{ are linearly independent at } (0, 0,
y_1^*, y_2^*, v_1^*, v_2^*), \ilabel{psi-nd}\end{equation}
and such that $\Ghat$ is given by
\begin{equation}
\Big\{ d \big( \frac{\psi}{x_1 x_2} \big) \mid d_{s,v_1, v_2} \psi = 0 \Big\}.
\ilabel{ghat-codim2}\end{equation}
Furthermore we require that $d_{s,v_1, v_2} \psi(0, 0, y_1^*, y_2^*, v_1^*,
v_2^*) = 0$, and that the point on $\Ghat$ corresponding to $(0, 0, y_1^*,
y_2^*, v_1^*, v_2^*)$ is $q$.

\begin{remark} The non-degeneracy conditions imply that the subset
$$ C_\psi = \{ (s, u, y_1, y_2, v_1, v_2) \mid d_s \psi = 0, d_{v_1}(\psi_1
+ u \psi_2) = 0, d_{v_2} \psi_2 = 0,\ u = \frac{x_1}{s} \}
$$ is a submanifold\footnote{The partial derivative $d_s \psi$ in the
equation above is taken keeping $x_1$ fixed, not keeping $u$ fixed} and
that \eqref{ghat-codim2} defines a diffeomorphism between $C_\psi$ and
$\Ghat$ locally near $(0, 0, y_1^*,
y_2^*, v_1^*, v_2^*)$, so this indeed corresponds to the usual notion of non-degenerate
paramet\-rization. Notice that under this correspondence $s$ is a boundary
defining function for $\partial_\sharp \Ghat$  and $u$ is a boundary defining
function for $\partial_1 \Ghat$.
\end{remark}

In the second case, since we are away from the lift of $\{ x_1 = 0 \}$,
given by $x_1/s = 0$, we do not need the special variable $s \geq 0$, and
we obtain the following: a local parametrization of $(G, G_1^\sharp)$ near
$q$ is a phase function of the form
\begin{equation}
\psi(x_1, y_1, y_2, v) = 1 + x_1 \psi(x_1, y_1, y_2, v)
\ilabel{psi-simple}\end{equation}
defined in a neighbourhood of $(0, y_1^*, y_2^*, v^*)$, satisfying the
non-degeneracy condition
\begin{equation}
 d_{y_1, y_2, v} \big( \frac{ \partial \psi}{\partial v^j} \big) \text{ are
linearly independent at } (0, y_1^*, y_2^*, v^*),
\ilabel{psi-simple-nd}\end{equation} such that $\Ghat$ is given by
\begin{equation}
\Big\{ d \big( \frac{ \psi}{x_1 x_2} \big) \mid
d_{v} \psi = 0 \Big\}.
\ilabel{ghat-codim2-2}\end{equation}
Furthermore we require that $d_{v} \psi(0, y_1^*, y_2^*, v^*) = 0$, and
that the point on $\Ghat$ corresponding to $(0, y_1^*, y_2^*, v^*)$ is $q$.

\begin{remark} This is very similar to the parametrization of a smooth Legendrian, but with respect to a different fibration on $H_1$,
where the base of the fibration is a point. This can also be seen by noting
that blowing up $\{ \mu_1 = 0, x_1 = 0 \}$ amounts to introducing the
variable $M_1 = \mu_1/x_1$ as a smooth coordinate. This is dual to
$dy_1/x_2$ and so corresponds to a coordinate along the fibre of the
fibration rather than on the base. This is related to the blowup of the submanifold $W$  in Section~\ref{sec:prop}.
\end{remark}

\begin{remark} Notice that, if we localize the phase function in \eqref{psi1psi2} to the region $x_1/s \geq \epsilon > 0$, then it can be expressed in the form 
$$ 1 + x_1 \big( w \psi_1(y_1, x_1 w, v_1) + \psi_2(x_1 w, 1/w, y_1, y_2,
v_1, v_2) \big), \quad w = \frac{s}{x_1},
$$ and is therefore of the form \eqref{psi-simple}. So these two forms of
parametrization are consistent on their overlapping regions of validity.
\end{remark}

\subsubsection{Codimension three corners}\ilabel{conic-param-3}

Now the lifted submanifold $\Ghat$ is a manifold with corners of
codimension three. The three boundary hypersurfaces are  denoted $\partial_1 \Ghat$ (at $x_1 = 0$), $\partial_2 \Ghat$ (at
$x_2 = 0$ and away from $\Gsharp_2$), and $\partial_\sharp \Ghat$ (at $\{ x_2 = 0 \}
\cap \Gsharp_2$). Locally near a point on the interior of $\partial_1 \Ghat$ or $\partial_2 \Ghat$ the situation is as for a smooth Legendrian
distribution, so consider a point $q$ on $\partial_\sharp \Ghat$. If $q$ is not also in $\partial_1 \Ghat$ then the situation is (locally) the codimension two situation described
above, so we assume that $q \in \partial_\sharp \Ghat \cap \partial_1 \Ghat$.  We need to distinguish two
cases: the first is that $q$ is on the codimension three corner $\partial_1 \Ghat \cap \partial_2 \Ghat \cap \partial_\sharp \Ghat$ and the second is that $q$ is on the interior of $\partial_\sharp \Ghat \cap \partial_1 \Ghat$.

To make things concrete we shall assume that coordinates have been chosen
so that $\Gsharp_1$ is the Legendrian $\{ \nu_1 = 1, \mu_1 = 0 \}$, that
$\Gsharp_2$ is the Legendrian $\{ \nu_1 = \nu_2 = 1, \mu_1 = 0, \mu_2 = 0
\}$, so that $1 + x_1$ parametrizes $\Gsharp_2$, and that $\mu_2^1$ is a
local boundary defining function for the third boundary hypersurface of
$\Ghat$.  Let $q \in \Ghat$ lie on the codimension three corner.  A
non-degenerate parametrization of $(G, \Gsharp_2)$ near $q \in \Ghat$ is
then a smooth phase function $\Psi$ of the form
\begin{equation}\begin{gathered}
\Psi(s, x_1, x_2, y_1, y_2, y_3, v_2, v_3) = 1 + x_1 + sx_1 \psi_2(s, x_1,
y_1, y_2, v_2) \\ + x_1 x_2 \psi_3(s, x_1, x_2/s, y_1, y_2, y_3, v_2, v_3)
, \quad s \geq 0, v_i \in \RR^{k_i},
\end{gathered}\ilabel{Psi1}\end{equation}
where $\psi_2$ and $\psi_3$ are smooth, with $\Psi$ non-degenerate in the
sense that such that
\begin{equation}
d_{y_2} \psi_2, \ d_{y_2, v_2} \big( \frac{\partial \psi_2}{\partial v_2^i}
\big) \text{ and } d_{y_3, v_3} \big( \frac{\partial \psi_3}{\partial
v_3^i} \big) \text{ are linearly independent at } q'
\ilabel{Psi-nd}\end{equation} with
\begin{equation} \Ghat = \{ d\big( \frac{\Psi}{x_1x_2x_3}\big)(q'') \mid q'' \in C_\Psi
\} \text{ (lifted to $[\Tsfstar[\partial X] X; J]$) near } q,
\ilabel{Ghat-leg}\end{equation}
and such that $q'$ corresponds to $q$ under this correspondence.

The non-degeneracy condition implies that there is a local diffeomorphism
between the set
$$ C_\Psi = \{ (s, x_1, u, y_1, y_2, y_3,v_2,v_3) \mid d_s \Psi = d_{v_2}
\Psi = d_{v_3} \Psi = 0 \text{ at } (s, x_1, su, y_1, y_2, y_3, v_2, v_3)
\}
$$ and $\Ghat$.

In the second case, as we are localizing away from the boundary of $\{ x_2
= 0 \}$, given by $x_2/s = 0$, we do not need the special variable $s$. In
this case, a non-degenerate parametrization of $(G, \Gsharp)$ near $q \in
\Ghat$ is a smooth phase function $\Psi$ of the form
\begin{equation}
1 + x_1 x_2 \psi(x_1, x_2, y_1, y_2, y_3, v) \ilabel{Psi2}\end{equation}
defined on a neighborhood of $q'=(0, 0, y_1^*, y_2^*, y_3^*, v^*)$ with
$\Psi$ non-degenerate in the sense that such that
\begin{equation}
d_{y_2} \psi, \ d_{y_2, v} \big( \frac{\partial \psi}{\partial v^i} \big)
 \text{ are linearly independent at } q' \ilabel{Psi2-nd}
\end{equation}
with
\begin{equation}\ilabel{Ghat-leg-2} \Ghat = \{ d\big( \frac{\Psi}{x_1x_2x_3}\big)(q'') \mid q'' \in C_\Psi
\} \text{ (lifted to $[\Tsfstar[\partial X] X; J]$) near } q,
\end{equation} and such that $q'$ corresponds to $q$ under this correspondence.

\subsection{Existence of parametrizations}
For brevity we only show the existence of paramet\-rizations in the
codimension 3 setting. The construction is analogous (and simpler) in the
codimension 2 setting. We use coordinates as in the proof of
Proposition~\ref{prop:param} above, in which we have $G_1 = \{ \nubar_1 = 1, \mu_1 =
0 \}$ and $G_2^\sharp = \{ \nubar_1 = (1+x_1) \nubar_2, \mu_1 = 0, \mu_2 =
0 \}$.

First let $q \in \Ghat$ lie on the codimension three corner of $\Ghat$.
Recall that $\Ghat$ fibres over $\Ghat_2$ with fibres that are Legendrian
submanifolds of $\Tscstar[\partial F] F$; therefore we can find a
splitting of the $y_3$ coordinates, $y_3 = (y_3^\flat, y_3^\sharp)$, so
that $(y_3^\sharp, \mu_3^\flat)$ form coordinates on the fibre over $\pi(q)
\in \Ghat_2$. Also, as in \cite{hasvas1}, Proposition 3.5, we can find a
splitting of the $y_2$ coordinates, $y_2 = (y_2^1, y_2^\flat, y_2^\sharp)$,
where $y_2^\flat = (y_2^2, \dots, y_2^j)$, so that, with $\muhat^\flat =
(\mu_2^2/\mu_2^1, \dots, \mu_2^j/\mu_2^1)$, $\muhat_2^\sharp =
\mu_2^\sharp/\mu_2^1$, the functions $(y_2^\sharp, \mu_2^1, \muhat^\flat)$
form coordinates $\Ghat_2$ near $\pi(q)$.  It follows that
$$\mathcal{Z} = (x_1, \mu_2^1, x_2/\mu_2^1, y_1, y_2^\sharp,
\muhat_2^\flat, y_3^\sharp, \mu_3^\flat)$$ form coordinates on $\Ghat$ near
$q$.

We now follow the proof of Proposition~\ref{prop:param} as closely as possible.
Writing $\nubar_i$, $y_i^\flat$ and $\mu_i^\sharp$ in terms of the
coordinates $\cZ$ on $\Ghat$, we have
\begin{equation}\begin{aligned}
\nubar_1 &= N_1(\cZ) \\ \nubar_2 &= N_2(\cZ) \\ \nubar_3 &= N_3(\cZ) \\
y_i^\flat &= Y_i^\flat(\cZ), i = 2, 3 \\ \mu_1 &= M_1(\cZ) \\
\muhat_2^\sharp &= M_2^\sharp(\cZ) \\ \mu_3^\sharp &= M_3^\sharp(\cZ)
\end{aligned} \qquad  \text{ on } \Ghat.
\end{equation}
Since $G$ is Legendrian, we have
\begin{equation}\begin{gathered}
dN_1 + N_2 dx_1 + x_1 N_3 dx_2 - M_1 \cdot dy_1 \\ - x_1 \mu_2^1 \big(
dY_2^1 + \muhat_2^\flat \cdot dY_2^\flat + M_2^\sharp \cdot dy_2^\sharp
\big) - x_1 x_2 \big( \mu_3^\flat \cdot dY_3^\flat - M_3^\sharp \cdot
dy_3^\sharp \big) = 0.
\end{gathered}\ilabel{contact00}\end{equation}
We claim that the function (where we substitute $s$ for $\mu_2^1$, $v_2$ for
$\muhat_2^\flat$ and $v_3$ for $\mu_3^\flat$)
\begin{equation}
\Psi(x_1, x_2, s, y_1, y_2, y_3, v_2, v_3) = N_1 + x_1 s \big( (y_2^1 -
Y_2^1) + (y_2^\flat - Y_2^\flat) \cdot v_2 \big) + x_1 x_2 \big( (y_3^\flat
- Y_3^\flat) \cdot v_3 \big)
\ilabel{Psidefn}\end{equation}
is a local parametrization of $G$.
%To avoid confusion, let us write $v_i$ instead of $\mu_i^\flat$ for the corresponding arguments of $\Phi$. 

First, observe that $N_1$ is equal to $1$ at $x_1 = 0$ and is equal to $1 +
x_1 + O(s)$ at $s=0$ since the value of $\nubar_1$ on $\Gsharp_2$ is equal
to $1+x_1$. Hence it has the form \eqref{Psi1}.

Second, suppose that $d_{s} \Psi = 0$. This means that
$$ d_{s} N_1 + x_1(y_2^1 - Y_2^1) - x_1 s d_s (Y_2^1 + Y_2^\flat \cdot v_2
\big) - x_1 x_2 d_s Y_3^\flat \cdot v_3 = 0.
$$ Using the $ds$ component of \eqref{contact00} and dividing by an overall
factor of $x_1$ we now obtain $y_2^1 = Y_2^1$.  In a similar way, the
conditions that $d_{v_i} \Psi = 0$ imply that $y_i^\flat = Y_i^\flat$, $i =
2, 3$. This also shows the non-degeneracy condition, since $d (d_s \psi_2)
= dy_2^1$, $d(d_{v_2} \psi_2) = y_2^\sharp$, $d(\partial_{v_3} \psi_3) =
dy_3$ at $q;$ these are manifestly linearly independent differentials.

To see that the set
$$ G' = \{ d \big( \frac{\Psi}{x_1 x_2 x_3} \big) \mid d_{s, v_2, v_3} \Psi
= 0 \}
$$ coincides with $G$ locally near $q$, consider the value of $\mu_2^1$ on
$G'$; it is given by $d_{y_2^1} \Psi/x_1 = s$. Similarly, the value of
$\mu_2^\flat$ is given by $s v_2$, and the value of $\mu_3^\flat$ is given
by $v_3$. So we can re-identify these values. Next consider the value of
$\nubar_1$ on $G'$. It is given by the value of $\Psi$, that is, by
\eqref{Psidefn}. This simplifies to $N_1$ when $d_{s,v_i} \Psi = 0$, since
we have $y_2^1 = Y_2^1$ when $d_s \Psi = 0$ and $y_i^\flat = Y_i^\flat$
when $d_{v_i} \Psi = 0$. Next consider the value of $\nubar_2$. This is
given by $d_{x_1} \Psi$ which is equal to
$$ d_{x_1} N_1 - x_1 s d_{x_1} Y_1^\flat \cdot \mu_1^\flat - x_1 s d_{x_1}
Y_2^\flat \cdot \mu_2^\flat - x_1 x_2 d_{x_1} Y_3^\flat \cdot \mu_3^\flat
$$ (again using $y_i^\flat = Y_i^\flat$ when $d_{s,v_i} \Psi = 0$). Since
the $dx_1$ component of \eqref{contact00} vanishes, this is equal to
$N_2$. So $\nubar_2 = N_2$ on $G'$. In a similar way we deduce that
$\nubar_3 = N_3$, and $\mu_i^\sharp = M_i^\sharp$ on $G'$. It follows that
$G'$ coincides with $G$.

\subsection{Equivalence of phase functions}\jw{More needed?}

We sketch the proof of equivalence of parametrizations only in the
codimension three case.

Two phase functions $\Psi$, $\tilde{\Psi}$ are said to be \emph{equivalent}
if they have the same number of phase variables of each type $v_1,v_2$ and
there exist maps
$$ V_1(x_1,\vecy,\vecv,s),\ V_2(x_1,\vecy,\vecv,s),\ S(x_1,\vecy,\vecv,s)
$$ such that
$$ \tilde{\Psi}(\vecx,\vecy,V_1,V_2,S) = \Psi.
$$

\begin{proposition}
The phase functions $\Psi = 1 + x_1 + sx_1 \psi_2+ x_1 x_2 \psi_3,$
$\tilde\Psi = 1 + x_1 + sx_1 \tilde\psi_2+ x_1 x_2 \tilde\psi_3$ are
locally equivalent iff
\begin{enumerate}
\item
They parametrize the same Legendrians,
\item
They have the same number of phase variables of the form $v_2,$ $v_3$
separately,
\item
\begin{align*}
\sgn d_{v_2}^2 (\psi_2) &= \sgn d_{v_2}^2 (\tilde\psi_2),\\ \sgn d_{v_3}^2
(\psi_3) &= \sgn d_{v_3}^2 (\tilde\psi_3),
\end{align*}
\end{enumerate}
\end{proposition}
By using the codimension two result from \cite{hasvas1}, we reduce to the
case
$$ \Psi = 1 + x_1 + sx_1 \psi_2+ x_1 x_2 \psi_3,\ \tilde\Psi =1 + x_1 +
sx_1 \psi_2+ x_1 x_2 \tilde \psi_3.
$$ As usual, we can arrange that the two functions agree to first order
along $C := \{d_{s,v_2}(s \psi_2+x_2\psi_3), d_{v_3}\psi_3 =0\}.$ Thus
$$ \tilde \psi_3-\psi_3 = \frac 12 (\nabla'_{\vecv,s}\Psi)^t B
(\nabla'_{\vecv,s}\Psi)
$$ where we define $\nabla' \Psi = (\pa_{v_2} (s\psi_2+\psi_3),
\pa_{v_3}\psi_3, \pa_s(s\psi_2+\psi_3)).$ We now expand
$$ \Psi(x_1,\vecy,\tilde{\vecv},\tilde s) - \Psi(x_1,\vecy,\tilde{\vecv},s)
= (\tilde{\vecv}-\vecv ) \cdot \pa_{\vecv}\Psi+(\tilde{s}-s ) \cdot
\pa_{s}\Psi + O((\tilde{\vecv}-\vecv )^2+(\tilde s-s)^2).
$$ Set
$$ (\tilde v_1, \tilde v_2,\tilde s) -(v_1,v_2, s) = (x_2 w_1, w_2, x_2
w_3)\cdot \nabla'_{\vecv}\psi
$$ for $w_i=w_i(x_1,\vecy,\vecv,s).$ Thus
$$ \Psi(x_1,\vecy,\tilde{\vecv},\tilde{s}) - \Psi(x_1,\vecy,\vecv,s) = x_1
x_2 (\nabla'_{\vecv,s}\Psi)^t (w+O(w^2)) (\nabla'_{\vecv,s}\Psi).
$$ We want
\begin{align*} \Psi(x_1,\vecy,\tilde{\vecv},\tilde{s}) - \Psi(x_1,\vecy,{\vecv},s) &=
\tilde\Psi(x_1,\vecy,\vecv,s) -\Psi(x_1,\vecy,\vecv,s)\\ &=x_1 x_2
(\tilde\psi_{3}-\psi_{3})
\end{align*}
We thus need to solve
$$ x_1 x_2 (\nabla'\Psi)^t (w+O(w^2)) (\nabla'\Psi) = \frac{x_1 x_2}2
(\nabla'\Psi)^t B (\nabla'\Psi)
$$ for $w.$ This can always be accomplished for $B$ small by the inverse
function theorem, and extended to the general case by using the condition
on signatures.

\subsection{Legendre distributions associated to a conic pair}\ilabel{Legdist-conic}

\subsubsection{Codimension 2 corners}
Let $X$ be a scattering fibred manifold with codimension 2 corners, let $N
= \dim X$ and let $(G, \Gsharp_1)$ be a conic Legendrian pair.  Let $m, p$
and $r$ be real numbers, and let $\nu$ be a smooth nonvanishing
scattering-fibred half-density. A Legendre distribution of order $(m, p;
r)$ associated to $(G, \Gsharp_1)$ is a half-density distribution of the
form $u_1 + (u_2 + u_3 + u_4 + u_5) \nu$, where 
\begin{itemize}

\item
$u_1$ is a
Legendre distribution of order $(m; r)$ associated to $G$ and microsupported
away from $J$, 

\item $u_2$ is given by an finite sum of local expressions
\begin{equation}\begin{gathered}\ilabel{eq:dist-1-c2}
u_2(x_1, x_2,y_1, y_2) = \int_{\RR^{k_2}} \int_{\RR^{k_1}} \int_0^\infty
e^{i\psi(s, x_1, y_1, y_2, v_1, v_2)/x_1 x_2} a(s, \frac{x_1}{s}, x_2, y_1,
y_2, v_1, v_2)\\ x_2^{m-(1+k_1+k_2)/2+N/4} \big(\frac{x_1}{s}\big)^{r - (1
+ k_1)/2 - f_1/2 + N/4} s^{p-1-f_1/2+N/4} \,ds \, dv_1 \, dv_2,
\end{gathered}\end{equation}
where $a$ is a smooth compactly supported function of its arguments, $f_1$ is the dimension of the fibres of $H_1$, and $\psi = 1 + s\psi_2 + x_1 \psi_2$ is a phase function locally parametrizing $(G, G_1^\sharp)$ near a point $q \in \partial_1 \Ghat \cap \partial_\sharp \Ghat$, as in  \eqref{psi1psi2}, 

\item $u_3$ is given by an finite sum of local expressions
\begin{equation}\begin{gathered}\ilabel{eq:dist-11-c2}
u_2(x_1, x_2,y_1, y_2) = \int_{\RR^{k}} e^{i\psi(x_1, y_1,
y_2, v)/x_1 x_2} \tilde a(x_1, x_2, y_1, y_2, v)\\ x_2^{m-k/2+N/4}
x_1^{p-1-f_1/2+N/4} \,dv,
\end{gathered}\end{equation}
where $\tilde a$ is smooth and compactly supported, and  $\psi$ is 
a local parametrization of $(G, G_1^\sharp)$ near a point $q \in   \partial_\sharp \Ghat \setminus \partial_1 \Ghat$ as in \eqref{psi-simple},

\item $u_4$ is given by 
\begin{equation}\begin{gathered}
u_4(x_1, y_1, z) = \int e^{i(1 + s \psi_1)/x_1 x_2} b(x_1, s,
 \frac{x_1}{s}, y_1, v, z) \\ \big( \frac{x_1}{s}\big)^{r - (1 + k_1)/2 -
 f_1/2 + N/4} s^{p-1-f_1/2+N/4} \, dv_2
\end{gathered}\ilabel{eq:dist-2-c2}
\end{equation}
where $\psi_1$ is as above and $b$ is smooth and $O(x_2^\infty)$ at $\mf = \{ x_2 = 0 \}$, 
and 

\item $u_5 \in x_1^{p - f_1/2 + N/4} x_2^\infty e^{i/x_1 x_2} \CI(X)$ (which always contains $\CIdot(X)$ as a subset). 
\end{itemize}

The set of such distributions is denoted $I^{m,p;r}(X, (G, \Gsharp_1); \sfOh)$.

\subsubsection{Codimension 3 corners}\ilabel{6.5.2}

We now assume that $X$ is a scattering-fibred manifold with codimension 3
corners. Let $N = \dim X$, let $m, r_1, r_2$ and $p$ be real numbers, and
let $\nu$ be a smooth nonvanishing scattering-fibred half-density on $X$. A
Legendre distribution of order $(m, p; r_1, r_2)$ associated to $(G,
\Gsharp_2)$ is a half-density distribution of the form $ u_1 + u_2 +  (
u_3 + u_4 + u_5 + u_6) \nu$, where 
\begin{itemize}

\item $u_1$ is a Legendre distribution of
order $(m; r_1, r_2)$ associated to $G$ and microsupported away from $J$, 

\item $u_2$ is a Legendre distribution of order $(m, p; r_2)$ associated to $(G, G_2^\sharp)$ and supported away from $H_1$, as defined above, 

\item $u_3$
is given by an finite sum of local expressions
\begin{equation}\begin{gathered}\ilabel{eq:dist-1a}
u_2(x_1, x_2,x_3,y_1, y_2,y_3) = \int\limits_{\RR^{k_3}} \int\limits_{\RR^{k_2}}
\int\limits_0^\infty e^{i\Psi(x_1, x_2, y_1, y_2,y_3, s, v_2, v_3)/\xx} a(x_1, s,
\frac{x_2}{s}, x_3,y_1, y_2,y_3, v_2, v_3)\\ x_3^{m-(1+k_2+k_3)/2+N/4}
\big(\frac{x_2}{s}\big)^{r_2 - (1 + k_2)/2 - f_2/2 + N/4} s^{p-1-f_2/2+N/4}
x_1^{r_1-f_1/2 + N/4}\,ds \, dv_2 \, dv_3,
\end{gathered}\end{equation}
where $a$ is a smooth compactly supported function of its arguments, $f_i$ are the dimension of the fibres on $H_i$, and $\Psi = 1 + x_1 + sx_1 \psi_2 + x_1 x_2 \psi_3$ is a local parametrization of $(G, G_2^\sharp)$ near a corner point $q$ as in \eqref{Psi1}, 

\item $u_4$ is given by an finite sum of local expressions
\begin{equation}\begin{gathered}\ilabel{eq:dist-1b}
u_2(x_1, x_2,x_3,y_1, y_2,y_3) = \int_{\RR^{k}} \int_0^\infty e^{i\Psi(x_1,
x_2, y_1, y_2,y_3, v)/\xx} \tilde a(x_1, x_2, x_3,y_1, y_2,y_3, v)\\
x_3^{m-k/2+N/4} x_2^{p-1-f_2/2+N/4} x_1^{r_1-f_1/2 + N/4}\,dv,
\end{gathered}\end{equation}
where $\tilde a$ is smooth and compactly supported,  $\Psi$ is  a local parametrization of $(G, G_2^\sharp)$ near a point $q \in \partial_1 \Ghat \cap \partial_\sharp \Ghat \setminus \partial_2 \Ghat$ as in \eqref{Psi2},  

\item $u_5$ is given by 
\begin{equation}\begin{gathered}
u_4(x_1, x_2, y_1, y_2, z_3) = \int e^{i(1 + sx_1 \psi_2)/\xx} b(x_1, s,
 \frac{x_2}{s}, y_1, y_2, z_3, v_2) \\ \big( \frac{x_2}{s}\big)^{r_2 - (1 +
 k_2)/2 - f_2/2 + N/4} s^{p-1-f_2/2+N/4} x_1^{r_1-f_1/2 + N/4}\,dv_1 \,
 dv_2
\end{gathered}\ilabel{eq:dist-2a}
\end{equation}
where $\psi_2$ is as above,  $b$ is smooth and $O(x_3^\infty)$ at  $\mf$, and

\item $u_6 \in x_1^{r_1 - f_1/2 + N/4} x_2^{p - f_2/2 + N/4} x_3^\infty e^{i(1+x_1)/\xx} \CI(X) $ (which includes $\CIdot(X)$ as a subset).
\end{itemize}
 The set of such distributions is denoted $I^{m,p;r_1, r_2}(X, (G, \Gsharp_2); \sfOh)$. 

\subsection{Symbol calculus}

\subsubsection{Codimension 2 corners}
For a conic pair of Legendre submanifolds $\Gt = (G, G^\sharp_1)$, with $\Gh$ the
desingularized submanifold obtained by blowing up $J = \phit_{12}^{-1}(\spnn  G^\sharp_1)$, the symbol
calculus takes the form

\begin{proposition}\ilabel{prop:ex-conic-2} Let $s$ be a boundary
defining function for $\partial_\sharp \Ghat \subset \Ghat$. Then there is an exact sequence
\begin{multline}
0 \to I^{m+1,p; r}(X, \Gt; \sfOh) \to I^{m,p; r}(X, \Gt; \sfOh)
\\ \to x_1^{r-m} s^{p-m} \CI(\Gh, \Omega^\half_b \otimes
S^{[m]}(\Gh)) \to 0.
\ilabel{ex-conic-2}\end{multline}

If $P
\in \Diffsf(X; \sfOh)$ has principal symbol $p$ and $u \in  I^{m,p; r}(X, \Gt; \sfOh)$, then $Pu
\in  I^{m,p; r}(X, \Gt; \sfOh)$ and
$$ \sigma^m(Pu) = \big( p \rest \Ghat \big) \sigma^m(u).
$$ Thus, if $p$ vanishes on $\Ghat$, then $Pu$ is an element of
$ I^{m+1,p; r}(X, \Gt; \sfOh)$ by \eqref{ex-conic-2}. The
symbol of order $m+1$ of $Pu$  in this case is given by
\eqref{transport}.

\end{proposition}

\subsubsection{Codimension 3 corners} Let $\Gt = (G, G^\sharp_2)$ now be a conic pair of Legendre submanifolds in the codimension three setting. Then we have 

\begin{proposition}\ilabel{prop:ex-conic-3} Let $s$ be a boundary
defining function for $\partial_\sharp \Ghat \subset \Ghat$, and let $\rho$ be a boundary defining function for $\partial_2 \Ghat$ (for example, $\rho = x_2/s$). Then there is an exact sequence
\begin{multline}
0 \to I^{m+1,p; r_1, r_2}(X, \Gt; \sfOh) \to I^{m,p; r_1, r_2}(X, \Gt; \sfOh)
\\ \to x_1^{r_1-m} \rho^{r_2 - m} s^{p-m} \CI(\Gh, \Omega^\half_b \otimes
S^{[m]}(\Gh)) \to 0.
\ilabel{ex-conic-3}\end{multline}

If $P
\in \Diffsf(X; \sfOh)$ has principal symbol $p$ and $u \in  I^{m,p; r_1, r_2}(X, \Gt; \sfOh)$, then $Pu
\in  I^{m,p; r_1, r_2}(X, \Gt; \sfOh)$ and
$$ \sigma^m(Pu) = \big( p \rest \Ghat \big) \sigma^m(u).
$$ Thus, if $p$ vanishes on $\Ghat$, then $Pu$ is an element of
$ I^{m+1,p; r_1, r_2}(X, \Gt; \sfOh)$ by \eqref{ex-conic-2}. The
symbol of order $m+1$ of $Pu$  in this case is given by
\eqref{transport}.

\end{proposition}

\subsection{Residual space}\ilabel{resid-conic} In the codimension two case, consider the case where $X = Y \times [0, \epsilon]_{x_2}$ where $Y$ is a manifold with boundary. In this case, the residual space
$$
I^{\infty, p; r}(X, (G, G^\sharp_1); \sfOh) = \cap_{m} I^{m, p; r}(X, (G, G^\sharp_1); \sfOh)
$$
may be identified with
$$
x_2^\infty C^\infty \big( [0, \epsilon]; I^{r-1/4, p-1/4}(X, (x_2^{-1}G_1, x_2^{-1}G_1^\sharp), \sfOh) \big).
$$

In the codimension three case, if $X = Y \times [0, \epsilon]_{x_3}$ where $Y$ is a scattering-fibred manifold with codimension two corners, then the residual space is
$$
I^{\infty, p; r_1, r_2}(X, (G, G^\sharp_2); \sfOh) = \cap_{m} I^{m, p; r_1, r_2}(X, G, G^\sharp_2); \sfOh)
$$
and this may be identified with 
$$
x_3^\infty C^\infty \big( [0, \epsilon]; I^{r_2 - 1/4, p-1/4; r_1 - 1/4}(X, (x_3^{-1}G_2, x_3^{-1}G_2^\sharp), \sfOh) \big) \otimes |dx_3|^{1/2}.
$$

%%%%%%%%%%%%%%%%%%%%%%%%%%%%%%%%%%%%%%%%%%%%
\section{Legendrian-Lagrangian distributions}

\subsection{Legendrian-Lagrangian submanifolds}
The final type of distribution we shall introduce are
`Legendrian-Lagrangian distributions' associated to the scattering
cotangent bundle $\Tscstar X$ of a manifold with boundary $X$. We shall
restrict attention to $X$ of the form $X = Y \times [0, h_0)$. We ignore the noncompactness of $Y \times [0, h_0)$ as $h \to h_0$
since we will only be interested in distributions supported near the
boundary at $h=0$.

 Let $\Tscstarbar X$ be the compactification of $\Tscstar X$ via radial
 compactification of each fibre. This is a manifold with corners of
 codimension two; its boundary hypersurfaces are the fibrewise radial compactification of $\Tscstar[\partial X] X$, which we denote $\sinfty$ (`semiclassical limit'), and the new hypersurface at
 `fibre-infinity', which we shall denote $\finfty$. Fibre-infinity has a
 natural contact structure given by  $\rho \sum_i \eta_i dy_i $ in local coordinates $y$ on $Y$ (where $\eta$ are the dual cotangent coordinates), where $\rho$ is a boundary defining
 function  for $\finfty$ (e.g. $\rho = 1/|\eta|$). If $\eta_1/|\eta| > 0$ locally then we may take $\rho = 1/\eta_1$ and then the contact form takes the form $dy_1 + \sum_{i \geq 2} \eta_i/\eta_1 dy_i$. 
 
 There is a natural subbundle $S$ of $\Tscstar[\partial X] X$ given by the annihilator of $h^2 \partial_h$, or equivalently, spanned by the one-forms $dy_i/h$. Let $\partial S \subset \finfty \cap \sinfty$ denote the boundary of $S$ after radial compactification. 
 
 \begin{defn} A Legendrian-Lagrangian submanifold on $X$ is a
 Legendre submanifold with boundary $L \subset \sinfty$ that meets the corner $\sinfty \cap \finfty$ transversally, and such that $\partial L \subset \partial S$. 
 \end{defn}
 
Recall that given local coordinates $y$ on $Y$, we have
 coordinates $h, y, \nu, \mu$ on $\Tscstar X$ near $\{ h = 0 \}$ given by
 writing any element of $\Tscstar X$ relative to the basis $d(1/h)$ and
 $dy_i/h$:
 $$ \Tscstar X \ni p = \nu d\big( \frac1{h} \big) + \sum_i \mu_i
\frac{dy^i}{h} .
$$ The coordinates $(\nu, \mu)$ are linear coordinates on each fibre, and $S$ is given by $\{ h = 0, \nu = 0 \}$. Now
let $q$ be a point on the corner of $\partial L$ after radial
compactification of the fibres. Let us assume for a moment that $\mu_1/|\mu| > 0$ at $q$ (which can always be arranged after a linear change of $y$ variables), so that we can use $\rho =1/\mu_1$ as a boundary defining
function for $\finfty$ near $q$. Let $\sigma = \nu/\mu_1$ and $M' =
\mu'/\mu_1$, where $\mu' = (\mu_2, \dots, \mu_{m})$, $m = \dim Y$. Then
$(h, y, \rho, \sigma, M')$ are local coordinates for $\Tscstarbar X$ near
$q$.  

At $\sinfty$, the contact structure is given by the form $d\sigma - dy_1 -
M' \cdot dy' - \sigma d\rho/\rho$. This form vanishes on $L$. Therefore at
$\partial L,$ which is contained in $\{ \sigma = 0 \}$ by assumption, we
have $dy_1 + M' \cdot dy' = 0$. Thus, $\pa L$ can be naturally identified
with a Legendrian at fibre-infinity on the fibrewise compactification of
$T^* Y$, and hence with a conic Lagrangian $\Lambda$ in $T^* Y \setminus
0$. We shall soon see that a Legendrian-Lagrangian distribution is, for
fixed $h > 0$, a Lagrangian distribution associated to $h^{-1} \Lambda$.

\subsection{Parametrization}
Let $q \in \partial L$. We shall use coordinates $(h, y, \rho, \sigma, M')$
as above.  We recall that $\sigma = 0$ at $q$, indeed everywhere on
$\partial L$.

A local parametrization of $L$ near $q$ is a function $\Phi/\overline{\rho}$,  where
 $\Phi = \Phi(y, \overline{\rho}, v)$ is a smooth function of $y$, $\overline{\rho}$ and $v \in
 \RR^k$, defined in a neighbourhood of $(y_0, 0, v_0)$ so that
$$ d_{\overline{\rho},v}(\Phi/\overline{\rho})|_{y_0,0,v_0} =0
$$ and
$$ q = \Big(0, y_0; d_{h,y} \big( \frac{\Phi}{\overline{\rho} h} \big) \Big),
$$ such that $\Phi$ is non-degenerate in the sense that
$$ d \big( \frac{\partial \Phi}{\partial v_i} \big), \quad d\Phi \text{ and
} d\overline{\rho} \text{ are linearly independent at } (y_0, 0, v_0),
$$ and so that
\begin{equation}
L = \big\{ \big(0, y, d_{h,y} \big( \frac{\Phi}{\overline{\rho} h} \big):
d_{\overline{\rho},v}\big( \frac{\Phi}{\overline{\rho}}\big) =0 \big) \big\}
\ilabel{LL}\end{equation} locally near $q$.

This is a parametrization using `compact coordinates'. We may also use
noncompact or homogeneous coordinates by introducing $w \in \RR^{k+1}$
given in terms of $(\overline{\rho}, v)$ by $w = (w_1, w')$ with $w_1 = 1/\overline{\rho}$ and
$w' = v/\overline{\rho}$. Also write $\Phi = \Phi_1(y, v) + \overline{\rho} \Phi_0 (y, \overline{\rho},
v)$. Then changing to the $w$ variables we have a parametrization of the
form
$$ \Psi_1(y, w) + \Psi_0(y, w)
$$ where $\Psi_1 = w_1 \Phi_1$ is homogeneous of degree $1$ in $w$ and
$\Psi_0 = \Phi_0$ is a symbol of order zero in $w$. Then $\Psi_1$
parametrizes the Lagrangian $\Lambda$.

\subsection{Existence of parametrizations} This is proved in the usual way. Let $y = (y_1, y', y'')$, $\mu = (\mu_1, \mu', \mu'')$, $\xi' = \mu'/\mu_1$ and $\xi'' = \mu''/\mu_1$.  We choose coordinates so that $(\rho, \xi', y'')$ are coordinates on $L$ near $q$. (Note that $\rho$ always has nonzero differential on $L$ at $q$ since $L$ is assumed transverse to $\finfty$.) We can therefore express the other coordinates on $L$ as smooth functions of $(\rho, \xi', y'')$:
\begin{equation*}\begin{aligned}
y_1 &= Y_1(\rho, \xi', y'') \\ y' &= Y'(\rho, \xi', y'') \\ \sigma &=
\Sigma(\rho, \xi', y'') \\ \xi'' &= \Xi''(\rho, \xi', y'')
\end{aligned}\end{equation*}
We claim that
\begin{equation*}
\frac{\Phi(y, \overline{\rho}, \overline{\xi'})}{\overline{\rho}} = \frac{ \Sigma(y, \overline{\rho}, \overline{\xi'}) + (y_1 - Y_1(y, \overline{\rho}, \overline{\xi'})) + (y' - Y'(y, \overline{\rho}, \overline{\xi'})) \cdot
\xi'}{\overline{\rho}}
\end{equation*}
parametrizes $L$ near $q,$ with $v=\xi'.$ In fact, since $L$ is Legendrian,
the form
$$ -d\Sigma + \Sigma \frac{d\rho}{\rho} + dY_1 + \xi' \cdot dY' + \Xi''
\cdot dy''
$$ vanishes on $L$. Setting the coefficients of $d\rho$, d$\xi'$ and $dy''$
to zero we find that
\begin{align}
- \dbyd{\Sigma}{\rho} + \frac{\Sigma}{\rho} + \dbyd{Y_1}{\rho} + \xi' \cdot
\dbyd{Y'}{\rho} &= 0 \ilabel{Phi1} \\ 
\dbyd{\Sigma}{\xi'} + \dbyd{Y_1}{\xi'} + \xi'
\cdot \dbyd{Y'}{\xi'} &= 0 \ilabel{Phi2} \\ 
\dbyd{\Sigma}{y''} + \dbyd{Y_1}{y''} + \xi'
\cdot \dbyd{Y'}{y''} + \Xi'' &= 0 \ilabel{Phi3}
\end{align}
Using \eqref{Phi1} and \eqref{Phi2}, one finds that $d_{\overline\rho, v} (\Phi/\overline{\rho}) = 0$ implies that $Y_1 = y_1$
and $Y' = y'$ on $L$, while equating $d_{y_1} (\Phi/\overline\rho)$ with $1/\rho$ and $dy' (\Phi/\overline\rho)$ with $\xi'/\rho$ gives $\overline\rho = \rho$ and $\overline\xi' = \xi'$. Finally one obtains \eqref{LL}  near $q$.

\subsection{Equivalence of phase functions}
Acceptable changes of variables for our phase function are smooth
coordinate changes of the form $(\rho,v) \mapsto (\tilde\rho, \tilde v)$
where $\tilde\rho =\rho f,$ with $f \in \CI.$ In the noncompact model
described above, this is equivalent to $w \mapsto \tilde w$ where $\tilde
w$ is a polyhomogeneous symbol of order $1$ in the $w$ variables.  We
therefore declare two phase functions to be equivalent if such a
transformation maps one to the other.  We continue to employ the noncompact
phase variable description of parametrization in what follows.

\begin{proposition}
The phase functions $\Psi$, $\tilde\Psi$ are locally equivalent near $q$ iff
\begin{enumerate}
\item
They parametrize the same Legendrian,
\item
They have the same number of phase variable.
\item
$ \sgn d_{w}^2 \Psi = \sgn d_{w}^2 \tilde \Psi$ at $q$.
\end{enumerate}
\end{proposition}
\begin{proof}
%% Let $C$ and $\tilde C$ denote the respective
%% sets where $d_{w} \Psi=0,$ $d_{w} \tilde\Psi=0.$

%% By nondegeneracy and the implicit function theorem there exists a map
%% $\Upsilon$ such that
%% $$w= \Upsilon(y,d_y\Psi,d_w\Psi),
%% $$ and a corresponding $\tilde \Upsilon.$  The map $\Upsilon$ necessarily
%% has an asymptotic expansion in the $d_y\Psi,d_w\Psi$ variables, beginning
%% with degree $1$ and descending in degree.\jw{For reasons which completely
%%   elude me.}

%% We now form the map
%% $$
%% \Gamma: (y,w) \mapsto (y,\tilde\Upsilon(y,d_y\Psi, d_w\Psi)+
%% A\pa_{w}\Psi).
%% $$
%% where $A$ is a linear transformation to be determined later.  Note that
%% this map is a polyhomogeneous symbol in the fiber variables.

%% The map $\Gamma$ takes $C$ to $\tilde C$ and commutes with the map to
%% $L.$  We may now choose $A$ so that $\Gamma$ is a diffeomorphism
%% in a conic neighborhood of $C$ for $\abs{w}$ sufficiently large.

%% We now improve this result to exact equivalence of $\Psi$ and $\tilde \Psi$
%% under the assumption that the function agree to second order on
%% $C.$

We begin as usual by arranging to have $\Psi$ and $\tilde\Psi$ in agreement
to first order along $C = \{d_w\Psi=0\}.$

We may thus expand in a Taylor series
$$ \tilde \Psi - \Psi = \frac 12 (\nabla_{w}\Psi)^t B (\nabla_{w}\Psi)
$$ for some matrix $B=B(y,w).$ As both $\tilde\Psi$ and $\Psi$ are symbols of order $1$
in $w$, $B$ is also symbolic of order $1$.

The non-degeneracy assumptions on $\Psi, \tilde \Psi$ means precisely that
$\det (I+B\pa_{w}^2\psi_2)\neq 0.$ We now expand
$$ \Psi(y,\tilde{w}) - \Psi(y,w) = (\tilde{w}-w ) \cdot \pa_{w}\Psi +
O((\tilde{w}-w)^2).
$$ Set
$$ \tilde w-w = z \cdot \nabla_w\Psi,$$ 
where $z$ is a matrix depending on $w$; note that this is a change of variables of the required form.  We thus
have
$$ \Psi(y,\tilde{w}) - \Psi(y,w) = (\nabla_{w}\Psi)^t (z+O(z^2)) (\nabla_{w}\Psi).
$$ We want
$$\Psi(y,\tilde{w}) - \Psi(y,w) = \tilde\Psi(y,w) -\Psi(y,w);
$$ we thus need to solve
$$ (\nabla_{w}\Psi)^t (z+O(z^2)) (\nabla_{w}\Psi)=\frac 12
(\nabla_{w}\Psi)^t B (\nabla_{w}\Psi)
$$ for $z.$ This can be accomplished for $B$ small by the inverse function
theorem, with a result that is symbolic in $w.$  Lemma~3.1.7 of
\cite{Ho:FIO1} enables us to extend to the case of arbitrary $B.$
\end{proof}

\subsection{Legendrian-Lagrangian distributions}\label{llds}
Let $L$ be a Legendrian-Lagrangian submanifold as above. Let $N = \dim X$.  A
Legendrian-Lagrangian distribution $u$ of order $(m,r)$ associated to $L$
on $X$, denoted $u \in I^{m,r}(X, L)$, is a half-density $u = u_1 + u_2 + u_3 + u_4$, where 
\begin{itemize}

\item $u_1$ is in $h^\infty C^\infty \big( [0, h_0); I^{-r-1/4}(Y, h^{-1}\Lambda; \Omega^\half) \big) \otimes |dh|^{1/2}$,

\item $u_2$ is a Legendrian
distribution of order $m$ associated to $L$ and microsupported away from
fibre-infinity, 

\item $u_3$ is a sum of terms of the form
\begin{equation}
h^{m - (k+1)/2 + N/4} \int_0^\infty \int_{\RR^k}  e^{i\Phi(y, \rho, v)/\rho h} \rho^{r - k/2 - 1 - N/4}
 a(h, y, \rho, v) \, dv \, d\rho \big|
\frac{dy dh}{h^{N+1}} \big|^{1/2}
\ilabel{LLd}\end{equation}
where $\Phi$ is a local parametrization of $L$ and $a$ is smooth, and

\item $u_4 \in \CIdot(X)$.
\end{itemize} 
We remark
that \eqref{LLd} is an oscillatory integral unless $r$ is sufficiently positive.

If we rewrite this using the homogeneous parametrization, we get
\begin{equation}
h^{m - (k+1)/2 + N/4} \int e^{i\Psi_1(y, w)/h} e^{i\Psi_0(y,w)/h} \tilde
a(h, y, w) \, dw \big| \frac{dy dh}{h^{N+1}} \big|^{1/2}
\end{equation}
where $\tilde a$ is a (classical) polyhomogeneous symbol of order $-r - (k+1)/2 +N/4-1/2$ in $w$. For $h > 0$ the
$e^{i\Psi_0/h}$ factor is a symbol of order zero and so $e^{i\Psi_0/h} \tilde a$ is a symbol of the same order as $\tilde a$.  Hence for fixed $h > 0$ this is a Lagrangian distribution, of order $-r - 1/4$, associated to $\Lambda$
depending smoothly on $h$ for $h > 0$, and whose symbol is itself oscillatory as $h \to 0$. 

\subsection{Symbol calculus and residual space} 

\begin{proposition} The symbol map for Legendre distributions, defined in the interior of $G$ \cite{MZ}, extends by continuity to 
give  an exact sequence\ah{these need to be checked carefully}
$$ 0 \to I^{m+1, r}(X, L; \sfOh) \to I^{m, r}(X, L; \sfOh) \to
\rho^{r-N/4} \CI(L, \Omega^\half_b \otimes S^{[m]}(G)) \to 0. 
$$

The residual space $I^{\infty, r}(X, L; \sfOh) \equiv \cap_m I^{m, r}(X, L; \sfOh)$ may be identified with
$$
h^\infty C^\infty \big( [0, h_0); I^{-r-1/4}(Y, h^{-1}\Lambda; \Omega^\half) \big) \otimes |dh|^{1/2}.
$$
\end{proposition}

\subsection{Distributional limits of Legendrian conic pairs}
We now consider a situation which leads to a Legendrian-Lagrangian
distribution. Let $M$ be a compact manifold with boundary. We may view $M
\times [0, h_0)$ as a scattering-fibred manifold (where we again ignore
the noncompactness at $h=h_0$) with main face $H_2 = M \times \{ 0 \}$
and other boundary hypersurface $H_1 = \partial M \times [0, h_0)$ with
fibration $H_1 \to \partial M$ given by projection onto the $\partial M$
factor. Suppose that we have a distribution $u \in I^{}(G, \Gsharp_1)$
associated to a Legendrian conic pair $(G, \Gsharp_1)$ as described in the
previous section (codimension 2 case). Coordinates in this case can be
taken to be $x, h, y$ near $H_1$, where $y$ is a local coordinate for
$\partial M$, extended to a collar neighbourhood of $\partial M$, and $x$
is a boundary defining function for $M$. Corresponding scattering-fibred
cotangent coordinates are $\nu_1, \nu_2, \mu$ given by expressing covectors
$$ q = \nu_1 d \big( \frac1{xh} \big) + \nu_2 d \big( \frac1{h} \big) + \mu
\cdot \frac{dy}{xh}.
$$ We suppose as in Section~\ref{conic-param-2} that $\Gsharp_1$ is given by $\{ \mu = 0, \nu_1 = 1 \}$.

We consider the problem of restricting $u \in I^{m,p;r}(X, (G, \Gsharp_1))$
to $H_1$. The first issue is that $u$ is a half-density, so to restrict we
must divide by the half-density $|dx/x|^{1/2}$ to obtain a half-density on
$H_1$.  The second issue is that $u$ is oscillatory as $x \to 0$, so we
must first multiply by $e^{-i/xh}$ to have any hope of being able to
restrict to $x=0$. Thirdly we must divide by a power of $x$, depending on
the order $p$ at $\Gsharp_1$, in order to get a finite, nonzero limit.  If
we do all this, and if the Legendrian $G$ intersects $\{ \mu = 0 \}$ only
at $\Gsharp_1$, then it turns out that $u$ has a restriction in the
distributional sense. In the non-semiclassical case (no $h$ variable) this
was proved by Melrose-Zworski \cite{MZ}.

Let $\Ghat$ be the blowup of the singular Legendrian $G$ at $J = \{ x = 0,
\mu = 0 \} \subset \Tsfstar[\mf] (M \times [0, h_0))$. We write $\tilde J$
for the new boundary hypersurface created by the blowup. Recall that
$\Ghat$ is a manifold with corners of codimension 2, with one boundary
hypersurface $\partial_\sharp \Ghat$ at $\tilde J$ and the other, $\partial_1 \Ghat$,  at $x = 0$ but away from $\tilde
J$.

\begin{lemma}\ilabel{LL2} The submanifold $\tilde J \cap \{ \nu_1  = 1 \}$
  of $\tilde J$ is naturally diffeomorphic to the fibrewise
  compactification of $\Tscstar[\partial M \times \{ 0 \}] (\partial M)
  \times [0, h_0))$, and under this identification, the boundary
    hypersurface $\partial_\sharp \Ghat$  (which lies inside  $\tilde J
    \cap \{ \nu_1  = 1 \}$) is a Legendrian-Lagrangian submanifold $L$. 
\end{lemma}

\begin{proof} Scattering-fibred covectors in $\Tsfstar(M \times [0,h_0))$
  are represented by forms of the form $d\big((f(y) + x g)/xh\big)$ where $g$ is
  smooth. If we restrict to the set $\nu_1 = 1, \mu = 0$ then these are of
  the form $d((1 + xg)/xh) = d(1/xh) + d(g/h)$, and the coordinates are given
  by $y, \nu_1 = 1, \nu_2 = g, \mu = xdg$. Thus on the interior of $\tilde
  J$, where we may take $x$ as a boundary defining function, the
  coordinates are given by $y, \nu_1 = 1, \nu_2 = g, \mu/x = d g$. It is now
  clear that the map from the point on $\tilde J \cap \{ \nu_1 = 1 \}$
  specified by $(y, \nu_2 = g, \mu/x= d g)$ to $d(g/h) \in \Tscstar[\partial M
    \times \{ 0 \}](\partial M \times [0,h_0))$ is a natural
  diffeomorphism. This identifies the interior of $\tilde J \cap \{ \nu_1 =
  1 \}$ with $\Tscstar[\partial M \times \{ 0 \}](\partial M \times
  [0,h_0))$. The fibres of the blowdown map $\tilde J \to J$ are radial
  compactifications of vector spaces coordinatized by $\mu/x$, since we
  have $x/|\mu|$ as a boundary defining function for $\partial \tilde J$ and
  $\mu/|\mu|$ as a coordinate along the boundary. Hence the natural
  diffeomorphism extends from $\tilde J \cap \{ \nu_1 = 1 \}$ to the
  fibrewise radial compactification of $\Tscstar[\partial M \times \{ 0
    \}](\partial M \times [0,h_0))$. 

The contact form on $\Tsfstar[\mf](M \times [0,h_0))$ is given by $d\nu_1 +
x d\nu_2 + \mu \cdot dy$. Let $\eta = \mu/x$. Then this can be written
$d\nu_1 + x (d\nu_2 + \eta \cdot dy),$ which vanishes at $\Ghat$.  Taking
the differential and restricting to $\partial_\sharp \Ghat$ we get $dx \wedge
(d\nu_2 + \eta \cdot dy) = 0$, and since $dx \neq 0$ at the interior of
$\tilde J$ we conclude that $d\nu_2 + \eta \cdot dy = 0$ at $\partial_\sharp \Ghat$. Since $\nu = 1$ on $\Gsharp_1$, we have $\nu_1 = 1$ at $\partial_\sharp \Ghat$. Using our identification of the interior of $\tilde J \cap
\{ \nu_1 = 1 \}$ with $\Tscstar (\partial M \times [0, h_0))$, we see that
the image of $\partial_\sharp \Ghat$ under this identification, which we
denote $L$, is Legendrian. Also, by the transversality requirements in the
definition of a Legendrian conic pair, $L$ is transverse to the boundary at
fibre-infinity. Finally, since $\Ghat$ is a compact submanifold of
$[\Tsfstar[\mf](M \times [0,h_0)); J]$, and since $\nu_2$ is a continuous
function on this space, the value of $\nu_2$ is bounded on $\Ghat$. On the
other hand, the value of $\eta = \mu/x$ goes to infinity at the boundary of
$L$. Hence $\nu_2/\eta = 0$ at the boundary of $L$. It follows that $\partial L \subset \partial S$, so $L$ is
a Legendrian-Lagrangian submanifold.
\end{proof}

The analytic result corresponding to this geometric lemma is

\begin{proposition}\ilabel{distlim} Let $X$ be a scattering-fibred manifold
  with codimension 2 corners, let $\dim X = N$, and let $(G, \Gsharp_1)$ be
  as above. Suppose that $u \in I^{m,p;r}(X, (G, \Gsharp_1))$, and assume
  that $G \cap \{ \mu = 0 \}$ is contained in $\Gsharp_1$. Then
\begin{equation}
x^{-p + N/4} e^{-i/xh} \big| \frac{dx}{x} \big|^{-1/2} u
\ilabel{limq}\end{equation} has a distribution limit as $x \to 0$. The
limit is an element of $I^{m-1/4, p-r-(N-1)/4}(H_1, L)$, where $L$ is as in
Lemma~\ref{LL2}.
\end{proposition}

\begin{proof} By definition, $u$ is a sum of terms $u_i$ as in the
  definition above \eqref{eq:dist-1-c2}.  Clearly we can ignore any
summands which are rapidly decreasing at $x = 0$.  We next note that, if we
microlocalize $u$ to any region where $\mu \neq 0$, then (non-)stationary
phase (involving repeated integrations by parts in $y$) shows that the
pairing of $u$ with any function of $y$ is rapidly decreasing in $x$ as $x
\to 0$. Hence we can restrict attention to the microlocal region where
$\mu$ is close to zero, which by assumption is near $\Gsharp_1$, i.e.\ near
the conic singularity of $G$.

In this region, we have seen that $u$ can be written as a sum of terms of
the form \eqref{eq:dist-1-c2} and \eqref{eq:dist-11-c2} (with $y_2$ and
$v_2$ absent and $x_1$ replaced by $x$). Consider an integral of the form
\eqref{eq:dist-1-c2}:
\begin{equation*}\begin{gathered}
\int_{\RR^k} \int_0^\infty e^{i(1 + s\psi_1(s,y,v) + x \psi_2(s, x/s, y,
v))/xh} a(s, \frac{x}{s}, y, v, h) \\ \times \  h^{m - (1+k)/2 + N/4} s^{p-1/2
+ N/4} \big( \frac{x}{s} \big)^{r - (1+k)/2 - 1/2 + N/4} \, \frac{ds}{s} \,
dv \, \mu
\end{gathered}\end{equation*}
where $\mu$ is a scattering-fibred half-density. We may take $\mu$ to be
$$ \mu = \Big| \frac{dy}{(xh)^{N-2}} \frac{dx}{x^2 h} \frac{dh}{h^2}
\Big|^{1/2}.
$$ We want to express this in terms of a scattering half-density $\nu = |dy
dh/h^{N}|^{1/2}$; we see that
$$ \mu = h^{-1/2} x^{-(N-1)/2} \big| \frac{dx}{x} \big|^{1/2} \nu.
$$ It follows that \eqref{limq} is given by
\begin{equation*}\begin{gathered}
\int_{\RR^k} \int_0^\infty e^{i(s\psi_1(s,y,v) + x \psi_2(s, x/s, y,
v))/xh} a(s, \frac{x}{s}, y, v, h) h^{m - (1+k)/2 + N/4 - 1/2} \\ \times \
\big( \frac{x}{s} \big)^{p-1/2 + N/4} \big( \frac{x}{s} \big)^{r - (1+k)/2
- 1/2 + N/4} \, \frac{ds}{s} \, dv \, \nu.
\end{gathered}\end{equation*}
Let us introduce the variables $\eta_1 = s/x$ and $\eta' = v\eta_1 \in
\RR^k$. We can write this
\begin{equation}\begin{gathered}
\int_{\RR^k} \int_0^\infty e^{i(\eta_1\psi_1(x\eta_1,y,\eta'/\eta_1) +
\psi_2(x\eta_1, \eta_1^{-1}, y, \eta'/\eta_1))/h} a(x\eta_1,
\frac{1}{\eta_1}, y, \eta'/\eta_1, h) \\ \times h^{m - (1+k)/2 + N/4 - 1/2}
\eta_1^{p-r+(1+k)/2} \eta_1^{-k} \frac{d\eta_1}{\eta_1} \, d\eta' \, \nu.
\end{gathered}\ilabel{wtt}\end{equation}
If we set $x=0$ in the integrand then the integral becomes
\begin{equation*}\begin{gathered}
\int_{\RR^k} \int_0^\infty e^{i(\eta_1\psi_1(0,y,\eta'/\eta_1) + \psi_2(0,
\eta_1^{-1}, y, \eta'/\eta_1))/h} a(0, \frac{1}{\eta_1}, y, \eta'/\eta_1,
h) \\ \times h^{m - (1+k)/2 + N/4 - 1/2} \eta_1^{p-r-(1+k)/2} \, d\eta_1 \,
d\eta' \, \nu.
\end{gathered}\end{equation*}
It is straightforward to check that $\eta_1 \psi_1(0,y,\eta'/\eta_1) +
\psi_2(0, \eta_1^{-1}, y, \eta'/\eta_1)$ is a non-degenerate
parametrization of $L$.\ah{Prove this?} Therefore this is a
Legendrian-Lagrangian distribution of order $(m-1/4, p-r-(N-1)/4)$
associated to $L$. It remains to prove that this is indeed the
distributional limit of \eqref{wtt} as $x \to 0$. This is clear if the
exponent of $\eta_1$ is sufficiently negative, since then the integral is
absolutely convergent, uniformly in $x$. In general, we can exploit the
fact that $d_{y^1} \psi \neq 0$ according to \eqref{psi-nd} (where $y^1$ is the first component of $y$) and integrate
by parts repeatedly in $y^1$, using
$$ e^{i\eta_1 \psi_1/h} = \frac{h}{i\eta_1} \frac1{\pa_{y^1} \psi_1} \pa_{y^1}
e^{i\eta_1 \psi_1/h}.
$$ Doing this sufficiently many times eventually reduces the exponent of
$\eta_1$ to the point of absolute integrability. We can then take the limit
$x \to 0$ and perform the integrations by parts in reverse, which gives the
desired conclusion.

A similar argument applied to an integral of the form \eqref{eq:dist-11-c2}
gives the same result (in this case, we only get a Legendrian distribution
microlocalized to a compact part of the interior of $L$ since this part is
away from the corner of $\Ghat$).
\end{proof}

%%%%%%%%%%%%%%%%%%%%%%%%%%%%%%%%

\section{Quadratic scattering-fibred structure}\ilabel{qsfs}
In order to describe precisely the microlocal structure of the Schr\"odinger propagator, we need to introduce the \emph{quadratic-scattering fibred structure} on manifolds with codimension three corners. 
This structure is a variant of the scattering-fibred structure, in which
we have an extra order of vanishing of the Lie algebra at some of the
boundary hypersurfaces. The basic example, on a manifold with boundary, is
the quadratic scattering structure, which we now review.

\subsection{The basic structure} 
Recall that the quadratic scattering structure on a manifold with boundary, $X$,  is given by the quadratic scattering Lie algebra $\Vqsc(X) \equiv x \Vsc(X)$. Locally near the boundary, using coordinates $(x,y)$, $x \geq 0$ a boundary defining function, $\Vqsc$ is the $\CI(X)$-span of the vector fields
$$
x^3 \partial_x, \quad x^2 \partial_{y_i}.
$$
This structure was used to analyze the propagation of singularities at
infinity of solutions to the time-dependent Schr\"odinger equation
\cite{Wunsch1}, \cite{R-Z}. 

In the quadratic scattering-fibred structure on a manifold with codimension 3 corners, we start with a manifold $X$ with fibrations $\phi_i$, as in Definition~\ref{2.3}. However, instead of a distinguished total boundary defining function $\xx$ we require a distinguished function $\qx$ which vanishes to second order at the $H_1$ and $H_2$ boundary hypersurfaces; in other words $\qx = x_1^2 x_2^2 x_3$ for some boundary defining functions $x_i$ of $H_i$. Correspondingly, we consider a different Lie algebra of vector fields. In place of Definition~\ref{sfvf}, we make

\begin{defn}The Lie algebra of quadratic scattering-fibred vector fields $\Vqsf$ is defined by
\begin{equation}
V \in \Vqsf(X) \text{ iff }  V = x_1 x_2 W, \quad W(\qx) \in x_1^3 x_2^3 x_3^2 \CI(X) \text{ and } W
\text{ is tangent to } \Phi.
\end{equation}
\end{defn}

%As before we may check that this is a Lie algebra. 
An analogue of Proposition~\ref{X-coords} applies, where we replace the last condition $\Pi x_i = \xx$ by $x_1^2 x_2^2 x_3 = \qx$. In terms of such coordinates, the Lie algebra is given locally by arbitrary linear combinations (over $\CI(X)$) of vector fields of the form
\begin{equation}\begin{aligned}
& -(x_1^2x_2^2x_3) x_1 \partial_{x_1}, & (x_1^2 x_2^2 x_3) \partial_{y_1}, \\ & 
x_1 x_2^2 x_3 (x_1 \partial_{x_1} - x_2 \partial_{x_2}), & (x_1x_2^2 x_3) \partial_{y_2},
\\ & x_1 x_2 x_3 (x_2 \partial_{x_2} - 2x_3 \partial_{x_3}), & x_1 x_2 x_3 \partial_{y_3}.
\end{aligned}
\ilabel{sf-qvfs}
\end{equation}

It
follows, as in Section~\ref{two}, that $\Vqsf(X)$ is the space of sections of a vector bundle over $X$. The dual bundle, denoted $\Tqsfstar(X)$, is spanned by one-forms of the form $d(f/(x_1^2 x_2^2 x_3))$ where $f\in\CIsf(M)$. 

The dual basis to the vector fields \eqref{sf-qvfs} is
\begin{equation}
d \big( \frac{1}{x_1^2 x_2^2 x_3} \big), \quad d \big(\frac1{x_1x_2^2 x_3} \big),
\quad d\big(\frac1{x_1 x_2 x_3}\big), \quad \frac{dy_1}{x_1^2 x_2^2 x_3}, \quad
\frac{dy_2}{x_1x_2^2x_3}, \quad \frac{dy_3}{x_1x_2x_3}.  \ilabel{qsf-cvs1}
\end{equation}
Here $dy_i$ is shorthand for a $k_i$-vector of 1-forms, if $y_i \in
\RR^{k_1}$.  An alternative basis is given by
\begin{equation*}
d \big( \frac{1}{x_1^2 x_2^2 x_3} \big), \quad\frac{dx_1}{x_1^2x_2^2 x_3} ,
\quad \frac{dx_2}{x_1 x_2^2 x_3}, \quad \frac{dy_1}{x_1^2 x_2^2 x_3}, \quad
\frac{dy_2}{x_1x_2^2x_3}, \quad \frac{dy_3}{x_1x_2x_3}.  \ilabel{qsf-cvs2}
\end{equation*}
Any element of $\Tsfstar X$ may therefore be written uniquely as
\begin{equation}
\tilde \nu_1 d\big(\frac{1}{x_1^2 x_2^2 x_3}\big) + \tilde\nu_2
d\big(\frac1{x_1 x_2^2 x_3} \big)
+ \tilde\nu_3 d\big(\frac1{x_1 x_2 x_3}\big) + \tilde\mu_1 \cdot \frac{dy_1}{x_1^2 x_2^2 x_3} +
\tilde\mu_2 \cdot \frac{dy_2}{x_1 x_2^2 x_3} + \tilde\mu_3 \cdot
\frac{dy_3}{x_1 x_2 x_3}
\ilabel{qcan}\end{equation} 
or, alternatively, as
\begin{equation}
\overline {\tilde \nu}_1 d\big(\frac{1}{x_1^2 x_2^2 x_3}\big) + \overline
	  {\tilde\nu}_2 \frac{d x_1}{x_1^2 x_2^2 x_3}
+ \overline {\tilde\nu}_3 \frac{dx_2}{x_1 x_2^2 x_3} +  \tilde\mu_1 \cdot \frac{dy_1}{x_1^2 x_2^2 x_3} +
 \tilde\mu_2 \cdot \frac{dy_2}{x_1 x_2^2 x_3} +  \tilde\mu_3 \cdot
	  \frac{dy_3}{x_1 x_2 x_3}
\ilabel{qcana}\end{equation} 
The function $\tilde\nu_1$, regarded as a linear form
on the fibres of $\Tsfstar X$, can be identified with the vector field
$(x_1^2 x_2^2 x_3) x_1 \partial_{x_1}$, and similarly for the other fibre
coordinates. The same expression can be viewed as the canonical one-form on
$\Tqsfstar X$. Taking $d$ of \eqref{can} therefore gives the symplectic form
on $\Tqsfstar X$.

The same reasoning as in Section~\ref{two}, but considering differentials of the form $d(f/(x_1^2 x_2^2 x_3)$ where $f \in \CIsf(X)$,  leads to the definition of the bundles $\Tqsfstar(F_i, H_i)$ and $\Nqsfstar Z_i$ as well as the induced fibrations $\phit_i$. 

The contact form on $\Tqsfstar [\mf] X$ is defined by contracting the symplectic form $\omega$ with $(x_1^2 x_2^2 x_3) x_3 \partial_{x_3}$ and restricting to $\mf$. This gives 
\begin{multline}
\chi = d\tilde \nu_1 + x_1 d\tilde\nu_2 + x_1 x_2 d\tilde\nu_3 - \tilde\mu_1 \cdot dy_1 - x_1 \tilde\mu_2
\cdot dy_2 - x_1 x_2 \tilde\mu_3 dy_3 \\
= d\overline{\tilde \nu}_1 - \overline{\tilde\nu}_2 dx_1 - x_1\overline{\tilde\nu}_3 dx_2 - \tilde\mu_1 \cdot dy_1 - x_1 \tilde\mu_2
\cdot dy_2 - x_1 x_2 \tilde\mu_3 dy_3
\end{multline}
which is exactly the same expression as the contact form in the scattering-fibred case. In a similar way we get induced contact forms on $\Nqsfstar Z_i$ and on the fibres of $\phit_i$. 

Given a scattering-fibred manifold, $Y$, with codimension two corners, 
we can form the product $X_h = Y \times [0, h_0)$ and endow it with the structure of a scattering-fibred manifold with codimension three corners, as in Section~\ref{sec:legendrian}, or we can form the product 
$X_t = Y \times [0, t_0]$ and endow it with the structure of a
  \emph{quadratic} scattering-fibred manifold with codimension three
  corners. It turns out that there is a contact transformation $Q$ between
  $\Tsfstar_{\mf} X_h \setminus \Nullset$ and $\Tqsfstar_{\mf} X_t
  \setminus \Nullset$, where $\Nullset$ denotes the subbundle of
  $\Tsfstar_{\mf} X_h$, resp.\ $\Tqsfstar_{\mf} X_t$,  spanned, at $p \in
  \mf$, by elements of the form $d(f/(x_1 x_2 x_3))$, resp. $d(f/(x_1^2
  x_2^2 x_3))$, where $f \in \CIsf (X)$ vanishes at $p$. This contact
  transformation is very useful in relating the semiclassical resolvent and
  the propagator (in the case that $Y = \MMb$). The map $Q$ is defined by 
\begin{equation}
Q \big( \frac{df}{x_1 x_2 h} \big) = \big( \frac{d(f^2)}{2x_1^2 x_2^2 t} \big), \quad f \in \CIsf(X_h), \ f \neq 0.
\ilabel{J}\end{equation} 
The proof of this is very straightforward if we use the coordinates
$\nubar_i$, $\mubar_i$ from \eqref{nubar} on $\Tsfstar_{\mf} X_h$ and the
analogous coordinates on $\Tqsfstar_{\mf} X_t$. Then, with $\tilde \chi$
the contact form on $\Tqsfstar_{\mf} X_t$, we find that $Q^* (\tilde \chi) =
\nubar_1 \chi$, showing that $Q$ is a contact transformation away from
$\nubar_1 = 0$, which is the set denoted $\Nullset$ above. We remark that
such contact transformations can be defined far more generally; the
point of the transformation $f \mapsto f^2/2$ is that shows up when we obtain the propagator from the resolvent via an integral over the spectral measure.

\subsection{Legendre distributions} 

The theory of Legendre distributions on quadratic scattering-fibred
manifolds proceeds in parallel to that of scattering-fibred manifolds.

\begin{defn} A Legendre
submanifold of a quadratic scattering-fibred manifold $X$ of dimension $N$ is a submanifold $G$ of dimension $N-1$ of $\Tqsfstar[\mf] X$ on
which the contact form $\chi$  vanishes, and such that $G$ is
transverse to each boundary $\Tqsfstar[\mf \cap H_i] X $ of $\Tqsfstar[\mf] X$. 
\end{defn}

A parametrization of a quadratic scattering-fibred Legendre submanifold can be defined much as in the scattering-fibred case. 
In fact, near a point $q \in G \cap \Tqsfstar_{\mf \cap H_1 \cap H_2} X$, the definition is identical to that of Section~\ref{leg-param} except that we replace \eqref{Gleg} by 
\begin{equation}
G=\{d\big( \frac{\psi}{(x_1^2 x_2^2 x_3)} \big) \mid (\vecx,\vecy,\vecv) \in C_\psi\}.
\ilabel{Gqleg}\end{equation}
We also give the definitions for a local parametrization near a point $q \in G$ lying in the interior of $G$, or in the interior of one of the boundary hypersurfaces of $G$. In the former case this is just a standard Legendre parametrization locally. In the latter case the definition is analogous to the codimension 3 case except that, near the interior of $H_1$ we split the coordinates as $(x_1, x_3, Y_1 = y_1, Y_3 = (y_2, x_2, y_3))$, our phase function is of the form 
$$
\psi(x_1, Y_1, Y_3, v_1, v_3) = \psi_1(Y_1, v_1) + x_1 \psi_2(x_1, Y_1, Y_3, v_1, v_3)
$$
and we ignore the variables with a `2' subscript. Near $H_2$ we split the coordinates $(x_2, Y_2 = (x_1, y_1, y_2), x_3, Y_3 = y_3)$, our phase function is of the form 
$$
\psi(x_2, Y_2, Y_3, v_2, v_3) = \psi_1(Y_2, v_2) + x_2 \psi_2(x_2, Y_2, Y_3, v_2, v_3)
$$
and we ignore the variables with a `1' subscript.

Let $m,r_1,r_2$ be real numbers, let $N = \dim X$, let $G \subset \Tqsfstar[\mf] X$ be a quadratic Legendre submanifold, and let $\nu$ be a smooth
nonvanishing quadratic scattering-fibred half-density. The set of quadratic  Legendre distributions of
order $(m; r_1, r_2)$ associated to $G$, denoted $I^{m, r_1, r_2}(X, G; \qsfOh)$, is the set of half-density distributions that can be written in 
the form $(u_1 + u_2 + u_3 + u_4 + u_5) \nu$, such that  
\begin{itemize}
\item $u_1 \cdot \nu$ is a quadratic Legendre distribution of order $(m; r_1)$ associated to $G$ and
supported away from $H_2$, i.e. $u_1$ is given by a finite sum of expressions of the form 
\begin{equation}\begin{gathered}
\int e^{i\psi(x_1, x_2, \vecy,\vecv)/(x_1^2 x_2^2 x_3)} a(\vecx,\vecy,\vecv)\\
x_3^{m-(k_1+k_2)/2+N/4} 
x_1^{r_1-k_1 - k_2/2-f_1/2 + 3N/4}\,dv_1 \, dv_2,
\end{gathered}\end{equation}
where $\psi$ parametrizes $G$ locally and $a$ is smooth and supported away from $x_2 = 0$, 
\item $u_2 \cdot \nu $ is similarly a quadratic Legendre distribution of order $(m; r_2)$ associated to $G$ and supported away from $H_1$,  

\item $u_3$ is given by an finite sum of
local expressions of the form
\begin{equation}\begin{gathered}\ilabel{eq:dist-1-q}
\int e^{i\psi(x_1, x_2, \vecy,\vecv)/(x_1^2 x_2^2 x_3)} a(\vecx,\vecy,\vecv)\\
x_3^{m-(k_1+k_2+k_3)/2+N/4} x_2^{r_2 - (k_1 +  k_2) - k_3/2 - f_2/2 + 3N/4}
x_1^{r_1-k_1 - (k_2 + k_3)/2 -f_1/2 + 3N/4}\,dv_1 \, dv_2 \, dv_3,
\end{gathered}\end{equation}
with $v_i \in \RR^{k_i},$ $a$ smooth and compactly supported,  $f_i$ the dimension of the
fibres of $H_i$ and $\psi = \psi_1 + x_1 \psi_2 + x_1 x_2 \psi_3$ a phase function locally parametrizing $G$ near a corner point $q \in \partial_{12} G$, 
\item $u_4$ is  given by  a finite sum of terms of the form
\begin{equation}\begin{gathered}
  \int e^{i(\psi_1 + x_1 \psi_2)/(x_1^2 x_2^2 x_3)} b(x_1,
 x_2, x_3, y_1, y_2, y_3) \\ x_2^{r_2 - (k_1 + k_2) - f_2/2 + 3N/4}
 x_1^{r_1-k_1 - k_2/2-f_1/2 + 3N/4}\,dv_1 \, dv_2
\end{gathered}\ilabel{eq:dist-2-q}
\end{equation}
with $\psi_1, \psi_2$ and $f_i$ as above, $b$ smooth with support compact and $O(x_3^\infty)$ at $\mf$, and

\item $u_5 \in \CIdot (X)$. 
\end{itemize}

If $X = X_t = Y \times [0, t_0]$ as in the previous subsection and $G$ is disjoint from $\Nullset$ then we can locally write $G = Q(G')$ for some Legendrian $G' \subset \Tsfstar_{\mf} X_h$; then if $\phi/x_1 x_2 h$ is a local parametrization of $G'$, $\phi^2/2(x_1^2 x_2^2 t)$ is a local parametrization of $G$.

\begin{proposition}\ilabel{sph} Suppose that $u_h \in I^{m, r_1, r_2}(X_h, G; \sfOh)$
  is a Legendre distribution associated to the Legendrian $G$ which does
  not intersect $\Nullset$. Also suppose that $\chi(t)$ is a smooth function of $t \in \RR$ that vanishes on $[0,R]$ and is identically equal to $1$ on $[2R, \infty]$, for some $R > 0$. Then the integral in $h$
\begin{equation}
\int_0^\infty e^{-it/2h^2} \chi(\sqrt{t}/h)  u_h \, \frac{|dh dt|^{1/2}}{h\phantom{a^{1/2}}}
\ilabel{intoverh}\end{equation}
is in $$I^{m+1/2,r_1 + m + 1/2 - N/4, r_2 + m + 1/2 - N/4}(X_t, Q(G); \qsfOh),$$
i.e.\ is a quadratic Legendre distribution associated to  $Q(G)$, with orders shifted by $1/2$ at $\mf$ and $m + 1/2 - N/4$ at $H_1$ and $H_2$. 
\end{proposition}

\begin{remark} Our interest in this lemma is for the following reason: if $u_h$ is $(2\pi i)^{-1}$ times the difference of the limit of the semiclassical resolvent on the spectrum, taken from above and below, 
$$
u_h^\pm = \frac1{2\pi i} 
\Big( \big( h^2 \Delta +2 h^2 V- (1 + i0) \big)^{-1} -\big( h^2 \Delta +2
h^2 V- (1 - i0) \big)^{-1} \Big) \otimes |dh|^{1/2}, 
$$
then the integral above gives the Schr\"odinger propagator $e^{-it
  (\Delta/2+V)}$ (times $|dt|^{1/2}$).   Note that the $V$ term is of a higher
semiclassical order in this setting than in the usual semiclassical
resolvent, hence $V$ does not affect the Legendrian geometry of the
Schr\"odinger propagator.
\end{remark}

\begin{proof} Locally $u$ may be written in the form
\begin{equation}\begin{gathered}
\int e^{i\psi/(x_1 x_2 h)} a(\vecx,\vecy,\vecv)\\
x_3^{m-(k_1+k_2+k_3)/2+N/4} x_2^{r_2 - (k_1 + k_2)/2 - f_2/2 + N/4}
x_1^{r_1-k_1/2-f_1/2 + N/4}\,dv_1 \, dv_2 \, dv_3 \cdot \nu
\end{gathered}\end{equation}
The condition on $G$ means that $|\psi| \geq \epsilon > 0$ at $G$; by cutting off the symbol close to $G$ we may assume that $|\psi| \geq \epsilon$ everywhere on the support of the symbol. Then in the integral \eqref{intoverh}  we get a phase function of the form
$-t/2h^2 + \psi/x_1 x_2 h$. Changing variable to $k = t x_1 x_2/ h$ this becomes $(-k^2/2 + k\psi)/(x_1^2 x_2^2 t )$, while the symbol becomes a function of $t x_1 x_2 /k$. Due to the $\chi$ cutoff, the integral in $k$ is supported in $k \geq R x_1 x_2 \sqrt{t}$. 

Let us insert cutoff functions $1 = \chi_1(k) + \chi_2(k)$, where $\chi_1$ is supported in $\{ k \leq \epsilon/2 \}$ and $\chi_2$ is supported in $\{ k \geq \epsilon/4 \}$. 

With $\chi_1$ inserted, there are no stationary points in the integral in $k$ since the phase is stationary when $k = \psi$. This term is in $\CIdot(X)$, as follows by writing
$$
e^{i(-k^2/2 + k\psi)/x_1^2 x_2^2 t} = \Big( i \frac{x_1^2 x_2^2 t}{k - \psi} \partial_k \Big)^N e^{i(-k^2/2 + k\psi)/x_1^2 x_2^2 t}
$$
and integrating by parts $N$ times, for arbitrary $N$. We gain at least
$\sqrt{t} x_1 x_2$ with each integration-by-parts.  

With $\chi_2$ inserted, we avoid the singularity caused by the argument $t x_1 x_2 /k$ in the symbol, and this term is a quadratic Legendre distribution associated to $Q(G)$ since the phase function $\psi k - k^2/2$ parametrizes $Q(G)$. Collecting powers of $t, x_1$ and $x_2$ (bearing in mind that the number $k_1$ of $v_1$ variables has increased by $1$ due to the appearance of $k$) completes the proof. 
\end{proof}

\subsection{Conic pairs}  

We now give an analogous sketch of the theory of 
Legendre distributions associated to conic Legendrian pairs on quadratic scattering-fibred manifolds.

\begin{defn} Let $G_2^\sharp$ be a projectable Legendrian in $\Nqsf Z_2$, and let $G \subset \Tqsfstar[\mf] X$ be a Legendre submanifold that is singular at  the boundary. Let $J_2$ be the span of $G_2^\sharp$ in $\Nqsf Z_2$, and $J$ the preimage of $J_2$ in $\Tqsfstar[\mf \cap H_2] X$ under the map $\phit_2$. We say that $(G, G_2^\sharp)$ form an conic Legendrian pair if $G$ has conic singularities at $J$, i.e.\ lifts to $[\Tqsfstar[\mf] X; J]$ to a smooth manifold $\Ghat$ transverse to $\tilde J$ as well as  to the lifts of $\Tqsfstar[\mf \cap H_1] X$ and $\Tqsfstar[\mf \cap H_2] X$. 
\end{defn}

A parametrization of a quadratic scattering-fibred Legendre submanifold can be defined much as in the scattering-fibred case. 
In fact, near a point $q \in G \cap \Tqsfstar_{\mf \cap H_1 \cap H_2} X$, the definition is identical to that of Section~\ref{conic-param-3} except that we replace \eqref{Ghat-leg} and \eqref{Ghat-leg-2} by 
\begin{equation}
\Ghat=\{d\big( \frac{\Psi}{(x_1^2 x_2^2 x_3)} (q'') \big) \mid q'' \in C_\Psi\} \text{ (lifted to $[\Tqsfstar[\mf] X; J]$) near $q$ }. 
\ilabel{Gqleg2}\end{equation}

If $X = X_t = Y \times [0, t_0]$ as in the previous subsection and $G$ is disjoint from $\Nullset$ then we can locally write $G = Q(G')$ for some Legendrian $G' \subset \Tsfstar_{\mf} X_h$; then if $\phi/x_1 x_2 h$ is a local parametrization of $G'$, $\phi^2/2(x_1^2 x_2^2 t)$ is a local parametrization of $G$. To simplify the definition of a distribution associated to a quadratic conic Legendrian pair, assume that this is the case. Then we can use the parametrizations of $G'$ from Section~\ref{Legdist-conic}. 

Let $N = \dim X$, let $m, r_1, r_2$ and $p$ be real numbers, and
let $\nu$ be a smooth nonvanishing quadratic scattering-fibred half-density on $X$. A quadratic 
Legendre distribution of order $(m, p; r_1, r_2)$ associated to $(G,
\Gsharp_2)$ is a half-density distribution of the form $ u_1 + u_2 +  (
u_3 + u_4 + u_5 + u_6) \nu$, where 
\begin{itemize}

\item $u_1$ is a Legendre distribution of
order $(m; r_1, r_2)$ associated to $G$ and microsupported away from $J$, 

\item $u_2$ is a Legendre distribution of order $(m, p; r_2)$ associated to $(G, G_2^\sharp)$ and supported away from $H_1$, 
\item $u_3$
is given by an finite sum of local expressions
\begin{equation}\begin{gathered}\ilabel{eq:dist-1a-q}
u_2(\vecx, \vecy) = \int\limits_{\RR^{k_3}} \int\limits_{\RR^{k_2}}
\int_0^\infty e^{i\Psi^2(x_1, x_2, \vecy, s, v_2, v_3)/2\qx} \ a(x_1, s,
\frac{x_2}{s}, x_3,\vecy, v_2, v_3) 
x_3^{m-(1+k_2+k_3)/2+N/4}  \\  \times \ 
\big(\frac{x_2}{s}\big)^{r_2 - (1 + k_2) -  k_3/2 - f_2/2 + 3N/4} s^{p-1-f_2/2+3N/4}
x_1^{r_1-f_1/2 -(k_2 + k_3)/2+ 3N/4}\,ds \, dv_2 \, dv_3,
\end{gathered}\end{equation}
where $a$ is a smooth compactly supported function of its arguments, $f_i$ are the dimension of the fibres on $H_i$, and $\Psi = 1 + x_1 + sx_1 \psi_2 + x_1 x_2 \psi_3$ is a local parametrization of $(G', (G')_2^\sharp)$ near a corner point $q$ as in \eqref{Psi1}, 

\item $u_4$ is given by an finite sum of local expressions
\begin{equation}\begin{gathered}\ilabel{eq:dist-1b-q}
u_2(\vecx,\vecy) = \int\limits_{\RR^{k}} \int_0^\infty e^{i\Psi^2(x_1,
x_2, \vecy, v)/2\qx} \ a(\vecx,\vecy, v)\\ \times \ 
x_3^{m-k/2+N/4} x_2^{p-1-k/2-f_2/2+3N/4} x_1^{r_1-k/2-f_1/2 + 3N/4}\,dv,
\end{gathered}\end{equation}
where $\tilde a$ is smooth and compactly supported,  $\Psi$ is  a local parametrization of $(G', (G')_2^\sharp)$ near a point $q \in \partial_1 \Ghat \cap \partial_\sharp \Ghat \setminus \partial_2 \Ghat$ as in \eqref{Psi2},  

\item $u_5$ is given by 
\begin{equation}\begin{gathered}
u_3(x_1, x_2, y_1, y_2, z_3) = \int e^{i(1 + sx_1 \psi_2)^2/2\qx} \ b(x_1, s,
 \frac{x_2}{s}, y_1, y_2, z_3, v_2) \\ \big( \frac{x_2}{s}\big)^{r_2 - (1 +
 k_2)/2 - f_2/2 + 3N/4} s^{p-1-f_2/2+3N/4} x_1^{r_1-k_2/2-f_1/2 + 3N/4}\,dv_1 \,
 dv_2 
\end{gathered}\ilabel{eq:dist-2aa}
\end{equation}

where $\psi_2$ is as above,  $b$ is smooth and $O(x_3^\infty)$ at  $\mf$, and

\item $u_6 \in x_1^{r_1 - f_1/2 + 3N/4} x_2^{p - f_2/2 + 3N/4} x_3^\infty e^{i(1+x_1)^2/2\qx} \CI(X) $ (which includes $\CIdot(X)$ as a subset).
\end{itemize}
 The set of such distributions is denoted $I^{m,p;r_1, r_2}(X, (G, \Gsharp_2); \qsfOh)$.

\begin{proposition}\ilabel{sph-conic} Suppose that $u_h \in I^{m, p; r_1, r_2}(X_h, G; \sfOh)$ is a Legendre distribution associated to the Legendrian $G$ which does not intersect $\Nullset$. Then
$$
\int e^{-it/2h^2} u_h \, \frac{|dh dt|^{1/2}}{h\phantom{{}^{1/2}}}
$$
is in $I^{m+1/2, p+m+1/2-N/4; r_1 + m + 1/2 - N/4, r_2 + m + 1/2 - N/4}(X_t, Q(G); \qsfOh)$, i.e.\ is a quadratic Legendre distribution associated to  $Q(G)$, with orders shifted by $1/2$ at $\mf$ and $m + 1/2 - N/4$ at $J$, $H_1$ and $H_2$. 
\end{proposition}

The proof is identical to that of Proposition~\ref{sph}.

%%%%%%%%%%%%%%%%%%%%%%%%%%%%%%%%

\part{Resolvent}

\section{The example of Euclidean space}\ilabel{Euc}

In this section we look at the structure of the resolvent, Poisson operator, scattering matrix and propagator 
on Euclidean space, with the flat metric and no potential, from a Legendrian point of view, and show explicitly that these kernels obey the claims made in Theorems~\ref{main1}, \ref{sc-thm} and \ref{prop-thm}, and Corollaries~\ref{sp-cor} and \ref{Poisson-cor}. 

We begin with the outgoing resolvent kernel. We may identify functions and half-densities via the Riemannian half-density, and regard our kernel as acting on half-densities. The kernel itself is then a half-density on $\RR^{2n}$. In order to fit into the framework here we need to multiply by a half-density in $h$, so that the kernel becomes a half-density on $X$. Which power of $h$ to
include with this half-density factor is an arbitrary choice. We will adopt the convention that the semiclassical outgoing ($+$)/incoming ($-$) resolvent is  
$$ (h^2 \Delta - (1 \pm i0))^{-1} |dh|^{1/2}.
$$
The difference of these, multiplied  by $(2\pi i)^{-1}|dh/h^2|^{1/2}$, is  the spectral measure $dE(\lambda^2)$ ($\lambda = h^{-1}$). 
%making the semiclassical order $3/4$. This has the advantage that it
%restricts to the canonical Poisson operator. Other choices are also
%possible.

The kernel of the outgoing resolvent is then 
\begin{equation}R_+ = h^{-n} e^{i|z-z'|/h} f_n(|z-z'|/h) |dz dz'|^{1/2} |dh|^{1/2} =
e^{i|z-z'|/h} f_n(|z-z'|/h) h \mu
\ilabel{fres}\end{equation} 
where $\mu$ is a nonvanishing scattering-fibred half-density, and
$f_n(t) \sim c_n t^{-(n-1)/2}$ as $t \to \infty$. 

Let us compute the orders of this as a Legendrian distribution at $N^* \Delta_b$ (see \eqref{ndb}), at the propagating Legendrian $L_+$, and at the b-face and the left and right boundary. We recall that the order convention is that the order gets larger as the distribution becomes smaller, i.e.\ more regular, and that the order is $N/4$ if the distribution is borderline $L^2$. To determine the order at $N^* \Delta_b$, we microlocalize away from $L_+$ by inserting a cutoff function $\chi$ that vanishes in a neighbourhood of $\{ |\zeta| = 1 \}$ and write the kernel as an oscillatory integral:
$$ h \int e^{i(z-z') \cdot \zeta/h} (|\zeta|^1 - 1)^{-1} \chi(\zeta) \,
d\zeta \mu.
$$
On the one hand, this is a semiclassical pseudodifferential operator of
order $(-2, 0, 0)$; on the other hand, it is  a Legendre distribution
associated to $N^* \Delta_b$, of semiclassical order $m=1 + n/2 - (2n+1)/4 =
3/4$  and order at $\bfc$ equal to $r_2 = n/2+1/2 - (2n+1)/4 = 1/4$ (by
\eqref{eq:dist-1}). To determine the order at $L_+$ we use the expression
\eqref{fres}; then \eqref{ild}  gives the semiclassical order as $m=(n+1)/2
- (2n+1)/4 = 1/4$ and the order at $\bfc$ equal to $r_2 = (n-1)/2 - (2n+1)/4 +
1/2 = -1/4$. Note that both these orders are $1/2$ less than the
corresponding order at $N^* \Delta_b$ as required for an intersecting
Legendre distribution (see Section~\ref{ilds}).
The order at $H_1$, i.e.\ the left or right
boundaries, is calculated from \eqref{fres} to be  $r_1 = s_1 + k/2 + f_1/2 - N/4 = (n-1)/2 + (n+1)/2 - (2n+1)/4 =
n/2 - 1/4.$ 

In the case of the free resolvent, the Legendrian $L_+$ is smooth at
$L^\sharp$. However, when this is true, by writing the phase and the symbol
in polar coordinates around the intersection $L_+ \cap L^\sharp$ we can
regard an element of $I^{m; r_1, r_2}(X,L; \sfOh)$ as an element of $I^{m,
  r_2+d; r_1, r_2}(X,(L, L^\sharp); \sfOh)$ where $d$ is the codimension of
the intersection; here, $d = (n-1)/2$. (This is explained in more detail in
section 14 of \cite{MZ}.) Thus we see that the free outgoing resolvent
kernel is an element of $$\Psi^{2, 0, 0}(X) + I^{1/4; -1/4}(X,(N^* \Delta_b,
L_+); \sfOh) + I^{1/4, n/2 - 3/4; n/2 - 1/4, -1/4}(X,(L_+, L^\sharp);
\sfOh).$$ 

We recall that for a fixed $h > 0$, the semiclassical order has no meaning while the other orders must be adjusted by adding $1/4$, reflecting the fact that the orders are `zeroed' using $N/4$ where $N$ is the total dimension. We see then that the orders agree with those claimed in \cite{hasvas1} for the resolvent at a fixed energy. 

The Poisson kernel has a natural normalization: we can ask that the family 
$\{ P(\lambda) \} $, $\lambda = h^{-1} \in (0, \infty)$, form a unitary operator mapping from $M$ to $L^2(\partial M
\times \RR_+; \lambda^{n-1} d\lambda d\omega)$ with measure corresponding to the conic metric $d\lambda^2 + \lambda^2 d\omega$ (i.e.\
a scattering metric) (see \cite{hasvas1}, section 9). To do this we need to multiply the Poisson operator of \cite{MZ} and \cite{hasvas1} by the half-density $|dh/h^2|^{1/2}$. 

To obtain the Poisson kernel at $\rb$ we divide $R_+$ by $|dr'|^{1/2}e^{-i|z'|/h}$, where $r' = |z'|$, and restrict at $r' = \infty$, i.e.\ at $\rb$, 
to get a half-density at $\rb$: we get
$$ P(h^{-1}) = \lim_{|z'| \to \infty} e^{-i|z'| h} e^{i|z-z'|/h} f_n(|z-z'|/h) h \mu
|dr'|^{-1/2}.
$$ The limit of $e^{-i|z'| h} e^{i|z-z'|/h}$ as $\abs{z'} \to \infty$ is
$e^{-iy' \cdot z/h}$, where $y' = z'/|z'|$. Also, $\mu \abs{d {r'}}^{-1/2}$ is equal to $h^{-1/2}$
$\times (|z'|/|z|)^{(n-1)/2}$ times a nonzero scattering-fibred half-density $\nu$
on the Poisson space $M \times \pa M \times [0, h_0),$ since
\begin{align*}
\Big| \frac{dh}{h^2} \frac{dx}{x^2 h} \frac{dx'}{(x')^2 h}
\frac{dy}{(xh)^{(n-1)}} \frac{dy'}{(x'h)^{(n-1)}} \Big|^{\half} =& \ \Big|
\frac{dh}{h^2} \frac{dx}{x^2 h} \frac{dy}{(xh)^{(n-1)}}
\frac{dy'}{(xh)^{(n-1)}} \Big|^{\half} \abs{d r'}^{\half} h^{-\half} \big(
\frac{x}{x'} \big)^{\frac{n-1}{2}} \\
=& \ \nu  \abs{d r'}^{\half} h^{-\half} \big(
\frac{x}{x'} \big)^{\frac{n-1}{2}}.
\end{align*}
Therefore, using the asymptotic $f_n(t) \sim c_n t^{-(n-1)/2}$ as $t \to \infty$, we have 
$$P(h^{-1})= c_n (x')^{-(n-1)/2} h^{(n-1)/2} h e^{-iy' \cdot z/h} \nu h^{-1/2}
\big(\frac{x}{x'}\big)^{\frac{n-1}{2}} = c_n x^{(n-1)/2} h^{n/2} e^{-iy' \cdot z/h} \nu.
$$ This gives orders $m = n/2 - (2n)/4 = 0$ at the main face and $(n-1)/2 +
1/2 - (2n)/4 = 0$ at the b-face. The zero orders reflect the unitarity of
this operator. A possibly more natural way of writing the kernel is
$$
P(\lambda) = c_n e^{-i\lambda y' \cdot z} |\lambda^{n-1} d\lambda dy' dz|^{1/2},
$$
in which it is clear that $P(\lambda)$ is essentially the Fourier transform. 

To get the scattering matrix we again divide by $|dr|^{-1/2}$ and restrict at $r=\infty$. Let $\nu'$
be a scattering-fibred half-density on the scattering matrix space $\pa M
\times \pa M \times [0, h_0).$ Then
\begin{align*} \Big| \frac{dh}{h^2} \frac{dx}{x^2 h} \frac{dy}{(xh)^{(n-1)}}
\frac{dy'}{(xh)^{(n-1)}} \Big|^{1/2} &= \Big| \frac{dh}{h^2}
\frac{dy}{h^{(n-1)}} \frac{dy'}{h^{(n-1)}} \Big|^{1/2} |dr|^{1/2}
h^{-1/2} |z|^{n-1} \\ 
&= \nu' |dr|^{1/2}
h^{-1/2} |z|^{n-1}.
\end{align*}
Thus,
$$S(h^{-1}) =  \lim_{|z| \to \infty} c_n e^{-iy' \cdot z/h} |z|^{-(n-1)/2} h^{(n-1)/2} \nu'
|z|^{n-1}.
$$ This can be written
$$ S(h^{-1})= \lim_{r \to \infty} c_n e^{-iry' \cdot y/h} \big( \frac{r}{h}
\big)^{(n-1)/2} \big| \frac{dh}{h^2} dy dy' \big|^{1/2} = \delta(y - y')
\big| \frac{dh}{h^2} dy dy' \big|^{1/2}
$$ which is the scattering matrix times a scattering half-density in
$h$. We may also write
$$S(h^{-1}) = \int e^{i(y-y') \cdot \eta/h} d\eta \Big| \frac{dh dy dy'}{h^2 h^{n-1}
h^{n-1}} \Big|^{1/2}
$$ which has semiclassical order $m = 0 + (n-1)/2 - (2n-1)/4 = -1/4$ and
Lagrangian order $-1/4$. This implies that the order as a Lagrangian for a fixed positive $h$ is $0$ (see Section~\ref{llds}), so this again reflects unitarity of $S(h^{-1})$ for a fixed $h$. 

The free propagator is given, using the convention in \cite{HW} regarding the half-density factor in $t$,  by
$$
(2\pi i t)^{-n/2} e^{i|z-z'|^2/2t} \big| dz dz' dt \big|^{1/2}.
$$
We can write $t^{-n/2} |dz dz' dt|^{1/2}$ in the form $t^{n/2 + 1} \mu$,
where $\mu$ is a scattering fibred half-density, and this in turn can be written 
$$
t^{n/2 + 1} \big( \rho_{\lb} \rho_{\rb} \rho_{\bfc} \big)^{N/2} \mu_q,
$$
where $\mu_q$ is a nonvanishing quadratic scattering half-density. It follows from \eqref{eq:dist-1} that the orders of this distribution are $n/2 + 1 - (2n+1)/4 = 3/4$ at $Q(L)$, $(2n+1)/2 + (n+1)/2 - 3(2n+1)/4 = 1/4$ at $\lb$ and $\rb$ and $(2n+1)/2 + 1/2 - 3(2n+1)/4 = -n/2 + 1/4$ at $\bfc$. As with the resolvent, we may regard this as an especially simple case of a Legendrian associated to a pair of intersecting Legendre distributions $(Q(L), Q(L_2^\sharp))$ with conic points, where the distribution is in fact `smooth across $Q(L_2^\sharp)$'. We then find that the free propagator is an element of
$$
I^{3/4, n/2 + 1/4; 1/4, -n/2 + 1/4}(X,(Q(L), Q(L_2^\sharp)); \sfOh).
$$
These orders agree with those calculated in \cite{HW}. 
\ah{put into intro: we can verify with bare hands that these kernels
obey the theorem. In fact there is no conic singularity in this case. This is almost unique to the free geometry, so in this respect these kernels are very misleading as to the general case.}

%%%%%%%%%%%%%%%%%%%%%%%%%%%%%%%%%%%

\section{Pseudodifferential construction}\ilabel{sec:hpseudos}
In this section we show that the inverse $(h^2 \Delta + V - \evosq)^{-1}$
of our operator family lies in the algebra of \emph{semiclassical
pseudodifferential operators} when $\Re \evosr \neq 0$.

\subsection{$h$-pseudodifferential calculus}
The scattering calculus as described here was introduced by Melrose
\cite{Melrose:Spectral}, although its roots go back a good deal further: in
various guises on $\RR^n$ it has been examined by Shubin \cite{shubin},
Parenti \cite{parenti}, and Cordes \cite{cordes}; on manifolds, it has also
been considered by Schrohe \cite{schrohe:calculus}.  It is also the Weyl
calculus for the metric
$$
\frac{\abs{dz}^2}{1+\abs{z}^2}+\frac{\abs{d\xi}^2}{1+\abs{\xi}^2}.
$$ The semiclassical variant has been considered by Vasy-Zworski \cite{VZ}
and by the second author and Zworski \cite{Wunsch-Zworski}. See the
appendix of \cite{Wunsch-Zworski} for a summary of the properties of this
class of operators. Here we simply recall that the space
$\Psisch^{m,l,k}(X)$ of semiclassical scattering pseudos is indexed by the
differential order $m$, the boundary order $l$ and the semiclassical order
$k$. This space of operators can be expressed in terms of the space
$\Psisch^{m,0,0}(X)$ by
\begin{equation}
\Psisch^{m,l,k}(X) = x^l h^{-k} \Psisch^{m,0,0}(X).
\end{equation}
Following \cite{Wunsch-Zworski}, we shall restrict to operators with
polyhomogeneous symbols.  The symbols of such operators are functions $a$
on $(0,h_0) \times \Tbarscstar X$, having the property that $h^k x^{-l}
\rho^{-m} a \in \CI([0,h_0) \times \Tbarscstar X$, where $\rho$ is a
boundary defining function for the boundary hypersurface of $\Tbarscstar X$
at fibre-infinity. The (principal) symbol map is given by restriction of
$h^k x^{-l} \rho^{-m} a$ to the boundary of $[0,h_0) \times \Tbarscstar X)$
and denoted $\sigmasch^{m,l,k}(A)$, where $a$ is the symbol of $A$. Since
the boundary consists of three different faces, one at $\rho = 0$, one at
$x=0,$ and one at $h=0$, the principal symbol corresponding can be
decomposed into three parts (subject, of course, to compatibility
conditions where the faces intersect).  We shall call these parts the
interior symbol, the boundary symbol and the h-symbol respectively.  These
symbols lead to three separate exact sequences, with each symbol being the
obstruction to the operator being of lower order in the corresponding
sense: if the $h$-symbol vanishes, our operator is divisible by an
additional power of $h;$ if the $x$ symbol vanishes, by a power of $x;$ and
if the $\rho$ symbol vanishes, the operator is of lower order in the
(usual) sense of differentiation.

\ah{How much should be put here?}

\subsection{Resolvent (away from the spectrum)}
Let $\evosr$ have nonzero imaginary part. Then the principal symbol of $
h^2 \Delta + V - \evosq$ is equal to $g(z, \xi) + V(z) - \evosq$. This is
invertible on each boundary face, so by the symbol calculus there is an
operator $G_1(\evosr) \in \Psisch^{-2, 0, 0}(X)$ such that
$$(h^2 \Delta + V - \evosq) G_1(\evosr) = \Id - E(\evosr), \ E(\evosr) \in
\Psisch^{-1, 1, -1}.$$ 
Let $E_2(\evosr)$ be an asymptotic sum of the Neumann series $\Id + E(\evosr) + E(\evosr)^2 + \dots$. Then we have, with $G_2(\evosr) = G_1(\evosr) E_2(\evosr)$, 
$$(h^2 \Delta + V - \evosq) G_2(\evosr) = \Id - E_\infty(\evosr), $$ 
with $E_\infty(\evosr)$ in the `completely residual space'  $\Psisch^{-\infty, \infty, -\infty}$; equivalently, the kernel of $E_\infty(\evosr)$ is in $h^\infty \rho^\infty C^\infty(M^2)$, where $\rho$ is a product of boundary defining functions for $M^2$. 
The inverse of $\Id - E_\infty(\evosr)$ certainly exists as a bounded operator on $L^2(M)$, for small $h$, since the operator norm $\| E_\infty(\evosr) \|_{L^2 \to L^2}$ is $O(h^\infty)$. 
Let us write $(\Id - E_\infty(\evosr))^{-1}$ = $\Id + S(\evosr)$. We then have
$$
S(\evosr) = E_\infty(\evosr) + E_\infty(\evosr)^2 + E_\infty(\evosr) S(\evosr) E_\infty(\evosr).
$$
This identity shows that the kernel of $S(\evosr)$ is also in $h^\infty \rho^\infty C^\infty(M^2)$. Thus, we have $S(\evosr) \in \Psisch^{-\infty, \infty, -\infty}$. The resolvent is equal to 
$$
(h^2 \Delta + V - \evosq)^{-1} = G_2(\evosr) (\Id + S(\evosr))
$$
which is in  $\Psisch^{-2,0,0}$, as claimed.

%%%%%%%%%%%%%%%%%%%%%%%%%%%%%%%%%%%%%%%%%%%

\section{Structure of the propagating Legendrian}\ilabel{sec:prop}

We now consider the case where $\evo$ is real and positive, i.e.\ we are
on the spectrum. In the case where $\evo$ is in the resolvent set,
studied in the previous section, the singularities of $(H -
\evosq)^{-1}$ live on the conormal bundle of the diagonal. Here, by
constrast, singularities propagate off the diagonal. The reason for this is
that the characteristic variety of the operator $H - \evosq$, in either
the left or the right variable, intersects the conormal bundle of the
diagonal on $\mf$ (as well as $\bfc$). Moreover, the Hamilton vector field along the characteristic set is nonzero at
this intersection, which allows singularities to move into the
characteristic set away from the diagonal. In this section we analyze
the geometric structure of this flowout, along bicharacteristics, from the characteristic variety at the diagonal; we shall see that it
forms a Legendre submanifold in $\Tsfstar[\mf] (\XXbh)$ which
becomes smooth after certain blowups are performed.

The first step is to compute the left and right Hamilton vector fields for
the operator $H - \evosq$. First, we do this in the interior of
$\Tsfstar[\mf] \XXbh$. We may choose coordinates $z, z', \zeta, \zeta',
\tau$, corresponding to writing covectors
$$ \tau \cdot d \big( \frac1{h} \big) + \zeta \cdot \frac{dz}{h} + \zeta'
\cdot \frac{dz'}{h}.
$$ The left and right Hamilton vector fields take the form (where we divide by a factor of $2$ for convenience)
\begin{equation}\ilabel{Vl1}
V_l = h \Big( g^{ij}(z) \zeta_i \dbyd{}{z_j} -  \frac1{2} \big(
\frac{\partial g^{ij}(z)}{\partial z_k} \zeta_i \zeta_j + \frac{\partial
V}{\partial z_k} \big) \dbyd{}{\zeta_k} + g^{ij}(z) \zeta_i \zeta_j
\dbyd{}{\tau} \Big)
\end{equation}
and
\begin{equation}\ilabel{Vr1}
V_r = h \Big( g^{ij}(z') \zeta'_i \dbyd{}{(z')^j} -  \frac1{2} \big(
\frac{\partial g^{ij}(z')}{\partial (z')^k} \zeta'_i \zeta'_j +
\frac{\partial V}{\partial (z')^k} \big) \dbyd{}{\zeta'_k} + g^{ij}(z')
\zeta'_i \zeta'_j \dbyd{}{\tau} \Big).
\end{equation}
Let us write $V_l' = V_l/h$ and $V'r' = V_r/h$, restricted to $\{ h = 0
\}$.  These vector fields commute, and are tangent to the left and right
characteristic sets
\begin{equation}\ilabel{charsets} \Sigma_l = \{ g^{ij}(z)\zeta_i \zeta_j + V(z) = \evosq \}, \quad \Sigma_r =
\{ g^{ij}(z')\zeta'_i \zeta'_j + V(z') = \evosq \}.
\end{equation} Let $\Sigma=\Sigma_L\cap \Sigma_R$ denote the intersection
of the characteristic sets for $p_L$ and $p_R$, and let
\begin{equation}
N^* \Delta_b = \{ d \big( \frac{f}{\xx} \big)(p) \mid p \in \Delta_b, f \in \CIsf(X), f \restrictedto \Delta_b = 0 \};
\label{ndb}\end{equation}
in coordinates, $N^* \Delta_b = \{ h = 0, z = z', \zeta = -\zeta', \tau = 0 \}$. 
 Note that on
$N^*\Delta_b,$ $\Sigma_L$ and $\Sigma_R$ coincide; hence $N^* \Delta_b \cap
\Sigma_L = N^* \Delta_b \cap \Sigma_R = N^* \Delta_b \cap \Sigma$ and is
codimension $1$ in $N^* \Delta_b$.

Notice also that $V_l'$ and $V_r'$ are everywhere nontangential to $N^*
\Diag_b$. In fact, for $V_l'$ to be tangential we would need $\zeta = 0$
and $\nabla V = 0$, which means that $V_l' = 0$; but this contradicts the nontrapping
hypothesis.  Consider the flowout by $V'_l$ from the intersection of
$N^*\Delta_b \cap \Sigma_L $. It is at least locally a smooth Legendre
manifold (Legendre because the vector fields $V_l'$ and $V_r'$ are contact
vector fields and the initial hypersurface is isotropic of dimension
$2n-1$).  However, $V_l' - V_r'$ is tangent to $N^*\Delta_b$. Moreover,
$[V_l, V_r] = 0$, as follows directly from the commutation of the left and
right operators $h^2 \Delta_l$ and $h^2 \Delta_r.$ The two-plane
distribution spanned by $V_l,$ $V_r$ (or $V_l',$ $V_r'$) is therefore
integrable; as $V_l'-V_r'$ is tangent to $N^*\Delta_b,$ the integral
manifold consisting of all leaves through $N^* \Delta_b$ is thus
$2n$-dimensional (rather than $2n+1$-dimensional as one would expect
without this tangency).  It follows that the flowout from $N^*\Delta_b \cap
\Sigma$ by $V_l'$ coincides with the flowout by $V_r'.$

This geometry holds uniformly to the boundary of $N^* \Delta_b$. We now
work near $N^* \Delta_b \cap \bface$. Then we use coordinates (where
$\theta = x'/x$ is small, so $x' \ll x$)
\begin{equation}\ilabel{canonicaloneform}
\lambda' d\big( \frac1{x \theta h} \big) + \lambda d\big( \frac1{x h} \big)
+ \tau d\big( \frac1{h} \big) + \mu' \frac{dy'}{x \theta h} + \mu
\frac{dy}{x h} .
\end{equation}
In fact these coordinates are valid in the region $\theta \leq C$ for any
finite $C$, say $C=2$, a region which includes a neighbourhood of the
corner $\bface \cap \rb$.
 
The symbols of $h^2\Lap+V-\evosq$ acting on the left and right factors are respectively
$$ p_L = \lambda^2+h^{ij}(x,y)\mu_i \mu_j + V(x,y) -\evosq,
$$
$$ p_R = (\lambda')^2+h^{ij}(x',y')\mu_i' \mu_j' + V(x',y') -\evosq.
$$ The left and right vector fields thus take the form
\begin{multline}
V_l = xh \Big( -\lambda x \partial_x + \lambda \theta \partial_{\theta} +
h^{ij} \mu_i \partial_{y_j} + \big(h^{ij} \mu_i \mu_j + \frac{1}{2} x \partial_x
(h^{ij} \mu_i \mu_j  +  V)\big) \partial_{\lambda} \\ + \big(- \mu_k \lambda - \frac1{2}  \partial_{y_k} (h^{ij} \mu_i \mu_j +  V)\big) \partial_{\mu_k}  - \frac1{2} \partial_{x}
h^{ij} \mu_i \mu_j \pa_\tau \Big), \ilabel{Vl}\end{multline} and\ah{Potential term missing from these vector fields}
\begin{multline}\label{rightvf}
V_r = x\theta h \Big( -\lambda' \theta \partial_{\theta} + h^{ij}(x\theta, y') \mu'_i
\partial_{y'_j} + \big(h^{ij}(x\theta, y') \mu'_i \mu'_j + \frac{1}{2} \theta \partial_\theta ( h^{ij}(x\theta, y')
\mu'_i \mu'_j  + V(x\theta, y'))\big) \partial_{\lambda'}  \\ + \big(- \mu'_k \lambda' -\frac1{2} \partial_{y'_k} ( h^{ij}(x\theta, y')
\mu'_i \mu'_j + V(x\theta, y')) \big) \partial_{\mu'_k}  - \frac1{2}
(\partial_{x'} h^{ij})(x\theta, y') \mu'_i \mu'_j \pa_\tau \Big).\end{multline} In
this region let us write $V_l' = V_l/xh$ and $V_r' = V_r/x\theta h$. Then
we have
\begin{multline}\label{vfsum}
\frac{V_l}{xh} + \frac{V_r}{x\theta h} = (\lambda - \lambda') \theta
\pa_\theta -\h \left( \pa_x \abs{\mu}^2+ \pa_{x'} \abs{\mu'}^2\right)\pa_\tau
-\lambda x \pa_x  +h^{ij} \mu_i \partial_{y_j} \\ + \big(h^{ij} \mu_i \mu_j + \frac{1}{2} x \partial_x
(h^{ij} \mu_i \mu_j  +  V)\big) \partial_{\lambda}  + \big(- \mu_k \lambda - \frac1{2}  \partial_{y_k} (h^{ij} \mu_i \mu_j +  V)\big) \partial_{\mu_k} \\ 
+ {h'}^{ij} \mu'_i \partial_{y'_j}  + \big({h'}^{ij} \mu'_i \mu'_j + \frac{1}{2} \theta \partial_\theta ( {h'}^{ij}
\mu'_i \mu'_j  + V' )\big) \partial_{\lambda'}   + \big(- \mu'_k \lambda' -\frac1{2} \partial_{y'_k} ( {h'}^{ij}
\mu'_i \mu'_j + V') \big) \partial_{\mu'_k}
\end{multline}
Note that this vanishes only on $\{\lambda'=\lambda,\ \mu=\mu'=0,\ x=0\}$.

We define the sets $L_+, L_-$ and $L$ by
\begin{equation}\begin{gathered}
L_+ / L_- \text{ is the forward/backward flowout from } N^*\Delta_b \cap
\Sigma \text{ by } V_l' \\ L = L_+ \cup L_- .
\end{gathered}\ilabel{Ldefn}\end{equation}
Equivalently, we may define $L_+ / L_-$ as the forward/backward flowout
from $N^*\Delta_b \cap \Sigma$ by $V_r'$. By the arguments above, $L_\pm$
and $L$ are Legendrian submanifolds; moreover, the pairs $(N^* \Delta_b,
L_\pm)$ form intersecting pairs of Legendre submanifolds in the sense of
Section~\ref{section:ild}.

The main goal of this section is to determine the regularity of the
Legendrian $L$, which we call the `propagating Legendrian', as we
move far from $N^* \Delta_b$. By symmetry it suffices to consider just
$L_+$. It turns out that $L_+$ is smooth except for a conic singularity at
a submanifold $\Lsharp_2$ of $\Tsfstar[\bface] \XXbh$; when $J$, the span
of $\Lsharp_2$, is blown up, $L_+$ lifts to a smooth manifold with
codimension 3 corners.

First consider smoothness at $\bface = \{x = 0 \}$. Notice that the flows
$V_l'$, $V_r',$ when restricted to $\bface,$ are naturally identified with
the flows for a fixed positive value of $h$ from \cite{hasvas2}, so $L_+ \cap
\bface$\ah{inaccurate} can be identified with $L_+(\lambda)$ from
\cite{hasvas2}.  It was shown in \cite{hasvas2} that $L_+(\lambda)$ was smooth
after the space $\{ x=0, \lambda = \lambda', \mu = \mu' = 0 \}$ was blown
up\footnote{The coordinates $\lambda, \lambda'$ here correspond to $\tau',
\tau''$ from \cite{hasvas2}.}. Let us then define
$$ J = \{ x=0, \lambda = \lambda', \mu = \mu' = 0 \} = \spn \Gsharp_2 \subset \Tsfstar \XXb$$
with 
$$
\Gsharp_2 =  \{
x=0, \lambda = \lambda' = 1, \mu = \mu' = 0 \}
$$ and consider the space
\begin{equation}
[\Tsfstar \XXbh; J].  \ilabel{blowup-XXbh}\end{equation} We denote by $\Jt$
the lift of $J$ to this space, i.e.\ the new boundary hypersurface created
by the blowup.

\begin{proposition}\ilabel{prop:leaf} The closure of the lift of $L_+$ to
  the space \eqref{blowup-XXbh} is a smooth manifold with corners of
  codimension three. Consequently, in a neighbourhood of $\Gsharp_2$, the
  pair $(L_+, \Gsharp_2)$ forms a conic Legendrian pair of submanifolds in
  the sense of Section~\ref{Lconicpts}.  
\end{proposition}

\begin{proof} It suffices to show that $L = L_+ \cup L_-$ is a smooth
  manifold with corners of codimension three, since $L$ is transversal to $N^* \Delta_b$, which divides it smoothly into  two pieces $L_+$ and $L_-$. By standard ODE theory, $L$ is smooth at all points reachable from $N^* \Delta_b$ by the vector field $V_l'$ or $V_r'$ in a finite time. However,
we need to check the regularity of the closure of $L$ at the boundary of $\Tsfstar \XXbh$.

It has already been observed that the two-plane distribution $D$ determined
by $V_l'$ and $V_r'$ is integrable; therefore $L$ is foliated by
two-dimensional leaves, each of which intersects $N^* \Delta_b$ in a
one-dimensional set (since $V_l' - V_r'$ is tangent to $N^*
\Delta_b$). Consider a point $(q, \tilde q) \in N^* \Delta_b \cap \Sigma$,
where $q$ is a covector in the interior of $\Tscstar[z] X$ with $|q| =
\sqrt{\evosq - V(z)}$, and the tilde denotes negation of the fibre variables of
$q$. The leaf containing this point is the set of points $(q_1, \tilde
q_2)$, where $q_i$ lie on the same bicharacteristic $\gamma$ as $q$; we
shall denote it $\gamma^2 = \gamma_q^2$.

It is convenient to choose a `section' $S$ of $N^* \Delta_b$, by which we
mean a smooth submanifold of $N^* \Delta_b$ of codimension 1 that
intersects each $\gamma^2$ transversely at a unique point.  It is not
difficult to see that a section exists, using the following argument: From
\eqref{Vl}, under the flow along $V_l/(xh),$
$$
x' = -\lambda x, \quad \lambda' = h^{ij}\mu_i\mu_j + \half x
\partial_x h^{ij} \mu_i \mu_j + \half x
\partial_x V(x,y).
$$ The form $h^{ij} + \half x \partial_x h^{ij}$ is positive definite for
small $x$ and the potential term vanishes at $x=0,$ hence choosing $x_1>0$
and $\delta>0$ small, for $x \leq x_1$ and $\abs{\lambda}<\evosr-\delta,$
we have $\lambda'>0$ on the characteristic variety.  For small $x$, we can
now take $S$ to be $N^* \Delta_b \cap \{ \lambda = 0 \}$. The null
bicharacteristics corresponding to such points remain in the region where
$x$ is small, since on the flowout of $S\cap\Sigma,$ it is easy to verify
that $x'<0$ and $\lambda'>0$ except possibly when $\lambda>\evosr-\delta.$
In the region where $x$ is large---say $x \geq x_0$ where $x_0 <
x_1$---each entering geodesic meets the boundary $\{ x = x_0 \}$ in exactly
two points (by the same argument as above). We can take the section to be
that point on the diagonal of $\gamma^2$ corresponding to the point on the
geodesic which is halfway (with respect to arc length) between the two
intersection points with $\{ x = x_0 \}$. We interpolate between these two
prescriptions to obtain a smooth section $S$.  Then each leaf intersects
$S$ in a unique point.

The strategy of our proof is to first restrict attention to a single leaf
and analyze its closure; we shall show that it is a manifold with
codimension 2 corners. We shall then show that the union of these closed
leaves is the closure of $L$, and that this forms a submanifold with
codimension 3 corners.

We will, initially, have to work on a larger (i.e.\ more blown-up) space
than \eqref{blowup-XXbh}. Let $J_-$ be the submanifold
$$ \{ x = 0, \ \mu = 0 , \ \mu' = 0, \ \lambda = -\lambda' \} \subset
\Tsfstar[\bfc] \XXbh
$$ (the only difference in this definition and that of $J$ being a change
of sign in the equation $\lambda = \pm \lambda'$). We shall blow up the
submanifold $J_-$ as well as \footnote{Note that although these two
submanifolds intersect, the intersection is away from the closure of $L$,
since on $Z \cap Z_-$, $\mu = \mu' = 0$, $\lambda = \lambda' = 0$; on the
other hand,  $L$ is contained in $\Sigma_l \cap \Sigma_r$, so over $\bfc$ this is given by  $\lambda^2 +
|\mu|^2 = (\lambda')^2 + |\mu'|^2 = \evosq > 0$. We are only interested in a
neighbourhood of the closure of $L$, so $Z$ and $Z_-$ are disjoint in the
region of interest, hence they can be blown up independently.} $J$. Also
consider the submanifold
$$ W = \{ \theta = 0, \mu' = 0 \} \subset \Tsfstar[\rb](\XXbh).
$$ After $J \cup J_-$ is blown up, $W$ lifts to a new submanifold $W'$
which is transverse to $\Jt$ and $\tilde J_-$. Consider the space
\begin{equation}
\big[ [\Tsfstar \XXbh; J \cup J_-]; W' \big].
\ilabel{double-blowup}\end{equation} Denote the new boundary hypersurface
created by this blowup by $\Wt$. We shall work on the space
\eqref{double-blowup} for most of this proof, although eventually we shall
see that we can return to the space \eqref{blowup-XXbh}.

Consider a leaf $\gamma^2$ of the distribution $D$ which intersects $S$
at $(q,q)$, where $q$ lies in the interior of $\XXsc$. (Later we consider
$q$ lying in the boundary, i.e.\ at $x=0$.)  Let $y_{-\infty}$,
resp. $y_\infty$ be the points on $\partial X$ obtained as the
initial, resp.\ final end of the bicharacteristic through $q$.  Consider the
intersection of $\gtcl$ with the boundary of $\Tsfstar \XXbh$, i.e.\
with $\{ x = 0 \} \cup \{ x' = 0 \}$.  To get there we must send either $q_1$
or $q_2$ to infinity along the bicharacteristic. If we send $q_1$ to
infinity keeping $q_2$ fixed, we arrive at the set
$$ \{ (y_{-\infty}, 0, -1, z', \zeta', 0) \mid (z', \zeta') \in \gamma \}
\cup \{ (y_{\infty}, 0, 1, z', \zeta', 0) \mid (z', \zeta') \in \gamma \}
\subset \Tsfstar[\lb]\XXbh
$$ using coordinates $(y, \mu, \lambda, z', \zeta,
\theta^{-1})$. Similarly, if we send $q_2$ to infinity keeping $q_1$ fixed,
we arrive at the set
$$ \{ (z, \zeta, y_{-\infty}, 0, 1, 0) \mid (z, \zeta) \in \gamma \} \cup
\{ (z, \zeta, y_{\infty}, 0, -1, 0) \mid (z, \zeta) \in \gamma \} \subset
\Tsfstar[\rb]\XXbh
$$ using coordinates $(z, \zeta, y', \mu', \lambda',\theta)$. If $q_1$ and
$q_2$ are simultaneously sent to infinity, we end up at the set
$$ \{ (y_{-\infty}, 0, -1, y_{-\infty}, 0, 1, \theta ) \} \cup \{
(y_{\infty}, 0, 1, y_{\infty}, 0, -1, \theta ) \}
$$ if they go to infinity in the same direction, or
$$ \{ (y_{-\infty}, 0, -1, y_{\infty}, 0, -1, \theta ) \} \cup \{
(y_{-\infty}, 0, -1, y_{\infty}, 0, -1, \theta ) \}
$$ if they go to infinity in opposite directions.

We claim that the closure $\gtcl$ inside the space \eqref{double-blowup} is
a surface with corners, with eight edges as above, as in
Figure~\ref{fig:leaf}. Our analysis is based on the following lemma. Before stating this we need the following

\begin{defn} Let $X$ be on a manifold with corners, and let $\Nub(X)$ denote the smooth vector fields on $X$ tangent to each boundary hypersurface.  Let $\rho$ be a boundary defining function for a boundary hypersurface $H$ of $X$. 
We say that $V \in \Nub(X)$ is b-normal at $H$ if 
$$
V = c \rho \partial_\rho + \rho W \text{ for some } W \in \Nub(X)
$$
where the coefficient $c$ is never zero. We say that $V$ is incoming, resp. outgoing b-normal if c is positive, resp. negative. (We note that the vector field $\rho \partial_\rho \restrictedto H$ is  nonzero as a b-vector field, and independent of coordinates.)
\end{defn}

Notice that if a vector field $V$ is b-normal at $H$, then $V/\rho$ is smooth and transverse to $H$.

\begin{comment}
\begin{figure}
\input{leaf.eepic}
\end{figure}
\end{comment}
\begin{figure}\centering
\epsfig{file=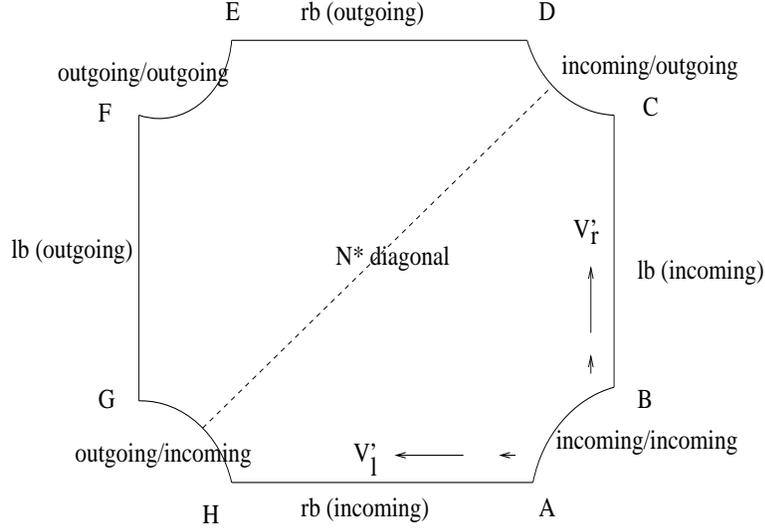,width=10cm,height=7cm}
\caption{The closure of a leaf. Here `incoming', resp. `outgoing' in the
left factor means at $\mu = 0, \lambda <0$, resp. $\lambda >0$, while
for the right factor it means $\mu' = 0, \lambda' <0$, resp. $\lambda' >0$.}
\label{fig:leaf}
\end{figure}

\begin{lemma}\ilabel{b}
On the space \eqref{double-blowup}, the vector field $V_r'$ is
incoming/outgoing b-normal at $\Wt \cap \Sigma_r$, $V_l' + V_r'$ is
incoming/outgoing b-normal at $\Jt \cap \Sigma_l \cap \Sigma_r$ and $V_l' -
V_r'$ is incoming/outgoing b-normal at $\tilde J_- \cap \Sigma_l \cap
\Sigma_r$. In all cases, the sign of $\lambda, \lambda' \in \{ \pm \evo \}$
determines whether the vector field is incoming or outgoing.
\end{lemma}

\begin{proof}
We first look at $V_l' + V_r'$. Since $L$ is contained in both the left and
right characteristic varieties, we have
\begin{equation}
\lambda^2 - (\lambda')^2 = |\mu'|^2 - |\mu|^2 \text{ on } L.
\ilabel{l}\end{equation} Thus $\lambda - \lambda' = o(|\mu'| + |\mu|)$ near
$p$, so we can take a boundary defining function for $\Jt$ in $L$ to be
$\rho_{\Jt} = \sqrt{x^2 + |\mu|^2 + |\mu'|^2}$.  By \eqref{l}, $\lambda -
\lambda'$ is $O(\rho_{\Jt}^2)$, so \eqref{vfsum} gives
$$ (V_l' + V_r')|_{L} = -\lambda (x \partial_x + \mu \cdot \partial_{\mu} +
\mu' \cdot \partial_{\mu'} + (\lambda - \lambda') (\partial_\lambda -
\partial_{\lambda'})) + \mu \cdot \partial_y + \mu' \cdot \partial_{y'} +
O(|\mu|^2) (\partial_\tau, \partial_\theta, \partial_\lambda, \partial_\mu,
\partial_{\lambda'}, \partial_{\mu'}).
$$ This implies that in our local coordinates,
\begin{equation}
V_l' + V_r' = -\lambda \rho_{\Jt} \partial_{\rho_{\Jt}} + \rho_{\Jt} \Nub,
\ilabel{V-prop}\end{equation} on the space \eqref{blowup-XXbh}. An
analogous argument applies to $V_l' - V_r'$ at the blowup of $J_-$.

We next analyze $V_r'$. In \eqref{double-blowup}, the submanifold $\Wt$ is
given by the equations $\{ \theta= 0, \mu'/\rho_{\Jt} = 0 \}$. Hence in a
neighbourhood of $\Wt$, $\mu' = o(\sqrt{x^2 + |\mu|^2})$, so we may take
$\rho_{\Jt} = \sqrt{x^2 + |\mu|^2}$ in this region, which we shall do from
now on.  We can switch to local coordinates on \eqref{blowup-XXbh}
\begin{equation}
y, \ y', \ \rho_{\Jt} = \sqrt{x^2 + |\mu|^2}, \ \rho_{\bfc} = \frac{x}{x +
\rho_{\Jt}}, \ \theta, \ M' = \frac{\mu'}{\rho_{\Jt}}, \ \lambda, \ \Lambda
= \frac{\lambda - \lambda'}{\rho_{\Jt}}, \ \hat \mu .
\ilabel{coords-1-1}\end{equation} From \eqref{rightvf}, in these coordinates,
$$ V_r' = -\lambda (\theta \partial_{\theta} + M' \cdot \partial_{M'}) +
\mu' \cdot \partial_{y'} + O(|\mu'|^2)(\partial_{\lambda}, \partial_{\mu'},
\partial_\tau),
$$ hence \ah{Need to think about whether this
is a complete proof}
$$ V_r' = -\lambda' \rho_{\Wt} \partial_{\rho_{\Wt}} + \rho_{\Wt} \Nub.
$$
\end{proof}

\emph{Continuation of the proof of Proposition~\ref{prop:leaf}.}  Now we
return to showing that the closure of $\gamma^2$ is a smooth 2-manifold
with corners. First consider the point $A$ in the figure. This lies on the
intersection of $\Jt$ and $\Wt$. We may set $\Vr = V_r'/\theta$ and $\Vc =
(V_l' + V_r')/\rho_{\Jt}$; then $\gamma^2$ is contained in the flowout from
$\gamma^2 \cap S$ by $\Vr$ and $\Vc$. Notice that these vector fields no
longer commute, but they still determine an integrable two-plane
distribution $D$. By Lemma~\ref{b} and the remarks above it, they are both
smooth vector fields on \eqref{double-blowup} such that $\Vr$ is
transversal to $\Wt \cap \Sigma_r$ and tangent to $\Jt$, while $\Vc$ is
transversal to $\Jt \cap \Sigma_l \cap \Sigma_r$ and tangent to $\Wt$. The
flowout by these vector fields therefore sweeps out a smooth, closed
2-dimensional manifold with corners meeting the boundary of
\eqref{double-blowup} transversally, and it is clear that this is the
closure of the leaf.

Since $L$ is invariant under the flow of $\Vr$, which is tangent to the
lift of $\rb$ and $\bfc$, the closure of the leaf is a smooth submanifold
which is disjoint from $\rb$ and $\bfc$ (assuming that $\gamma^2$ is a
leaf through $S^\circ$, the interior of the section $S$). It follows that $d\theta \neq 0$ at the
intersection of $L$ and $\Wt$, since $\theta$ can be taken as a boundary
defining function for $\Wt$ away from $\rb$.

Since $\Vc$ and $\Vr$ do not vanish in a neighbourhood of $A$, nearby
leaves also have this property, and they vary smoothly with their
intersection point $\ss \in S$ by standard ODE theory. This gives us smooth
coordinates on the closure of $L$ near the point $A$, namely $\theta$,
$\rho_{\Jt}$, and a coordinate on $S$.

Exactly the same argument gives smoothness near the point $E$. Indeed, a
similar argument applies to the corner points $D$ and $H$ in
Figure~\ref{fig:leaf}, since there is a symmetry of $L$ coming from the
involution $(q, q') \mapsto (q' q)$ on $\XXsc$. Moreover, essentially the
same argument also gives smoothness near the other corner points; the only
difference is that we are working near the blowup of $J_-$ rather than $J$,
but this makes no difference at all, because if we replace the minus sign
by a plus sign in the left hand side of \eqref{l} this makes no difference
to the argument.

We also need to check smoothness near a point on an edge. However, we have
effectively already done this, because our coordinates are valid for
$\theta \leq 2$, say, while for $\theta \geq 1/2$ we can perform the
involution above and use $V_l'$ instead of $V_r'$.

Notice that the closure of this leaf is disjoint from $\bfc$ (or more
precisely, the lift of $\bfc$ to \eqref{double-blowup}). In fact the vector
fields $\Vc$ and $\Vr$ are everywhere tangent to $\bfc$ so it is impossible
to reach $\bfc$ after flowing for a finite time along these vector
fields. There is another way of seeing this which gives more insight into
how these leaves fit together.  Notice that a boundary defining function
for $\bfc$ on the space \eqref{double-blowup} can be taken to be
$$ \frac{x+x'}{\sqrt{x^2 + (x')^2 + |\mu|^2 + |\mu'|^2}}
$$ in a neighbourhood of $L$. For an exactly conic metric, the quantity
$x/|\mu|$ is constant along the bicharacteristic and is equal to the
maximal value of $x$ that occurs along it. In a general scattering metric,
this quantity is approximately equal to the maximal value of $x$ along the
bicharacteristic, and this approximation is better and better (in the sense
that the ratio of these two quantities tends to 1 uniformly as the
bicharacteristic approaches the boundary uniformly); this follows from
\cite{HTW} for example. Hence, for a fixed interior leaf, the limiting
value of $x/|\mu|$ is nonzero, which says that the leaf is disjoint from
$\bfc$. On the other hand, the leaf will approach $\bfc$ uniformly as the
associated bicharacteristic approaches the boundary uniformly.

Now consider a leaf associated to a boundary point of $S$. In that case,
the bicharacteristic is a limiting geodesic contained in the boundary of
$X$, so the leaf is contained in $\bfc$. In this case, there is an explicit
formula for the leaf. Fix $(y_0, \mu_0) \in T_{y_0} \partial M$ with $|\mu_0| \leq 1$. Then the
leaf is given by
\begin{equation}\begin{aligned}\ilabel{eq:sp-1c}
{} &\{(\theta,y,y',\lambda,\lambda',\mu,\mu'):\
\exists(y_0,\hat\mu_0)\in S^* \partial X, \ s,s'\in(0,\pi),\ \text{s.t.}\\
&\quad \theta = \frac{\sin s'}{\sin s}, \ \lambda=-\abs{\evosr} \cos s, \ \lambda'=\abs{\evosr}\cos
s',\\ \quad(y,\mu)&=\abs{\evosr}\sin s\exp(sH_{\half h})(y_0,\hat\mu_0),
(y',\mu')=-\abs{\evosr}\sin s' \exp(s' H_{\half h})(y_0,\hat\mu_0)\}
\end{aligned}\end{equation}
This corresponds to the interior of the leaf in Figure~\ref{fig:leaf},
which we can think of in this boundary case as the square $(0, \pi)^2$ with
the $s$ axis horizontal and the $s'$ axis vertical; $V_l'$ is given by
$-\sin s \partial_s$ and $V_r'$ is given by $\sin s' \partial_{s'}$ in
these coordinates. The closure is given by the closed region in the figure,
where the parts over $\rb$ (i.e, the boundary lines AH and DE) now lie over
the intersection of $\rb \cap \bfc$ and the parts at $\Jt$ and $\tilde J_-$
(the boundary lines AB and EF, resp. CD, GH) now lie in the intersection of
those spaces with $\bfc$. So the closure of $L$ in the space
\eqref{double-blowup} is the disjoint union of the closed leaves, one for
every point in $S$. Since each leaf is contractible, this means that the
closure of $L$ on the space \eqref{double-blowup} is diffeomorphic to $S
\times \gamma^2$, for some fixed $\gamma^2$, and is therefore a smooth
manifold with corners of codimension 3.

Now we need to show that the closure of $L$ in the smaller space
\eqref{blowup-XXbh} is a smooth manifold of codimension 3\ah{There is some
redundancy here}. We use the following lemma:

\begin{lemma} Assume that $Z$ is a compact manifold with boundary, that $S
  \subset \partial Z$ is a submanifold and that $V$ is a smooth vector
  field on $Z$ that lifts to $[Z;S]$ to be b-normal at $\tilde S$, the lift
  of $S$ to $[Z;S]$. Suppose that $L \subset Z^\circ$ is a
  submanifold of the interior $Z^\circ$ of $Z$ such that $V$ is tangent to $L$,  the closure $\Lt$ of
  $L$ in $[Z;S]$ is transverse to $\tilde S$ and disjoint from $\partial Z
  \setminus \tilde S$. Finally, assume that at each point $s \in \tilde S
  \cap \Lt$,  
\begin{equation}
\text{the intersection of $T_s \tilde L$ with the tangent space to the
fibre of $\tilde S$ at $s$ is trivial.}  \ilabel{triv1}\end{equation} Then
the closure $\Lc \subset Z$ of $L$ in $Z$ is transverse to $\partial Z$,
and $V |_{\Lc}$ is b-normal to $\partial \Lc$.
\end{lemma}

\begin{proof} We can find coordinates $(x,y,z)$ locally near a point of $S$
  so that $x$ is a boundary defining function for $\partial Z$,  and $S$ is
  given by $\{ x = y = 0\}$. Then coordinates near an interior point $s \in
  \tilde S$ are $x, Y = y/x$ and $z$, and the fibres of $\tilde S$ are
  parametrized by $z$. Near $s$, due to condition \eqref{triv1}, there is a
  splitting of the coordinates $z = (z',z'')$ so that $(x, z')$ furnish
  coordinates on the submanifold $\Lt$. Thus, on $\Lt$, the other
  coordinates are given by smooth functions of $x$ and $z'$: 
$$ Y_i = \tilde Y_i(x,z'), \quad z'' = \tilde z''(x,z').
$$ It follows that on $Z$, the coordinates $(y, z'')$ are given by smooth
functions of $x$ and $z'$, namely $y = x\tilde Y_i(x,z')$, $z'' = \tilde
z''(x,z')$, and hence $L$ is smooth up to, and transverse to, the boundary
of $Z$. Finally, the vector field $V$, restricted to $\Lt$, has the form
$$ x (a \partial_x + \sum_i b'_i \partial_{z'_i} ),
$$ where $a$ and $b'_i$ are smooth functions of $x$ and $z'$. This remains
true when viewed as restricted to $\Lc \subset Z$, which proves the final
statement of the lemma.
\end{proof}

\begin{example} The following simple examples may help to illustrate the
  lemma. First consider the vector field $V = -(x \partial_x + y \partial_y
  + z \partial_z)$ on $Z = \{ (x, y, z) \mid z \geq 0 \}$, let $S \subset Z
  = \{ (0,0,0) \}$ and let $L$ be the flowout from $\{ (x, y, z) \mid x^2 +
  y^2 = 1, z = 1 \}$  via $V$. Then condition \eqref{triv1} is not
  satisfied, and $L$ has a genuine conic singularity at $S$ which is
  resolved by blowup of $S$. 

Second, consider the case where $Z = \{ (x, y, z, x', y', z') \mid z \geq 0
\}$, $V = -(x \partial_x + y \partial_y + z \partial_z)$, $S = \{ (0, 0, 0,
x', y', z') \}$ and $L$ is the flowout from $ \{ (x, y, z, x', y', z') \mid
x^2 + y^2 = 1, z = 1, x' = x, y' = y, z' = z \}$ via $V$. In this case,
\eqref{triv1} \emph{is} satisfied, and the closure of $L$ is a smooth
manifold with boundary with no blowup required. Indeed, we can take
coordinates on $L$ to be $\theta = \tan^{-1}(y'/x')$ and $z$.  
\end{example}

\emph{Completion of the proof of Proposition~\ref{prop:leaf}.}  We apply
the lemma to $L$, with $S$ equal to $W'$ and $V$ equal to $V_r'$. Condition
\eqref{triv1} holds because coordinates on $L$ near $\Wt$ can be taken to be
$\theta, x, y, \mu$ (away from the codimension 3 corner of $L$). The
functions $y^i$ and $\mu_j$ have linearly independent differentials since
this is true on $N^* \Delta_b$ and since $y$ and $\mu$ are invariant under
the flow. Near the codimension 3 corner of $L$, we can take the three
boundary defining functions together with $y$ and $\hat \mu$ and the same
argument goes through. Then the lemma shows that we may blow down $\Wt$ and
$L$ is still a manifold with codimension 3 corners, with $V_r'$ still
b-normal at $\{ \theta = 0 \}$.

At $\tilde J_-$ a totally different argument is needed. Note the asymmetry
between $J$ and $J_-$: the diagonal $N^* \Delta_b$ intersects $J_-$, while
it is disjoint from $J$. To understand the structure of $L$ near $J_-$ we
can start from $N^* \Delta_b \cap J_-$, which is codimension 1 in $L$, and
flow using either $V_l'$ or $V_r'$. In the region $\theta \leq 2$ it
suffices to use $V_r'$. Then since $V_r'/\theta$ is smooth and non
vanishing in this region, we deduce that $L$ is smooth at $J_-$ before
blowup. Therefore, (the closure of) $L$ is a smooth manifold with
codimension 3 corners on the space \eqref{blowup-XXbh}.
\end{proof}

\begin{remark} We emphasize that the blowup at $J$ is essential---it
  resolves genuine conic singularities of the Legendrian $L$---while the
  blowup at $J_-$ resolves no singularities and can be dispensed
  with. Nevertheless,  the blowup at $J_-$ has some good features; in
  particular, it separates all the leaves. On the space
  \eqref{blowup-XXbh}, the leaves join together at $J_-$ like the pages of
  a book joined at the binding.
\end{remark}

%%%%%%%%%%%%%%%%%%%%%%%%%%%%%%%%%%%%%%

%%%%%%%%%%%%%%%%%%%%%%%%%%%%%%%%%%%%%%%%%%%%%%

\section{Parametrix construction}

\subsection{Near the h-scattering diagonal} We begin by using a
semiclassical scattering pseudodifferential operator to remove the diagonal
singularities of the resolvent.  Let $P_C=h^2 \Lap+V+C$ with $C>-\inf
V$. Then in Section~\ref{sec:hpseudos} we showed that $P_C^{-1}$ is a semiclassical pseudodifferential operator of order $(-2,0,0)$. We have
$$ (h^2 \Lap+V-\evosq) P_C^{-1}=\Id-(\evosq+C)P_C^{-1}.
$$
Let $Q$ be an asymptotic sum of the Neumann series 
$$
Q = P_C^{-1} \sum_{j=0}^\infty \big( (\evosq+C)P_C^{-1} \big)^j \in \Psisch^{-2,0,0}(X),
$$
which exists since the differential order of $P_C^{-1}$ is negative. Then 
$$ (h^2 \Lap+V-\evosq) Q = \Id +E_1, \ E_1 \in \Psisch^{-\infty, 0, 0}(X).
$$
Notice that the error term $E_1$ is trivial except at the boundary of the diagonal $\Delta_b \times [0, h_0)$ on $X$, i.e., at  $\Delta_b \times \{ h=0 \}$ and at $\partial \Delta_b \times [0, h_0)$. 
It remains to solve away the error $E_1:$ we now seek a solution $Q'$ to
\begin{equation}
(h^2 \Lap+V-\evosq ) Q' = -E_1; \ilabel{error-1}\end{equation} then
adding $Q'$ to $Q$ will give the desired parametrix.

\subsection{Near the h-b diagonal} 

We begin by considering the kernel of $E_1$ on the double space $\XXbh =
\XXb \times [0, h_0)$.  The fact that $E_1$ is an h-pseudodifferential
operator of differential order $-\infty$ means that its kernel has an
oscillatory integral representation of the form
\begin{equation*}
h^{-n} \int e^{i(z-z') \cdot \zeta/h} e(z, \zeta, h) \, d\zeta \, |dz
dz'|^{1/2}
\end{equation*}
near $\Delta_b$ and away from $\bfc$, and of the form
\begin{equation*}
h^{-n} \int \tilde e(x,\theta,y,y',h,\xi,\eta) e^{i\frac{(\theta-1) \xi +
    (y-y')\cdot \eta}{xh}} \, d\xi \, d\eta \, |dz dz'|^{1/2} \ (\theta = \frac{x}{x'})
\end{equation*}
near $\bfc$. We multiply this by the half-density $|dh|^{1/2}$, to turn it into a density on $X$.   It may then be regarded as
a half-density Legendre distribution of order $(3/4, 1/4)$ associated to
the Legendre submanifold $N^* \Delta_b$, where $\Delta_b \subset X \times
\{ 0 \}$ is the b-diagonal at $h=0$.
Since we wish to solve the
equation \eqref{error-1}, we need to take into account the (left)
characteristic variety $\Sigma_l \subset \Tsfstar(\XXbh)$ of the operator
$h^2 \Lap+V-\evosq$.  The Legendrian $N^* \Delta_b$ is given in the
coordinates of \eqref{canonicaloneform} by\jw{Added $\tau=0$ here.}
$$ \{y=y',\ \theta=1,\ \lambda=-\lambda',\ \mu=-\mu',\tau=0\},
$$ and the left characteristic variety is given in the same coordinates by
$$ \Sigma_l = \{ \lambda^2 + h^{ij}(y)\mu_i \mu_j + V(z) = \evosq \}.
$$ These intersect transversely in a submanifold of dimension $2n-1$, as
proved in Section~\ref{sec:prop}.  Let $L_\pm$ be defined by \eqref{Ldefn};
recall that $L_\pm$ are Legendrian submanifolds with boundary, which
intersect $N^* \Delta_b$ cleanly at $\partial L_\pm$, and are both
transverse to the boundary $\bface$; hence $(N^* \Delta_b\cap \Sigma,
L_\pm)$ have the appropriate geometry for a pair of intersecting Legendre
submanifolds, at least in a neighbourhood of $N^* \Delta_b$.

We now seek to solve away the error term $E_1$ near $\Delta_b$ using an
intersecting Legendrian distribution associated to $(N^* \Delta_b\cap
\Sigma, L_+)$ ; in particular, we would like to find
$$ Q_1 \in I^{1/4;-1/4}(\XXbh; (N^*\Delta_b, L_+))
$$ such that $(h^2 \Lap+V-\evosq )Q_1 -E_1$ is microsupported only at
$L_+$, in a region disjoint from $N^* \Delta_b$. (We choose $L_+$ for the \emph{outgoing} resolvent kernel, and $L_-$ for the \emph{incoming} resolvent kernel; the reason for this is that the coordinate $\nubar_1$ is positive, resp. negative on $L_+$, resp. $L_-$ which implies  having a positive, resp. negative phase function in the oscillatory integral expression for our kernel.)
 To do this we solve away
the singularity at $N^* \Delta_b$ order by order. (This is a standard construction for intersecting Lagrangian or Legendrian distributions; see \cite{MU}.) 

We begin by choosing a $Q_{1,1}$ to solve away the principal symbol of $E_1$ at $N^* \Delta_b$. We do this by choosing the symbol of $Q_{1,1}$ at $N^* \Delta_b$ to be
$$
\sigma^{3/4}(Q_{1,1}) = \sigma^{3/4}(E_1) / \sigma(h^2 \Delta + V - \evosq).
$$
Of course this has a singularity at $N^* \Delta_b \cap L_+$,  but the simple vanishing of $\sigma(h^2 \Delta + V - \evosq)$ at $L_+$ means this is  eligible to be the $N^*\Delta_b$ piece
of the symbol of an intersecting Legendrian distribution with respect to
$(N^*\Delta_b, L_+)$. The compatibility relation \eqref{compat} then
determines the value of the symbol on $L_+$ at $\pa L_+ = L_+ \cap N^*
\Delta_b$; it is essentially given by the residue of the singularity (see Section~\ref{ildsc}). We then specify the symbol at $L_+$ to be that symbol which solves the transport equation \eqref{transport} along $L_+$. Since $V_l$ is transverse to $N^*\Delta_b,$ this is a regular ODE and there is a unique solution with our initial condition specified above. 
This gives a $Q_{1,1} \in I^{1/4;-1/4}(\XXbh; (N^*\Delta_b, L_+))$ such that 
\begin{equation*}
(h^2 \Lap+V-\evosq )Q_{1,1} -E_1 \in I^{5/4, -1/4} (X; N^*(\Delta_b), L_+) \ilabel{int-Leg-1}\end{equation*} near
$N^* \Delta_b$ with principal symbol vanishing at $L_+$. Using \eqref{ex-int-2}, we see the error term is actually in 
\begin{equation}
 I^{7/4, -1/4} (X; N^* \Delta_b)+
I^{9/4, -1/4}(X; N^*(\Delta_b), L_+) 
\ilabel{int-Leg-1a}\end{equation} 
The error will thus be more regular at $N^* \Delta_b$ than
$E_1.$ 

Now we iterate this construction. Assume inductively that we have found
$Q_{1,n}$ such that
\begin{equation}
(h^2 \Lap+V-\evosq )Q_{1,n} -E_1 \in I^{n+3/4, 1/4} (X; N^* \Delta_b)+
I^{n+5/4, -1/4}(X; N^*(\Delta_b), L_+)  \ilabel{int-Leg-11}\end{equation}
in a neighborhood of $N^* \Delta_b.$
We want to improve this by finding $Q_{1,n+1}$ satisfying
\eqref{int-Leg-11} with $n$ replaced by $n+1$. By \eqref{ex-int-1} and
\eqref{ex-int-2} we have to solve away the principal symbol of the first
error term $E_{1,n,1}$ in \eqref{int-Leg-11} at $N^* \Delta_b$, and the
principal symbol of the second error term $E_{1,n,2}$ at $L_+$. We do this
as above: we let $Q'_{1,n}$ have symbol at $N^* \Delta_b$ equal to
$$ \sigma(h^2 \Lap+V-\evosq )^{-1} \sigma(E_{1,n,1})
$$ and symbol at $L_+$ given by solving the transport equation on $L_+$ to
remove the principal symbol of $E_{1,n,2}$ there, using the initial
condition coming from the compatibility condition \eqref{compat}. We cut
off this symbol away from $N^*\Delta_b$ to make it supported in a
neighborhood of $N^*\Delta_b.$ Letting $Q_{1,n+1} = Q_{1,n} + Q'_{1,n}$
completes the inductive step. We can take an asymptotic limit of the
$Q_{1,n}$ obtaining a $Q_1 \in I^{1/4, -1/4}(\XXbh; (N^* \Delta_b, L_+))$
satisfying
\begin{equation}
(h^2 \Lap+V-\evosq ) Q_1 -E_1=E_2 + E_2' \ilabel{E2}\end{equation} with $E_2' \in I^{\infty, 1/4}(N^* \Delta_b)+ I^{\infty, -1/4}(N^* \Delta_b, L_+)$ and  $E_2 \in
I^{1/4, -1/4}(X; L_+)$, arising from the cutoff,  \emph{microsupported away from $N^*\Delta_b$}. % in aneighborhood of $N^* \Delta_b.$}
In fact, we can improve this statement to $E_2' \in I^{\infty, 1/4}(N^* \Delta_b)+ I^{\infty, 3/4}(N^* \Delta_b, L_+)$ since $h^2 \Delta + V - \evosq$ is characteristic at $h^{-1}\partial_{\bfc}L_+$ \emph{for every $h > 0$} which automatically gives us an extra order of vanishing at $L_+$, hence an improvement by $1$ in the order at $\bfc$. 
  
\subsection{At the propagating Legendrian}\ilabel{atpl}
As in the finite energy case, we now consider the Legendrian conic
pair $$ \tilde L(\evosr)= (L(\evosr), L_2^\sharp(\evosr)),
$$ from Proposition~\ref{prop:leaf}.  We set aside the error $E_2'$ until Section~\ref{h>0} and seek here to solve away the error $E_2$
from \eqref{E2} by adding a Legendre distribution $Q_2 \in I^{1/4,p;
\mathbf{r}}(\XXbh, \tilde L(\evosr))$, where $p$ is the order at $L_2^\sharp$
and $\mathbf{r}$ represents orders $(r_{\bfc}, r_{\rb}, r_{\lb})$ at the
other boundary hypersurfaces. We shall see that the orders are $p=n/2-3/4$,
$r_{\bfc} = -1/4$, $r_{\rb} = r_{\lb} = n/2 - 1/4$.  Our precise goal in
this step in the construction is to find $Q_2$ so that
\begin{equation}
(h^2 \Delta + V - 1) Q_2 - E_2 \in I^{+\infty, n/2+1/4; (3/4,
n/2-1/4,n/2+7/4)}(\XXbh, \tilde L(\evosr)); \ilabel{Q2w}\end{equation} that
is, the error has been completely solved away at $h=0$.  The space in which
the error lies is the same as $h^\infty I^{-1/2, (n-2)/2;(n-1)/2, (n+3)/2}
(\XXb, \partial_{\bfc} L, L_2^\sharp)$ (see Section~\ref{resid-conic}), that is, a family of Lagrangian
distributions associated to the boundary of $L$ at $\bfc$ and to
$L_2^\sharp$, and rapidly decreasing as $h \to 0$. This will reduce the
problem to a parametrized version of the problem already studied in
\cite{hasvas2}.

Again we solve away error terms, this time on $L_+$, order by order. The
first step is to find $Q_{2,1}$ solving \eqref{Q2w} with the order $\infty$
at $h=0$ replaced by $9/4$. The order of $Q_{2,1}$ must be $5/4$ at $L_+$
and $-1/4$ at $\bfc$.  By \eqref{transport}, to solve \eqref{Q2w} it
suffices to obtain $q_2$ satisfying the ODE
\begin{equation}
\Big(\mathcal{L}_{V_l'} - \big( \frac1{2} + m - \frac{2n+1}{4} \big)
\frac{\partial p_L}{\partial \lambda} + f_l \Big) \sigma^{1/4}(q_2) \otimes
\abs{d(h\sigma x')} = e_2, \quad m = \frac{1}{4}
\ilabel{l-tr-eq}\end{equation} and\footnote{The factor $d(h\sigma x')$ in
the equation above is a `formal factor' adjusting for the difference in the
symbol bundle \eqref{symbolbundle} when the order $m$  changes by $1$} with the
`initial condition' that the symbol $q_2$ vanishes near $N^* \Delta_b$,
reflecting the fact that we do not want to disturb our parametrix near $N^*
\Delta_b$.  Here we are using coordinates induced from the canonical 1-form
\begin{equation}\ilabel{canonicaloneform-l}
\lambda d\big( \frac1{x' \sigma h} \big) + \lambda' d\big( \frac1{x' h}
\big) + \tau d\big( \frac1{h} \big) + \mu \frac{dy}{x' \sigma h} + \mu'
\frac{dy'}{x' h}
\end{equation}
which are valid for $\sigma = \theta^{-1} \leq 2$, say, thus valid near the
corner $\lb \cap \bfc$. Also $f_l$ denotes the subprincipal symbol of the
(left) operator $h^2\Lap+V-\evosq$. Since $V_l$ is smooth and nonvanishing
in the interior of $L_+$, this has a unique smooth solution in the interior
of $L_+$. We proceed to analyze the regularity of the symbol at the
boundary of $L_+$. This will be done exploiting the b-normal vector fields
from Lemma~\ref{b}.  Consider $L_+ \cap \lb$. Here the ODE takes the form
\begin{equation}
\Big( \lambda \mathcal{L}_{(\rho_{\lb} \partial_{\rho_{\lb}} + \rho_{\lb}
\Nub)} - (m - \frac{2n-1}{4} ) \lambda + f_l \Big) q_2 = 0
\ilabel{l-tr-eq-2}\end{equation} where $\Nub$ denotes a vector field on
$L_+$ tangent to the boundary at $\rho_{\lb} = 0$. We recall that the
sub-principal symbol $f_l$ vanishes at $\mu = 0$, hence is $O(\rho_{\lb})$
at $\rho_{\lb} = 0$. So we may write $f_l = \rho_{\lb} \tilde f_l$.  Also
recall that $q_2$ is a half-density and it is convenient to write it as a
b-half-density, that is, $q_2 = \tilde q_2 |d\rho_{\lb} d\rho_{\bfc}
d\ss/\rho_{\lb}|^{1/2}$; note that this half-density is invariant under Lie
derivation by $\rho_{\lb} \partial_{\rho_{\lb}}$. We get an equation for
$\tilde q_2$ of the form
$$ \Big( \lambda \rho_{\lb} \partial_{\rho_{\lb}} + \rho_{\lb} \Nub - (m -
\frac{2n-1}{4} ) \lambda + \rho_{\lb} \tilde f_l \Big) \tilde q_2 = 0
$$
hence we obtain
$$\tilde q_2 \in \rho_{\lb}^{(2n-1)/4 - m} \CI(L)
$$ at least locally. Thus, using Proposition~\ref{prop:ex-conic-2}, the order at
$\lb$ is $(2n-1)/4$.

To show regularity at $\rb$, we use the fact that near $\rb$ the symbol
$q_2$ automatically satisfies the right transport equation as well; that
is, if we define $q_2$ using the right transport equation rather than the
left, then we get the same result. We shall not give the proof here since
it is essentially identical to the proof of the analogous statement in
\cite{hasvas2}, section 4.4. Then, reversing the left and right variables
in the argument above proves regularity at $\rb$ with the order also equal
to $(2n-1)/4$.

To show regularity at $L^\sharp(\evosr)$, we combine both vector fields.  By
Lemma~\ref{b}, the vector field $V_l' + V_r'$ is b-normal to $\Jt,$ which is
the blowup of $L_2^\sharp;$ thus we add together the left and the right
transport equations. The right transport equation, written with respect to
the variables in \eqref{canonicaloneform}, takes the form of
\eqref{l-tr-eq} with the left and right variables switched:
\begin{equation}
\Big(\mathcal{L}_{V_r'} - \big( \frac1{2} + m - \frac{2n+1}{4} \big)
\frac{\partial p_R}{\partial \lambda} + f_r \Big) \sigma^{1/4}(q_2) \otimes
\abs{d(h\theta x)} = e_2, \quad m = \frac{1}{4}
\ilabel{r-tr-eq}\end{equation} To compare the two symbols, we must express
them with respect to the same total boundary defining function. The total
boundary defining function used in \eqref{l-tr-eq} is $hx'$, while that
used in \eqref{r-tr-eq} is $hx$. The ratio is $\theta$; in view of the
presence of the factor $|d\mathbf{x}|^{m-N/4}$ in the symbol bundle (see
\eqref{symbolbundle}), the symbol gets multiplied by a factor of
$\theta^{m-N/4}$ when we switch (where $N = 2n+1$ here). Thus, with
respect to the total boundary defining function $hx$,
\begin{equation}
\Big(\mathcal{L}_{V_r'} - \big( \frac1{2} + m - \frac{2n+1}{4} \big)
\frac{\partial p_R}{\partial \lambda} + f_r \Big) \sigma^{1/4}(\theta^{m -
N/4}q_2) \otimes \abs{d(h\theta x)} = e_2, \quad m = \frac{1}{4}
\end{equation} We can multiply this equation by
$\theta^{N/4-m}$ and add it to \eqref{l-tr-eq}. The effect of this is that
the $-\lambda' \theta \partial_\theta$ term in $V_r'$ gives a contribution
of $-(m-N/4)\lambda'$. As a result (taking into account $\lambda = \lambda' +
O(\rho_{\Jt})$ at $\Jt$ and $\pa p_L/\pa \lambda=2\lambda,$ $\pa p_R/\pa \lambda'=2\lambda',$),
\begin{equation}
\Big(\mathcal{L}_{V_l' + V_r'} - \big(\frac 12 + m - \frac{2n+1}{4} \big)
\lambda +\frac 12 \lambda' + f_l + f_r \Big) \sigma^{-1/4}(q_2) \otimes
\abs{d(h\theta x)} = 0, \quad m = \frac{1}{4}.
\ilabel{c-tr-eq}\end{equation} Since $f_l$ vanishes at $\mu = 0$ and $f_r$
vanishes at $\mu' = 0$, they both vanish at $\Jt$. So we can write $f_l +
f_r = \rho_{\Jt} \tilde f$.  Also, of course
$\lambda=\lambda'+O(\rho_{\Jt}).$ Thus \eqref{c-tr-eq} amounts to an equation
of the form (again writing $q_2 = \tilde q_2$ times a b-half-density)
$$ \Big( \lambda \rho_{\Jt} \partial_{\rho_{\Jt}} + \rho_{\Jt} \Nub - (m -
\frac{2n-3}{4} ) \lambda + \rho_{\Jt} \tilde f \Big) \tilde q_2 = 0
\implies \tilde q_2 \in \rho_{\Jt}^{(2n-3)/4 - m} \CI(L)
$$ locally. This shows regularity of the symbol at $L^\sharp$, and that the
order $p$ at $L^\sharp$ is $n/2 - 3/4$. The error term when we apply the
operator is given by \eqref{Q2w} with $9/4$ replacing $\infty$. This is
because the operator is characteristic at $L_+$, and at the induced
Legendrians at $\bfc$ and at $\lb$ (but not at $\rb$); in addition we have
solved the transport equations at $L_+$ and, trivially, at the left
boundary (this because the transport operator is trivial at $\lb$ at order
$(2n-1)/4$) so we gain two orders in each of these two cases.

Now we iterate the procedure. Assume inductively that we have found
$Q_{2,k}$ satisfying
\begin{equation} 
(h^2 \Delta + V - \evosq) Q_{2,k} - E_2 \in I^{5/4+k, n/2+1/4; (3/4,
n/2-1/4,n/2+7/4)}(\XXbh, \tilde L(\evosr)).  \ilabel{Q2k}\end{equation} We
want to improve the error term to have order $5/4 + k + 1$ at $L_+$. To do
this, we solve the transport equation at order $5/4 + k$ along $L_+$, and
as above the main point is to determine the regularity of the solution at
the boundary of $\hat L_+$.  Consider the solution of \eqref{l-tr-eq}, with
$m$ replaced by $1/4 + k$, and with the right hand side replaced by the
error term in \eqref{Q2k}. Using Proposition~\ref{prop:ex-conic-2} the right
hand side is $O(\rho_{\lb}^{(2n-1)/4 - (1/4 + k) + 1})$. Therefore the
right hand side avoids the indicial root, in this case $(2n-1)/4 - (1/4 +
k)$ which would lead to possible log terms in the solution, and we see that
the solution is in $\rho_{\lb}^{(2n-1)/4 - (1/4 + k)} \CI(L_+)$
locally. Since, as noted above, we get the same parametrix if we solve via
the right transport equation instead of the left, the same result is true
at $\rb$.  Similar reasoning also shows that the symbol is in
$\rho_{\Jt}^{(2n-3)/4 - (1/4 + k)} \CI(\hat L_+)$ at $\rho_{\Jt} = 0$; it
is essentially the same argument as in \cite{hasvas2}, section 4.4, so we omit
it\ah{Adequate?}. This completes the inductive step.  Taking an asymptotic
limit of the $Q_{2,k}$ gives a correction term satisfying \eqref{Q2w}.

\begin{remark} If the potential $V$ is replaced by $h^2 V$, then $V$ does not appear in the principal symbol of $H - \evosq$ and therefore does not affect the
bicharacteristic flow or the Legendrian $L$; on the other hand, it contributes an additional  error term on the right hand side of \eqref{l-tr-eq}. Because of our assumption $V = O(x^2)$, this additional  error term is also $O(\rho_{\lb}^{(2n-1)/4 - (1/4 + k) + 1})$, and therefore the construction goes through as above.
\end{remark}

\subsection{At the boundary for $h > 0$}\ilabel{h>0}
Our error term is now of the form (using Section~\ref{int-resid} and \ref{resid-conic})
\begin{multline*}  E_2' + E_3 \in  
I^{\infty, 1/4}(X, N^* \Delta_b; \sfOh)+ I^{\infty, 3/4}\big( X, (N^* \Delta_b, L_+), \sfOh \big) \\ +
I^{+\infty, n/2+1/4; r_{\bfc} + 1/4, r_{\lb} + 1/4, r_{\rb} + 1/4}\big(X, \tilde L(\evosr), \sfOh)
\end{multline*}
where $r_{\bfc} = 1/2$, $r_{\lb} = (n+3)/2$ and $r_{\rb} = (n-1)/2$.  Equivalently, the error term is a smooth, $O(h^\infty)$ function of $h$ valued in 
\begin{multline*}
I^0(\MMb, N^* \Delta_b; \sfOh) + I^{1/2}(\MMb, (N^* \Delta_b, h^{-1}\partial_{\bfc}  L_+), \sfOh)  \\ + I^{r_{\bfc}, n/2; r_{\lb}, r_{\rb}}(\MMb, h^{-1}\partial_{\bfc}  \tilde L(\evosr), \sfOh).
\end{multline*}
We now use the results of \cite{hasvas2} to solve away these errors. The
main point here is to keep track of powers of $h$: our error terms are
rapidly decreasing in $h$ and we would like to find a correction term that is also rapidly
decreasing in $h$. Examining the construction in \cite{hasvas2}, we see that
the vector fields in the transport equations are linear in $\lambda =
h^{-1}$, while $\lambda$ appears polynomially in the right hand side due to
derivatives bringing down powers of $\lambda$ from the phase and from the factor $\lambda^2$ in front of the potential. It follows
that the correction term is $O(h^\infty)$ if the error terms are
$O(h^\infty)$. Thus, by \cite{hasvas2} we can solve away the error term $E_3$ above 
with a term $Q_3$ in the space  
\begin{multline*}
Q_3 \in h^{\infty} \CI \big( [0, h_0); I^{-1/2} (\XXb, \partial_{\bfc} N^*
\Delta_b, h^{-1}\partial_{\bfc} L_+) \big) \\ + h^\infty \CI \big( [0, h_0);
I^{-1/2, (n-2)/2;  (n-1)/2} (h^{-1}\partial_{\bfc} L_+, h^{-1}L_2^\sharp)
\big),\end{multline*} or equivalently
\begin{multline*}
Q_3 \in  I^{-\infty, -1/4} \big(\XXb, (\partial_{\bfc} N^*
\Delta_b, \partial_{\bfc} L_+), \sfOh \big) \\ 
+ I^{\infty, (2n-3)/4;  (2n-1)/4, -1/4} \big( X, (\partial_{\bfc} L_+, L_2^\sharp), \sfOh \big).
\end{multline*}
up to a new error term $E_4$ where the expansions at $\lb$, $\bfc$ are trivial, but the expansion at $\rb$ has not been improved. (We recall that when we act with the operator on the left variable, we can improve our parametrix at $\lb$ order by order using the symbol calculus, but to improve at $\rb$ we have to solve global problems of the form $(h^2 \Delta + V - \evosq) v=f$, which of course we cannot do until we have constructed the resolvent kernel! So it cannot be expected that we get any improvement at $\rb$.)  
Thus
$E_4 \in I^{\infty, \infty; \infty, \infty, r_{\rb}}(X, (L_+, L_2^\sharp), \sfOh)$, $r_{\rb} = (n-1)/2$,or more simply,
\begin{equation}
E_4 \in h^\infty x^\infty (x')^{(n-1)/2} e^{i\evo/x'h} \CI(\XXbh; \Omega_{\sf}^{1/2}).
\ilabel{E4}\end{equation} %(where $\mu$ is a scattering-fibredhalf-density).  
In summary, we have found a parametrix $G(h)$ in the space
\begin{multline}
\Psisch^{-2, 0, 0}(\XXbh) \otimes |dh|^{1/2} + I^{1/4;-1/4}(\XXbh;
(N^*\Delta_b, L_+)) \\ + I^{1/4, (2n-3)/4;  (2n-1)/4, -1/4}(L_+,
L_2^\sharp) \ilabel{parspace}\end{multline} such that
\begin{equation}
(h^2 \Delta + V - \evosq) G(h) - \Id = E_4 \in h^\infty x^\infty
(x')^{(n-1)/2} e^{i\evo/x'h} \CI(\XXbh; \Omega_{\sf}^{1/2}).
\end{equation}

%%%%%%%%%%%%%%%%%%%%%%%%%%%%%%%%%%%%%%%

\section{Resolvent from parametrix}

Using the parametrix $G(h)$ constructed in the previous section, which lies
in the space \eqref{parspace}, we can show that the resolvent kernel itself
lies in this space for small $h$.  The error term $E_4$ in the previous
section is compact on weighted $L^2$ spaces $x^s L^2(X)$, for $s >
1/2$. Moreover, the Hilbert-Schmidt norm of $E_4$, thought of as an
operator on $x^s L^2(X)$ parametrized by $h$, tends to zero. It follows
that $\Id + E_4$ is invertible for small $h$. Let the inverse be $\Id +
F(h)$. Then the identity
\begin{equation} -F = E_4 + E_4^2 + E_4FE_4
\ilabel{EFE}\end{equation}
shows that $F$ also lies in the space \eqref{E4}. Finally, the resolvent
kernel is
$$ R(h) = G(h) + G(h) F(h).
$$ Since $F(h)$ is rapidly decreasing as both $h \to 0$ and as $x \to 0$,
it follows from this that $R(h)$ is also in the space \eqref{parspace};
indeed the rapid decrease of $F(h)$ in $x$ wipes out all expansions of $G(h)$ at
$\bfc$ and at $\rb$ in this composition, and the rapid decrease of $F(h)$ as $h \to 0$ wipes out all expansions of $G(h)$
as $h \to 0$. We are left with the expansion of $G(h)$ at $\lb$. This takes the form
$e^{i/xh} x^{(n-1)/2}$ times smooth functions of the other variables
(ignoring density factors), and the result of the composition is an
operator of the form
$$
x^{(n-1)/2} (x')^{(n-1)/2} e^{i\evo/xh} e^{i\evo/x'h} h^\infty \CI(M^2 \times [0, h_0)),
$$
rapidly decreasing at $\bfc$, at $\rb$, and as $h
\to 0$. So $G(h) F(h)$ is a particularly simple example of an operator in
\eqref{parspace} (corresponding to the term $u_6$ in Section~\ref{6.5.2}). This
completes the proof of Theorem~\ref{main1}.

\part{Applications}

%%%%%%%%%%%%%%%%%%%%%%%%%%%%%%%%

\section{Spectral measure and Schr\"odinger propagator}

In this section we prove Corollary~\ref{sp-cor} and Theorem~\ref{prop-thm}. Let $H$ denote $\Delta +  V$ in this section, let $R_\pm = (h^2 \Delta + h^2 V - (1 \pm i0))^{-1}$, and let $\lambda = h^{-1}$. By the remark at the end of Section~\ref{atpl}, $R_\pm$ has the same structure as the semiclassical resolvent with no potential term. (The term $h^2 V$ vanishes to second order at $\partial X$ so it is not present in the principal symbol of the operator, and hence plays no role in determining the Legendrians $L$ or $L_2^\sharp$. It does, of course, affect the \emph{symbol} of the resolvent, but does not change its regularity properties.)

A direct consequence of Theorem~\ref{main1} is the structure of the spectral measure $dE(\lambda^2)$. By Stone's theorem, we have
\begin{gather*}
dE(\lambda^2) = \frac1{2\pi i} \big( (H - (\lambda + i0)^2)^{-1} - (H - (\lambda - i0)^2)^{-1} \big) 2\lambda d\lambda \\
= \frac{1}{\pi i} \big( R_+(h) - R_-(h) \big)\otimes \big| \frac{dh}{h^2} \big|^{1/2} . 
\end{gather*}
We then have immediately from Theorem~\ref{main1} that $d(E(\lambda^2)) \otimes |dh/h^2|^{-1/2}$ is in the sum of spaces 
\begin{gather*}
\Psi^{-2, 0, 0}(X) +  I^{-1/4; -1/4}((N^* \Delta_b, L_+), X; \sfOh) + I^{1/4, n/2 - 3/4; n/2 - 1/4, -1/4}((L_+, L_2^\sharp), X; \sfOh) \\
+ I^{-1/4; -1/4}((N^* \Delta_b, L_-), X; \sfOh) + I^{1/4, n/2 - 3/4; n/2 - 1/4, -1/4}((L_-, L_2^\sharp), X; \sfOh).
\end{gather*}
However, the kernel of $dE(\lambda^2)$ solves an elliptic equation
$$
( \Delta + V - \lambda^2) dE(\lambda^2) = 0.
$$
So there can be no singularity of $dE(\lambda^2)$ at $N^* \Delta_b$, except at the characteristic variety $N^* \Delta_b \cap \Sigma_l = N^* \Delta_b \cap L$. 
along the diagonal. Moreover, $dE(\lambda^2)$ must be Legendrian along $L = L_+ \cup L_-$ at the intersection $L_+ \cap L_- = L \cap N^* \Delta_b$, since it is Legendrian away from $N^* \Delta_b$ and Legendrian regularity propagates along the bicharacteristic flow, which is non-vanishing at $L \cap N^* \Delta_b$. Thus in fact
\begin{gather*}
dE(\lambda^2) \otimes \big| \frac{dh}{h^2} \big|^{-1/2} \in  I^{1/4, n/2 - 3/4; n/2 - 1/4, -1/4}((L, L_2^\sharp), X; \sfOh) ,
\end{gather*}
which is 
Corollary~\ref{sp-cor}. 

\

We now turn to the proof of Theorem~\ref{prop-thm}. We begin with some
preliminaries on the geometry of the b-double space $\MMb$. A total
boundary defining function for this space can be taken to be $\xxb = (r^2 +
(r')^2)^{-1/2}$.  We need to consider small neighbourhoods of the
b-diagonal $\Delta_b$ in $\MMb$.  A neighbourhood is given, for example, by
$$
\{ (z, z') \mid  d(z, z') < \epsilon /\xxb \} = \{ (z, z') \mid  d(z, z') < \epsilon  \sqrt{r^2 + (r')^2} \} 
$$ for $\epsilon > 0$. Let $\phi$ be a smooth function on $[0, \infty)$
equal to $1$ on $[0,1]$ and equal to $0$ on $[2, \infty)$. Then
$\phi(d(z,z') \xxb/\epsilon)$ is a smooth function on $\MMb$ equal to $1$
at $\Delta_b$ and supported near $\Delta_b$ (for small $\epsilon$). Abusing
notation somewhat, we shall denote this function on $\MMb$ simply by
$\phi$. The local injectivity radius on $M$ is bounded below by $c r$ for
some $c > 0$; we shall assume that $\epsilon > 0$ is chosen so that the
local injectivity radius is at least $10\epsilon r$.  Then the square of
the distance function $d(z,z')^2$ will be smooth on the support of $\phi$.

To obtain the kernel of the propagator $e^{-itH/2}$,  $H = \Delta + V$, consider the integral over the spectrum:
\begin{equation}
e^{-itH/2} = \int_0^\infty e^{-it\lambda^2/2} dE(\lambda^2) .
\ilabel{spec-int}\end{equation} 
To deal with this integral we
break it into several pieces. We first use a spectral cutoff. Let us insert
$1 = \chi_1(\lambda \sqrt{t}) + \chi_2(\lambda \sqrt{t})$ into the
integral, where $\chi_1$ is equal to $1$ on $[0,1]$ and equal to $0$ on
$[2, \infty)$. The $\chi_1$ term yields the operator $\chi_1(H \sqrt{t})
e^{-itH/2}$. Letting $s = \sqrt{t}$, this is a $C_c^\infty$ function of $s^2 H$
and is therefore a semiclassical pseudodifferential operator (in $s$) of
order $-\infty$ \cite{DSj}. 
%Letting $\MMsc$ denote the `scattering
%doublespace' $[\MMb; \partial \Delta_b]$ and $\Delta_{SC}$ the
%corresponding lift of the diagonal (see \cite{Melrose:Spectral}) we
%find that the kernel is therefore smooth on the further blown up space
%\begin{equation}
%[\MMsc \times [0, s_0)_s; \Delta_{\SC} \times \{ 0 \} ]
%\ilabel{par-space}\end{equation}
%and essentially supported on the blowup face at $s=0$, i.e.\ 
In particular, it is smooth away from the diagonal, and rapidly
decreasing %at the lift of $\MMsc \times \{ 0 \}$ (equivalently, it is rapidly decreasing 
as $d(z,z')/s \to \infty$. Let us write this kernel
$U_{\near, 1}$. Notice that $(1 - \phi) U_{\near, 1}$ is residual, i.e.\ in
$\CIdot(\MMb \times [0, t_0))$, for any function $\phi$ localized near
  $\Delta_b$ as above (i.e., for any $\epsilon > 0$).

Now consider the integral with $\chi_2(\lambda \sqrt{t})$ inserted. We now
localize based on the value of the phase function $\psi/\xx$, $\xx = \xxb h
= \xxb/\lambda$,  in the representation of the semiclassical resolvent as a
Legendre distribution. Let us write
$$
1 = \chi_n(\psi/\epsilon) + \chi_i(\psi/\epsilon) +\chi_f(\psi/\epsilon)
$$
where $\chi_n$ is supported in $[0, 1/2]$, $\chi_i$ is supported in $[1/4,
  3]$ and $\chi_f$ is supported in $[5/2, \infty)$. We obtain three
  kernels, denoted $U_{\near, 2}$, $U_\innt$ and $U_\far$, by inserting the
  cutoffs $\chi_2(\lambda \sqrt{t})$ $\chi_\bullet(\psi/\epsilon)$ into
  \eqref{spec-int}. Let us also define $U_\near = U_{\near, 1} + U_{\near,
    2}.$  Thus we may write the exact propagator
$$
e^{-itH/2} = U_{\near}+U_{\innt}+U_{\far}.
$$

\begin{lemma}\ilabel{Ufar} (i) The kernel $d\phi \cdot U_\near$ is in $\CIdot(\MMb \times [0, t_0))$. 

(ii) The kernel  
\begin{equation} 
(1 - \phi) (D_t + \half H) U_\innt + (D_t + \half H) U_\far 
 \ilabel{triv-kernels}\end{equation}
is in $\CIdot(\MMb \times [0, t_0))$.

(iii) $U_\innt$ is a quadratic Legendre distribution associated to $Q(L)$, and $U_\far$ is a quadratic Legendre distribution associated to $(Q(L), Q(L_2^\sharp))$. 
\end{lemma}

\begin{proof} 
(i)  We have already observed that this is true in the case of $U_{\near,
    1}$ so consider $U_{\near, 2}$. 
Observe that
$$
e^{i\lambda\psi/\xxb} = \frac{-i \xxb}{\lambda} \frac1{d_v \psi} d_v e^{i\lambda\psi/\xxb},
$$ and that on the support of $1- \phi$ and on the support of
$d\chi_n(\psi/\epsilon)$ we have $d_v \psi \neq 0$. (This is because $d_v
\psi = 0$ implies that $\psi/\xxb = d(z,z')$, yet $\xxb d(z,z') \geq
\epsilon$ on the support of $1 - \phi$ and $\psi \leq \epsilon/2$ on the
support of $d\chi_n$.)  Thus we can integrate by parts in $v$ as many times
as we like.\footnote{If there are no $v$ variables then we simply use the
  fact that $\psi/\xxb = d(z,z').$}
Each integration by parts gains us
$\xxb/\lambda$. This allows us to absorb any number of spatial or
$t$-derivatives, as well as any number of negative powers of $\xxb$ or $t$
(remembering that the combination $\lambda^{-2} t^{-1}$ is bounded on the
support of $\chi_2$).  This proves membership in $\CIdot(\MMb \times [0,
  t_0))$.  

(ii) Let us start with the first term, $(1 - \phi) (D_t + \half H)
  U_\innt$. $U_{\innt}$ is given by a finite sum of integrals of the form
$$ \int \int e^{-it\lambda^2/2} e^{i\lambda \psi(\cdot, v)/\xxb}
\chi_2(\lambda \sqrt{t}) \chi_i(\psi/\epsilon) a(\lambda, \cdot, v) \, dv
\, d\lambda.
$$ Here, $\cdot$ refers to the spatial variables on $\MMb$. If we apply
$(D_t + \half H)$ to the integral then the result vanishes if none of the
derivatives hits one of the cutoffs $\chi_i(\psi/\epsilon)$ or
$\chi_2(\lambda \sqrt{t})$, so $(1 - \phi) (D_t + \half H) U_\innt$ is a
sum of terms of the form
$$
(1 - \phi) \int \int  e^{-it\lambda^2/2} e^{i\lambda \psi(\cdot, v)/\xxb} \chi_2'(\lambda \sqrt{t}) \chi_i(\psi/\epsilon) \tilde a(\lambda, \cdot, v) \, dv \, d\lambda
$$
or
$$
(1 - \phi) \int \int  e^{-it\lambda^2/2} e^{i\lambda \psi(\cdot, v)/\xxb} \chi_2(\lambda \sqrt{t}) \chi_i^{(k)}(\psi/\epsilon) \tilde a(\lambda, \cdot, v) \, dv \, d\lambda.
$$
where $k$, the number of derivatives falling on $\chi_i$, is either $1$ or $2$. In the first case, we can integrate by parts in $\lambda$ as many times as we like, using the identity
$$
e^{i(-t\lambda^2/2 + \lambda \psi/\xxb)} = \frac{-i\xxb}{-t\lambda \xxb + \psi} \frac{\partial}{\partial \lambda} e^{i(-t\lambda^2/2 + \lambda \psi/\xxb)}
$$
and the fact that the denominator is bounded below since $\psi \geq
\epsilon/4$ on the support of $\chi_i(\psi/\epsilon),$ $\xxb$ is a bounded
function, and it suffices to consider only times $t\ll 1.$  This allows us to
reduce the order of the symbol in $\lambda$, and increase the order in
$x_1$ and $x_2$,  as much as we like. Using the same reasoning as in part
(i), the kernel is in $\CIdot(\MMb \times [0, t_0))$. Exactly the same
  arguments allows us to dispose of the terms coming from $(D_t + \half H)
  U_\far$ when a derivative hits $\chi_2$.  

In the case of the second integral, we need to further divide into two cases, according as the derivative $\chi_i^{(k)}$ is supported in $[1/4, 1/2]$ or in $[5/2, 3]$. In the first case, supported in $[1/4, 1/2]$, we can integrate by parts in $v$ as many times as we like, as in part (i), and we see that these terms are in $\CIdot(\MMb \times [0, t_0))$. In the second case, supported in $[5/2, 3]$,  we note that modulo $\CIdot(\MMb \times [0, t_0))$ we can replace the factor $1 - \phi$ by $1$, for exactly the same reason. 

Now we see that these terms, with $\chi_i^{(k)}(\psi/\epsilon)$ supported
in $[5/2, 3]$ and with $1 - \phi$ replaced by $1$, exactly cancel the
remaining terms from $(D_t + \half H) U_\far$, since
$\chi_i^{(k)}(\psi/\epsilon) = - \chi_f^{(k)}(\psi/\epsilon)$ when
restricted to the interval $(\psi/\epsilon) \in [5/2, 3]$. We conclude that
\eqref{triv-kernels} is in $\CIdot(\MMb \times [0, t_0))$.

(iii) This follows immediately from Propositions~\ref{sph} and \ref{sph-conic}. 
\end{proof}

It appears to be difficult to determine the microlocal nature of $U_{\near}$ using the integral \eqref{spec-int}. One reason is that the spectral cutoffs $\chi_1, \chi_2$, needed in order to apply Propositions~\ref{sph} and \ref{sph-conic} in part (iii) of the above lemma,  interfere with the microlocal nature of the pieces. In particular, the piece $U_{\near, 1}$ does \emph{not} lie in the space \eqref{propspace}. We shall see that this is an artifact of the spectral cutoffs  and the sum $U_{\near, 1} + U_{\near, 2}$ \emph{does} lie in \eqref{propspace}. To see this we need to change strategy. What we shall do is construct a parametrix in the near-diagonal region, and show that we can glue it together with the kernel constructed above to obtain the true propagator modulo a $\CIdot(\MMb \times [0, t_0))$ error. 

In the near-diagonal region we use the same ansatz as in Step 1 of \cite{HW}. For the reader's convenience we recall  that this takes the form
$$
(2\pi it)^{-n/2} e^{i\Phi(z,z')/t} \sum_{j=0}^\infty t^j a_j(z, z').
$$
We want this to be a formal solution, so we apply the operator $t^2 D_t + t^2/2 \Delta$ and solve the resulting equations to each order in $t$. The first is the eikonal equation $-\Phi +
g(\nabla_z \Phi, \nabla_z \Phi) = 0$ which has the exact solution $\Phi(z, z') = d(z,z')^2/2$. Thus we see that this is a Legendrian associated to the same Legendre submanifold, namely $Q(L)$, to which $U_\innt$ and $U_\far$ are associated.
The remaining equations are transport equations taking the form (in normal coordinates $z$ about $z'$)
\begin{equation*}\begin{gathered}
(z_i + O(\abs{z}^2)) \frac{\pa}{\pa z_i} a_0 = f\cdot a_0,  \\
 (z_i + O(\abs{z}^2)) \frac{\pa}{\pa z_i} a_j + j a_j = f \cdot a_j - \frac{i}{2} \Delta_z a_{j-1}  \ (j \geq 1)
\end{gathered}
\quad f = \half
\Delta \Phi + \frac n2
 = O(z), 
 \end{equation*}
where all terms are smooth. These equations have unique solutions with $a_j$ smooth and $a_0(0) = 1$. We cut this formal solution off by multiplying by a cutoff function $\phi(d(z,z')/\epsilon r')$ localizing near $\Delta_b$. 

This argument only applies away from the front face of $\MMb$ since the
analysis of \cite{HW} was only carried out there. However, the near-diagonal
ansatz above holds uniformly up to $\bfc \subset \MMb$, i.e.\ in a full
neighbourhood of  $\Delta_b \subset \MMb$. We proceed to show this. We
first note that the function $\Psi = d(z,z')^2/2(r')^2$ is a smooth
function on $\MMb$ in a neighbourhood of $\Delta_b$. In fact, if we take
coordinates $x', \sigma = x/x', y'$ and $y$ locally near $\Delta_b$, where
$y'$ is a local  coordinate on $\partial M$ and for a fixed $y'$, $y$ are
normal coordinates on $\partial_M$ centred at $y'$ (hence, $y$ is
\emph{not} a coordinate lifted from the left factor of $\partial M$), then
$\Delta_b$ is defined by $\{ \sigma = 1, y = 0 \}$ and near $\Delta_b$,
$$
\Psi = (\sigma - 1)^2 + \sum y_i^2 + \text{ terms vanishing to third order at } \Delta_b.
$$
On the other hand, the operator $t^2  (D_t + \half \Delta)$ takes the form
$$
t^2   D_t + (t x' \sigma^2 D_\sigma)^2 + (n-1) t^2 x' \sigma^3 \partial_\sigma + h^{ij}(x) \Big(  (t x' D_{y_i}) (t x' D_{y_j}) + \Gamma_{ij}^k(x) (t^2 (x')^2 D_{y_k}) \Big)
$$
where $\Gamma(x)$ is the Christoffel symbol for the metric $h(x)$. Let us seek a formal solution, as a series in $t$, near the boundary of $\Delta_b$. It takes the form
$$
(2\pi it)^{-n/2} e^{i\Psi/t(x')^2} \sum_{j=0}^\infty t^j b_j(x', \sigma, y, y'), \quad \text{ with } b_j \text{ smooth.}
$$
Since $h^{ij} = \delta_{ij}$ at $y = y'$ and $\Gamma = O(y)$ there, it follows then that we end up with transport equations for the $b_j$ of the form 
\begin{equation*}\begin{gathered}
\Big( y_i  \frac{\pa}{\pa y_i} 
+ (\sigma - 1)  \frac{\pa}{\pa \sigma} + W \Big)
b_0 = f\cdot b_0,  \\
 \Big( y_i  \frac{\pa}{\pa y_i} 
+ (\sigma - 1)  \frac{\pa}{\pa \sigma} + W \Big)
b_j + j b_j = f \cdot b_j - \frac{i}{2} \Delta b_{j-1}  \ (j \geq 1)
\end{gathered}
 \end{equation*}
where all terms are smooth, $f = \half
\Delta \Phi + \frac n2$ vanishes linearly at $\Delta_b$ and $W$ is a vector field vanishing quadratically at $\Delta_b$. These equations have unique smooth solutions $b_j$ in a neighbourhood of $\Delta_b$, with $b_0 = 1$ at $\Delta_b$. 
An asymptotic sum of this formal series is therefore a solution to the equation to order $t^\infty$, i.e.\ the error term after applying $t^2(D_t + \half \Delta)$ is in $t^\infty \CI(\MMb \times [0, t_0))$ near $\Delta_b$. 

We also need our near-diagonal parametrix to be good as $x' \to 0$. To improve the error term at $x' = 0$ we expand in a Taylor series in $x'$. The error term has a Taylor series
$$
 e^{i\Psi/t(x')^2} \sum_{k=0}^\infty (x')^j e_j(t, x', \sigma, y, y'),
$$
near $\Delta_b$ 
with each $e_j = O(t^\infty)$ and smooth. We try to solve this away with a series
\begin{equation}
 e^{i\Psi/t(x')^2} \sum_{k=0}^\infty (x')^j c_j(t, x', \sigma, y, y').
\ilabel{corr-term}\end{equation}
This gives us equations of the form
\begin{equation*}\begin{gathered}
\Big( y_i  \frac{\pa}{\pa y_i} 
+ (\sigma - 1)  \frac{\pa}{\pa \sigma} + t \frac{\pa}{\pa t} + W \Big)
c_j  = t e_j + P(c_0, c_1, \dots, c_{j-1}), 
\end{gathered}
 \end{equation*}
 where $W$ is as above and $P$ is a differential operator with smooth coefficients. Since $e_j = O(t^\infty)$ there is a unique solution $c_j$ which is $O(t^\infty)$. Adding the correction term \eqref{corr-term} yields a parametrix with an error term $O(t^\infty (x')^\infty)$ locally near $\Delta_b$. Let us denote this near-diagonal parametrix, defined in a neighbourhood of $\Delta_b$ by $V_\near$. 
 
We now claim that, on the support of $d\phi$, $V_\near$ is equal to
$U_\innt$ up to $\CIdot(\MMb \times [0, t_0))$. In the interior of $\MMb$,
this follows from \cite{HW} where we showed that $V_\near$ is equal,
microlocally, to the exact propagator modulo $t^\infty C^\infty(\MMb \times
[0, t_0))$.  Our construction is such that $U_\near$ is in $\CIdot(\MMb
\times [0, t_0))$ on the support of $d\phi$ (Lemma~\ref{Ufar}) while
$U_\far$ is microsupported where the phase function is relatively
large. (Using the cutoff $\psi_f$, and the contact transformation $Q$, we
have $\overline{\tilde \nu}_1 \geq (5\epsilon/2)^2 /2$ on the microsupport
of $U_\far$, while we have $\overline{\tilde \nu}_1 \leq (2\epsilon)^2/2$
on the microsupport of $V_\near$ and on the support of $d\phi$. Here
$\overline{\tilde \nu}_1$ is the coordinate from \eqref{qcana} and $Q$
defined by \eqref{J}.) Therefore, at least away from the boundary of
$\MMb$, $V_\near$ is equal to $U_\innt$ modulo $t^\infty C^\infty(\MMb
\times [0, t_0))$ on $\supp d\phi.$

However, both $V_\near$ and $U_\innt$ are Legendre distributions associated
to the same Legendrian, and their full symbol expansion at $t=0$ is smooth
up to the boundary of $\MMb$. Since they agree everywhere in the interior
of $\MMb$ on $\supp d\phi,$ they agree up to the boundary. Hence $V_\near$
is equal to $U_\innt$ modulo $t^\infty C^\infty(\MMb \times [0, t_0))$
globally on the support of $d\phi$. Finally, both $V_\near$ and $U_\innt$
solve the Schr\"odinger equation microlocally, and we saw above that the
Taylor series of $V$ at $x' = 0$ was \emph{uniquely} determined by this
condition, it follows that $V_\near$ and $U_\innt$ are equal to all orders
in both $t$ and $x'$ microlocally near the Legendrian $L$ and on the
support of $d\phi$.

We now construct an accurate global parametrix for the propagator. Define
$$
U =  \phi V_\near +  (1 - \phi) U_{\innt} + U_\far. 
$$
We claim that this parametrix $U$ satisfies the initial condition
$$
\lim_{t \to 0} U(t) = \Id
$$ distributionally (i.e.\ the distribution limit of $U(t)$ as $t \to 0$ is equal to the
delta function on $\Delta_b$), and satisfies the equation $(D_t + \half H)
U(t) = 0$ up to an error term in $\CIdot(\MMb \times [0, t_0))$, (i.e.\
smooth and vanishing to infinite order at $t=0$ and all boundary
hypersurfaces of $\MMb$). The initial condition follows from the stationary
phase lemma applied to Legendre distributions; in particular the delta
function on the diagonal comes from $V_\near$ while $U_\innt$ and $U_\far$
contribute nothing, since the phase function is always nonzero for all
Legendre distributions comprising $U_\innt$ and $U_\far$.

To prove the claim about $U$ satisfying the equation, we write 
\begin{equation}\begin{gathered}
(D_t + \half H) U(t) = \phi (D_t + \half H) V_\near  + (1 - \phi)  (D_t + \half H) U_\innt + (D_t + \half H) U_\far  \\
+ \nabla \phi \cdot \nabla (V_\near - U_\innt) + \half \Delta \phi (V_\near - U_\innt).
\end{gathered}\ilabel{errors}\end{equation}
We have arranged that $V_\near$ is an accurate parametrix on the support of
$\phi$, so the first term is in $\CIdot(\MMb \times [0, t_0))$. Next,
Lemma~\ref{Ufar} shows that the sum of the second and third terms is in
$\CIdot(\MMb \times [0, t_0))$. Third, we have seen that $V_\near$ is equal
to $U_\innt$ up to $\CIdot(\MMb \times [0,t_0))$ on the support of $d
\phi$. It follows that the last two terms on the right hand side of
\eqref{errors} are in $\CIdot(\MMb \times [0, t_0))$. This completes the
proof that $U$ is a parametrix up to $\CIdot(\MMb \times [0, t_0))$ errors.

Finally we correct the error term. It follows from a commutator argument due to Craig, \cite{MR99f:35028b}
Th\'eor\`eme~14, that
\begin{equation}
e^{-itH/2} : \CIdot(M) \to \CIdot(M) \ \text{ for all } t.
\label{Schw-map}\end{equation}
We can correct our parametrix $U$ to the exact propagator by adding to $U$ the kernel
$$
i \int_0^t e^{-i(t-s)H/2} \big( (D_s + \half H) U(s) \big) \, ds \in t^\infty \CIdot(\MMb \times [0,t_0)).
$$
Since $\phi V_\near$, $(1 - \phi) U_\innt$, $U_\far$  and elements of $\CIdot(\MMb \times [0,t_0))$ are all Legendre distributions associated to the conic pair $(Q(L), Q(L_2^\sharp))$, the proof of the theorem is complete. 

\begin{remark} One might wonder why it was necessary to use the cutoff $\chi_1(\lambda \sqrt{t})$, instead of a $t$-independent cutoff. 
The reason is that a $t$-independent cutoff will yield a term that is
smooth on $\MMsc$ down to $t=0$. This term does not lie in the space
\eqref{propspace} so it would have to be eliminated by an a posteriori
argument. In this respect it is not so different from the term $U_{\near,
  1}$; however $U_{\near, 1}$ is localized close to the diagonal so it
automatically becomes harmless when we glue in our $V_\near$ parametrix,
which is a little more convenient. 
\end{remark}

\begin{remark} Note that $U_\far$ need not be supported away from the
  diagonal. In fact, if there is a geodesic curve on $M$ that
  self-intersects, then there will be a corresponding part of $U_\far$
  supported over the diagonal, although microlocally it will be away from
  the zero section. It is for this reason that we introduce $U_\innt$: we
  arranged that $U_\innt$ be supported close to, but not at, the diagonal,
  and this allowed us to piece together $V_\near$ and $U_\innt$ using the
  cutoff $\phi$ in \eqref{errors}.
\end{remark}

%%%%%%%%%%%%%%%%%%%%%%%%%%%%%%%%%%%%%%%%

\section{Poisson operator and scattering matrix}\ilabel{subsec:sojourn}
  
Having constructed the semiclassical resolvent as a Legendrian
distribution, we can now easily determine the structure of the
semiclassical Poisson operator and scattering matrix, since the kernels of
these operators are related in a simple way to the resolvent
kernel.\ah{changed the eigenvalue to $\lambda_0^2$ throughout this section;
needs to be done globally.}

We recall that the outgoing resolvent kernel was normalized, as a
half-density in $h$, as $(h^2\Delta +V - (\evosq + i0))^{-1}
|dh|^{1/2}$\ah{refer to where this is discussed earlier}. The Poisson
operator $P(h^{-1})$ may be defined by the restriction of $e^{-i\evosr/x' h}
|dr'|^{-1/2}$ times the resolvent kernel to the right boundary $\rb = H_1$
of $X$ (see Remark 8.4 of \cite{hasvas1}).  This
may be regarded as the principal symbol of the resolvent kernel at the Legendrian
$L_1$ corresponding to the base of the fibration $\partial_1 L \to L_1$  at $\rb = H_1$ (see Proposition~\ref{prop:fibre}).

Since the kernel of $P(\lambda)$ is a function on $M \times \partial M
\times [0, h_0)$ it is natural to regard $M \times \partial M \times [0,
h_0)$ as a scattering-fibred manifold, with the main face being $M \times
\partial M \times \{ 0 \}$ and the other boundary hypersurface, $\partial M
\times \partial M \times [0, h_0)$ fibred over $\partial M \times \partial
M$ by projection off the $h$ variable. To determine the Legendrian
structure of $P(h^{-1})$ we start with the geometry of the propagating
Legendrian $L$, defined in \eqref{Ldefn}, near the right boundary $\rb$ of
$\Tsfstar[\mf] \XXb$. Working near the right boundary, we use coordinates
$(x, \theta = x'/x, y, y', h, \lambda, \lambda', \mu, \mu', \tau)$, as
defined in \eqref{canonicaloneform}.

Let $W$ be the set $\{ \theta = 0, \mu' = 0 \} \subset
  \Tsfstar \XXb$, and consider the blowup $[\Tsfstar \XXb; W]$ of $\Tsfstar
  \XXb$ at $W$. Let $\tilde W$ denote the new boundary hypersurface created
  by this blowup, and write $\mubar' = \mu'/\theta;$ this is a smooth
  function in the interior of $\tilde W$.
\begin{lemma}
$\tilde W \cap \{ \lambda' =
  \evosr \}$ is diffeomorphic to $\Tsfstar (M \times \partial M \times [0, h_0))$
  and hence $\tilde W \cap \{ \lambda' = \evosr, \, h = 0 \}$ has a natural
  contact structure (degenerating at $x=0$), contactomorphic to
  $\Tsfstar[\mf] (M \times \partial M \times [0, h_0)).$
 \end{lemma}

\begin{proof}
The contact form at $\Tsfstar[\mf] \XXb$ is given in coordinates $(\theta,
x, h, y', y; \lambda', \lambda, \tau, \mu', \mu)$ by
\begin{equation}
-d\lambda' - \theta d\lambda -x \theta d\tau + \mu' \cdot dy' + \theta \mu
\cdot dy.  \ilabel{contact-mf}\end{equation} Performing the blowup of $W,$
 i.e.\ introducing the new coordinate $\mubar',$ and restricting to $\lambda'=\evosr,$ we find that this contact form becomes 
$$
\theta(-d\lambda - x d\tau + \mubar' \cdot dy' +  \mu
\cdot dy).
$$
Dividing by $\theta,$ i.e.\ taking the leading
part at $\tilde W\cap\{\lambda'=\evosr\},$ yields the contact form
\begin{equation}
-d\lambda - x d\tau + \mubar' \cdot dy' +  \mu
\cdot dy.
\ilabel{lp}\end{equation}

On the other hand, we may write the canonical one-form on 
 $\Tsfstar (M \times \partial M \times [0,
  h_0))$ as
$$
\tilde\lambda d\big( \frac1{xh} \big) + \tilde\tau  d\big( \frac1{h} \big)
+ \tilde\mu \frac{dy}{xh} + \tilde\mu' \frac{dy'}{xh};
$$
in the induced canonical coordinates, the contact form on this space
becomes
$$
-d\tilde\lambda - x d\tilde\tau + \tilde\mu' \cdot dy' +  \tilde\mu
\cdot dy,
$$
hence identifying
$\lambda$ with $\tilde\lambda$, $\tau$ with $\tilde\tau$, $\mu$ with
$\tilde\mu$ and $\mubar'$ with $\tilde\mu'$ exhibits the desired
contactomorphism.
\end{proof}

\begin{lemma}\label{lemma:sr}
 The propagating Legendrian $L$ intersects $\tilde W \cap \{ \lambda' =
 \evosr, \, h = 0 \}$ transversally, hence using the identification above
 we may regard the boundary of $L$ at $\tilde W$, which we denote $\SR$
 (for `sojourn relation'), as a submanifold of $\Tsfstar[\mf] (M \times
 \partial M \times [0, h_0))$. Making this identification, then $\SR$ is a
 Legendre submanifold of $\Tsfstar[\mf] (M \times \partial M \times [0,
 h_0))$ which is smooth after further blowup of $\{ x = 0, \mu = \mubar' =
 0 \}$.
\end{lemma}

\begin{proof}
Since $L$ is Legendrian in $\Tsfstar \XXb,$ the form \eqref{contact-mf}
vanishes on it.  Near $\rb$, since $L$ is contained in $\Sigma_l$, the left
characteristic variety, we have $(\lambda')^2 = \evosq - h^{ij}(x',y')\mu'_i
\mu'_j-V(\theta x, y').$ Lemma~\ref{b} shows that $L$ meets $\{ \theta = 0
\}$ only in the interior of the blowup of $W$ and does so transversely, so
we can use the blow-up variable $\mubar'_i = \mu'_i/\theta$. In terms of
this we have
$$
(\lambda')^2 = \evosq - \theta^2 h^{ij} \mubar'_i \mubar'_j -V(\theta x, y') \implies d\lambda' =
\frac{\theta}{\lambda'} h^{ij} \mubar'_i \mubar'_j d\theta + O(\theta^2) 
$$ (recall that $V(x,y) = O(x^2)$).  Thus $d\lambda'/\theta$, which by \eqref{contact-mf} is equal to \eqref{lp} on $L$, vanishes when restricted to $L \cap \{ \theta = 0 \}$.

Now we consider the smoothness of $\SR$ at the boundary $\{ x = 0 \}$. 
By Lemma~\ref{b},
 $L$ is desingularized by blowing up first $Z = \{ \mu'  = \mu = 0, x = 0,
\lambda = \lambda' \}$ and then the lift of $W$. Thus away from $Z$, the
first blowup has no effect and $L$ is desingularized
by the blowup of $W$. We have to analyze the situation further near $L \cap
Z$. Here we can take advantage of the explicit formula for $L \cap \{ x = 0
\}$ given by \eqref{eq:sp-1c}. At $x = 0$, we have\jw{Rephrased this a bit.}
\begin{equation}
\theta = \frac{|\mu'|}{|\mu|}.
\ilabel{tmm}\end{equation}
It follows that in a neighbourhood of $L \cap Z$ we have $|\mu| > |\mu'|$. Similarly, we have 
\begin{equation}
\lambda' - \lambda = \sqrt{\evosq - |\mu'|^2-V(x,y)} - \sqrt{\evosq -
  |\mu|^2-V(x,y)} = O(x^2+|\mu|^2 + |\mu'|^2) = O(x^2+|\mu|^2)\quad
\text{on } L.\ilabel{llambda}\end{equation} It follows that after $Z$ is
 blown up, we may cover a neighborhood of the intersection of $L$ and the
 front face by coordinate charts in which either $x$ or $|\mu|$ is a
 boundary defining function. Thus, in place of $x, \mu, \mu', \lambda'
 - \lambda,$ we may take as coordinates  
$$
 \frac{\mu}{x}, \frac{\mu'}{x}, \frac{\lambda' - \lambda}{x} \text{ and } x,
$$
in the region where $dx \neq 0$, and 
$$
 \frac{\mu}{|\mu|}, \frac{\mu'}{|\mu|}, \frac{\lambda' - \lambda}{|\mu|}, \frac{x}{|\mu|} \text{ and } |\mu|.
$$
in the region where $d|\mu| \neq 0$. As for the second blowup, of $\{ \theta = 0, \mu/\rho_{\tilde Z} = 0 \}$, where $\rho_{\tilde Z}$ is a boundary defining function for the face $\tilde Z$ created by the $Z$ blowup,  \eqref{tmm} implies that $\theta$ may be taken as a boundary defining function for the new face in a neighbourhood of $L$. Thus coordinates replacing those above become
\begin{equation}
\frac{\mu}{x}, \frac{\mu'}{x\theta} = \frac{\mubar'}{x}, \frac{\lambda' - \lambda}{x} \text{ and } x; \text{ and } \theta, y, y', h, \lambda, \tau 
\ilabel{coord-1r}\end{equation}
in the region where $dx \neq 0$, and 
\begin{equation}
 \frac{\mu}{|\mu|}, \frac{\mu'}{|\mu|\theta} = \frac{\mubar'}{|\mu|},
\frac{\lambda' - \lambda}{|\mu|}, \frac{x}{|\mu|} \text{ and } |\mu|;
\text{ and } \theta, y, y', h, \lambda, \tau
\ilabel{coord-2r}\end{equation} in the region where $d|\mu| \neq 0$.  It
follows from this and from Lemma~\ref{b} that $\theta$, $x$, and $2n-2$ of
the remaining coordinates from \eqref{coord-1r} (in the first region), or
$\theta$, $x/|\mu|$, $|\mu|$ and $2n-3$ of the remaining coordinates from
\eqref{coord-2r} (in the second region) furnish coordinates on $L$ on this
space, and the remaining coordinates (restricted to $L$) can be written as
smooth functions of these coordinates on $L$. Restricting to $\{ \theta = 0
\}$, then, we see that $\SR$ is desingularized by blowing up $\{ x = 0, \mu
= 0, \mubar = 0, \lambda' - \lambda = 0 \}$.

We can also observe that $(\lambda' - \lambda)/x$ or $(\lambda' -
\lambda)/|\mu|$ cannot serve as a coordinate on $L$ at $\mu = 0$, since we
see from \eqref{llambda} that this function has vanishing differential
there. This implies that $\SR$ is also desingularized by blowing up  
$$
\{ x = 0, \mu = 0, \mubar = 0 \},
$$
which completes the proof of the lemma. 
\end{proof}

\begin{remark} This lemma shows that $\SR$ forms a Legendre conic pair with
  the Legendre submanifold $G^\sharp = \{ x = 0, \lambda = \evosr, \mu = 0,
  \mubar' = 0 \}$ which is contained in the contact manifold $\Nsfstar
  (\partial M \times \partial M)$, $\partial M \times \partial M$ being the
  base of the fibration at the hypersurface at $x=0$ of the
  scattering-fibred manifold $M \times \partial M \times [0, h_0)$.
\end{remark}

To interpret the Legendrian $\SR$ geometrically, we recall the definition
of the sojourn relation from \cite{HW} (in fact, we need to generalize it
to include the case of a nonzero potential). $\SR$ is the graph of a
contact transformation $S$ from $S^* M^\circ$ to $\Tscstar _{\partial M} M$
given as follows: given a unit covector $(z, \hat \zeta) \subset S^*
M^\circ$, we let $\gamma(s)$ be the bicharacteristic (geodesic, in the case
of no potential) emanating from $(z, \hat \zeta)$. By assumption,
$\gamma(s)$ tends to the boundary, and there is a well-defined final
`direction' $y$.  The \emph{action} $A(s)$ accumulated along the
bicharacteristic is the integral of $\evosq - V$ with respect to $s$ along
$\gamma$, with initial condition $A(0) = 0$. Since $x = O(s^{-1})$ along
$\gamma$ and $V = O(x^2)$, we see that $A(s) = \evosq s + O(1)$.  Moreover,
it follows from the regularity of the boundary of $\SR$
(Lemma~\ref{lemma:sr}) that $|\mu| = O(x)$, hence $\dot r = \evosr +
O(s^{-2})$ and so $r(s) = \evosr s + O(1)$.  We let $\nu$, the \emph{sojourn
time}, be defined by $\nu = \lim_{s \to \infty} A(s) - \evosr r(s),$ which is
well defined by the above considerations. We finally define $M = \lim_{s
\to \infty} \mu(\gamma(s))/s$. Then the sojourn relation is $S(z, \hat
\zeta) = (y, \nu, M) \subset \Tscstar[\partial M] M$.\footnote{The sojourn
relation $S$ actually depends on a choice of coordinates; it is invariantly
defined on a certain affine bundle identified in \cite{HW}.} If $V \equiv
0$ then $A(s)$ is $\evosr$ times the geodesic distance along $\gamma$.

\begin{lemma}\ilabel{SRlemma} The Legendrian $\SR$ is the (twisted) graph of the sojourn relation in the interior of 
$\Tsfstar[\mf] (M \times \partial M \times [0, h_0))$. 
\end{lemma}

\begin{proof}
Consider a local parametrization of the Legendrian $L$ near $\rb$ and away from $x = 0$. The Legendrian $L_{\rb} = \{ \lambda' = \evosr, \mu' = 0 \}$ is parametrized by the phase function $\evo/\theta x h$, so we can choose our phase function to be of the form $(\evosr + x' \psi)/x' h$, where 
$\psi = \psi(y', z, v)$ and it is non-degenerate in the sense that
\begin{equation}
d_{z, v} \big( \frac{\partial \psi}{\partial v_i} \big) \text{ are linearly independent}, \ i = 1 \dots k \text{ where } v \in \RR^k.
\end{equation}
Then $L$ is given locally by 
$$
L = \{ (x', y', z, \evosr + x' \psi + (x')^2 \psi_{x'}, \psi,  d_{y'} \psi, d_z \psi) \mid d_v \psi = 0 \}
$$
in coordinates $(x', y', z, \lambda', \tau, M', \zeta)$ given by writing covectors in the form
$$
\lambda' d\big( \frac1{x'h} \big) + \tau d\big( \frac1{h} \big) + M' \frac{dy'}{h} + \zeta \frac{ dz}{h}.
$$ 
By \eqref{charsets}, we have $\abs{\zeta}_g^2+ V=\evosq,$ hence under the
flow of $h^{-1}$ times the Hamilton vector field,
\eqref{Vr1} gives $\dot{\tau} = \evosq-V.$  In other words,
$$
\frac{\evosr}{x'}+\psi = \int \evosq-V \, ds = A(s).
$$ Thus $\psi(0, y', z, v) = \lim_{s \to \infty} A(s) - \evosr/x',$ which is
the sojourn time (when $d_v \psi = 0$). Moreover, $d_{y'} \psi = M' =
\mu'/x'$ where $\mu'$ is the variable dual to $dy'/x'h$. Finally $d_z \psi
= d_z (\evosr + x' \psi)/x'$ gives minus the covector $\hat \zeta$ at $z$
which is the initial condition $(z, \hat \zeta)$ for the bicharacteristic.

The boundary of the Legendrian $L$ at $\tilde W \cap \{ \lambda' = \evosr \}$ is given in these coordinates by
$$
\SR = \{ (y', z, \psi, d_{y'} \psi, d_z \psi) \mid d_v \psi = 0 \}
$$
and it is now evident from the interpretations of $\psi$, $d_{y'} \psi$ and $d_z \psi$ that this is a non-degenerate parametrization of the sojourn relation.\footnote{It would be more correct to say that we are `identifying' this with the sojourn relation; it is not exactly the same as the sojourn relation as defined in \cite{HW} since it lies in a different bundle, with different scalings as $x' \to 0$. This can be traced to the fact that the bicharacteristics in \cite{HW} tend to infinity quadratically, while here they move to infinity linearly, reflecting the different scalings in the two operators (propagator vs. resolvent). The two bundles are related via the identification $Q$ in \eqref{J}.}
\end{proof}

\begin{proposition} The Poisson operator  is a Legendrian conic pair  associated to the Legendre submanifold $\SR$ and the submanifold $G^\sharp$; in fact, $$P(h^{-1}) \in I^{0, (n-1)/2; \, 0}(M \times M \times [0, h_0); (\SR, \Gsharp)).$$
\end{proposition}

\begin{remark} The fact that the orders of $P(h^{-1})$ at $\mf$ and at $\bf$ are equal to zero reflects that the fact that the Poisson operator is a unitary operator mapping between $M$ and the space $\partial M \times \RR_+$ with a conic (i.e.\ scattering) metric, as proved in \cite{hasvas1}, section 9. 
\end{remark}

\begin{proof} 
The kernel of the resolvent is given by a finite sum of oscillatory
integrals, each giving a Legendre distribution associated to the
propagating Legendrian $L$, plus a smooth term vanishing at rapidly at each
boundary hypersurface of $\XXb \times [0, h_0)$. Consider a single
oscillatory integral expression involving a phase function parametrizing
some piece of $L$. There are four different types of such expressions to
consider, corresponding to regions of $L$ which are (i) away from $\{ x = 0
\}$, (ii) near $\{ x = 0 \}$ but away from $\{ \mu = \mubar' = 0 \}$, (iii)
near $\{ x = 0, \mu = \mubar' = 0 \}$ and near the codimension three corner
of $L$, (iv) near $\{ x = 0, \mu = \mubar' = 0 \}$ but away from the
codimension three corner of $L$.

In region (i), the result follows directly from the proof of
Lemma~\ref{SRlemma}. The proof in the other regions follows the same
pattern; we need only check that we can choose a non-degenerate phase
function of the form $(\evosr + \theta \psi)/x\theta h$ for $(L,
L^\sharp)$\ah{check notation} in each region, such that $\psi$ is a
non-degenerate parametrization of $(\SR, G^\sharp)$. This was explicitly
noted in \eqref{subpar-1}, which covers regions (i) and (ii)).  In the case
of region (iii), we can use a parametrization $\Psi$ as in \eqref{Psi1};
the corresponding function $\psi$ above is $\evosr + s \psi_2 + x_2
\psi_3$, using notation from \eqref{Psi1}. Comparison of \eqref{Psi-nd} and
\eqref{psi-nd} shows that $(\Psi - \evosr)/x_1$ is a non-degenerate phase
function (where we need to make the transformation $x_2 \to x_1, x_3 \to
x_2, (y_1, y_2) \to y_1$, $y_2 \to \{ \}$, $v_2 \to v_1, v_3 \to v_2$ to
make the comparison) in the sense of \eqref{psi1psi2}. Since we know that
it parametrizes $\SR$ for $x_2 > 0$, it follows that this is a
non-degenerate parametrization of $(\SR, \Gsharp)$. In region (iv) the
result follows from the analogous comparison of \eqref{Psi2} and
\eqref{psi-simple}.

To determine the orders, notice that we divided by the half-density factor
$|dr'|^{1/2}$ to obtain the Poisson kernel. In terms of the boundary
defining functions $x_1$ for $\rb$ and $x_2$ for $\bfc$, this is dividing
by $|dx_1/x_1^2 x_2|^{1/2}$. The semiclassical order is decreased by $-1/4$
accounting for the change in total dimension from $N$ to $N-1$, but the
orders at $\bfc$ \emph{increase} by $1/4$ in view of the power $x_2^{-1/2}$
in $|dx_1/x_1^2 x_2|^{1/2}$. This shows that the new orders are as stated
in the proposition.
\end{proof}

We now turn to the analysis of the scattering matrix $S(h^{-1})$. This is
defined on $f \in \CI(\partial M)$ by distributionally restricting the
outgoing part of $x^{-(n-1)/2} e^{-i\evosr/xh} P(h^{-1}) f$ to $\partial
M$. In terms of kernels, and taking into account the half-density factors,
it may be constructed from the Poisson operator by microlocalizing near
the intersection of $\SR$ and $\Gsharp$, multiplying by
$e^{-i\evosr/x\lambda} |dx/x^2|^{-1/2}$ and restricting to $x=0$.

Thus the only part of the Legendrian $\SR$ of importance for the scattering
matrix is the part in a neighbourhood of $\mu = 0$, i.e.\ at the blowup of
$Z$. Thus we make a further symplectic reduction  and restrict $\SR$ to the
face $Y$ created by the blowup of $\{ x = 0, \mu = 0, \mubar' = 0 \}
\subset \Tsfstar(M \times \partial M \times [0, h_0))$; let $T$ denote this
  set.

Lemma~\ref{lemma:sr} tells us that $T$ is a Legendrian-Lagrangian
submanifold of $\Tsfstar (\partial M \times \partial M \times [0,
h_0))$. Thus, the contact form, which may be written
 \begin{equation}
- d\tau + M'' \cdot dy' +  M \cdot dy .
\ilabel{contact-xx}\end{equation}
in terms of blowup coordinates $M = \mu/x$, $M'' = \mubar'/x$, vanishes at $T$.

Let us define the `total sojourn Legendrian' inside $\Tscstar[\partial M
\times \partial M \times \{ 0 \}] (\partial M \times \partial M \times [0,
h_0))$ as the set consisting of points $(y, y', \tau, M, M'')$ such that
there a point $(z, \hat \zeta)$ in the interior of $M$ with $(y, \tau_1, M)
= S(z, \hat \zeta)$, $(y', \tau_2, M'') = S(z, -\hat \zeta)$ and $\tau =
\tau_1 + \tau_2$. We can also express $\tau$ as the limit of $A(s_1, s_2) -
\evosr (1/x_2 + 1/x_1)$ where $s_1 \to -\infty$, $s_2 \to \infty$ and
$A(s_1, s_2)$ is the action accumulated along the bicharacteristic
determined by $(z, \hat\zeta)$.  If there is no potential then $\tau$ is
given by the limit of $\evosr \big( d(z_1, z_2) - 1/x(z_1) - 1/x(z_2) \big)
$ where $z_1$ goes to infinity along the geodesic in one direction and
$z_2$ goes to infinity along the geodesic in the opposite direction; this
is $\evosr$ times the original `sojourn time' defined by Guillemin
\cite{Gu}.

\begin{lemma} The Legendrian $T$ coincides with the total sojourn Legendrian.
\end{lemma}

\begin{proof}
The vector field $-V_l'$ is tangent to $\SR$ and b-normal to $Y$. Therefore, every point of $T$ is the endpoint of an integral curve of $-V_l'$ lying inside $\SR$. An arbitrary point of $T$ is therefore obtained from an interior point $(z_0, \hat \zeta_0; y', \tau, \mu')$ of $\SR$ by flowing 
along a $V_l'$ integral curve. This does not change the values of $y'$ or $M'$, while $(z, \hat \zeta)$ moves along the bicharacteristic with initial condition $(z_0, \hat \zeta_0)$. Thus when the bicharacteristic arrives at $Y$ the $y$ coordinate is the asymptotic direction of this bicharacteristic, while $M = \mu/x$ is the asymptotic `angular coordinate'. To work out an interpretation of the $\tau$ variable, notice that when we use coordinates on $\tilde W$ given by
$$
 \overline{\tau} d\big( \frac1{h}\big) + \frac{\zeta \cdot dz}{h} + M' \frac{dy'}{h} 
$$
then $\overline{\tau}$ has the interpretation of the sojourn time starting from $(z_0, \hat \zeta_0)$ (see the proof of Lemma~\ref{SRlemma}). 
Near $Y$ we change to variables given by 
$$
\lambda d\big( \frac1{xh} \big) + \tau d\big( \frac1{h} \big)  +  \mu \cdot \frac{dy}{xh} +  \mubar' \cdot  \frac{dy'}{xh}
= \lambda d\big( \frac1{xh} \big) + \tau d\big( \frac1{h} \big)  +  M \cdot \frac{dy}{h} +  M' \cdot  \frac{dy'}{h}.
$$
Comparing the two sets of coordinates gives $\tau = \overline\tau -
\lambda/x$. Since $\lambda = \evosr$ at $T$, this gives $\tau = \lim_{x \to
  0} \overline\tau - \evosr/x$ on $T$. Since $\overline\tau$ is the sojourn
time starting from $(z, \hat\zeta)$, i.e.\  the limit of $A - \evosr/x'$,
this shows that at $Y$, $\tau = \lim_{x, x' \to 0} (A - \evosr/x' -
\evosr/x)$ is the total sojourn time along the bicharacteristic determined
by $(y', M')$, or equivalently by $(y, M)$. This completes the proof that
$T$ is the total sojourn relation.
\end{proof}

\begin{proposition} The set $T$ is a Legendrian-Lagrangian submanifold of
  $\Tsfstar (\partial M \times \partial M \times [0, h_0)$, and the
  scattering matrix $S(h)$ is a Legendrian-Lagrangian distribution on
  $\partial M \times \partial M \times [0, h_0)$ associated to  $T$; indeed
  $S(h^{-1}) \in I^{-1/4, -1/4}(\partial M \times \partial M \times [0, h_0), T; \scOh)$.  
\end{proposition}

\begin{proof} 
This follows directly from Proposition~\ref{distlim}. 
\end{proof}

\end{document}